\documentclass[11pt,a4paper]{article}
\usepackage[all]{xy}
\usepackage{amsmath,amsfonts,amssymb,amscd,pb-diagram,pb-xy}
\usepackage{graphicx}
\usepackage{epsfig}
\usepackage{caption}
\usepackage{enumitem}

\newcounter{ExCounter}
\stepcounter{ExCounter}

\newcounter{RemCounter}
\stepcounter{RemCounter}

\newcounter{ResCounter}
\stepcounter{ResCounter}

\newenvironment{Remark}{\par\medskip\noindent{\sc Remark \arabic{RemCounter}\stepcounter{RemCounter}. }}{\hfill
$\Diamond$\par\medskip}
\newenvironment{Example}{\par\medskip\noindent{\sc Example. %\arabic{ExCounter}\stepcounter{ExCounter}.
}}{\hfill $\triangle$\par\medskip}
\newenvironment{Result}{\par\medskip\noindent{\sc Result \arabic{ResCounter}\stepcounter{ResCounter}. }}{\medskip}
\newenvironment{Proof}{\par\noindent{\sc Proof:}}{\hskip1.8em $\Box$\par\smallskip}

\newtheorem{Theorem}{Theorem}[section]

\newtheorem{Lemma}{Lemma}[section]
\newtheorem{Corollary}{Corollary}[section]

\newcommand{\Reali}{\ensuremath{\mathbb{R}}}
\newcommand{\Euclideo}{\ensuremath{\mathbb{E}}}
\newcommand{\Naturali}{\ensuremath{\mathbb{N}}}

\def\={\, = \,}

\def\veps{\varepsilon}
\def\vth{\vartheta}
\def\vph{\varphi}

\def\H{\mathcal H\/}

\def\M{\mathcal M\/}

\def\S{\mathcal S\/}
\def\Z{\mathcal Z\/}

\def\SO{\mathcal S_1\/}
\def\ST{\mathcal S_2\/}
\def\SOT{\mathcal S_{12}\/}

\def\J{J\/(\M)}
\def\JSO{J\/(\SO)}
\def\JST{J\/(\ST)}
\def\JSOT{J\/(\SOT)}

\def\velort{{\bf v}^{\perp}}
\def\uelortO{{\bf u}^{\perp}_1}
\def\uelortT{{\bf u}^{\perp}_2}

\begin{document}

\title{An algorithmic approach to the multiple impact of a disk in a corner.}

\author{
Claudia Fassino\footnote{Department of Mathematics, University of
Genova - Via Dodecaneso 35, 16146 GE\-NOVA - email:
fassino@dima.unige.it}, Stefano Pasquero\footnote{Department of
Mathematical, Physical and Computer Sciences, University of Parma.
Parco Area delle Scienze 53/a (Campus) 43124 PARMA. email:
stefano.pasquero@unipr.it} }
\date{\today}
\maketitle

\begin{abstract}
\noindent We present the algorithmic procedure determining the
impulsive behavior of a rigid disk having a single or possibly
multiple frictionless impact with two walls forming a corner. The
algorithmic procedure represents an application of the general
theory of multiple impacts as presented in
\cite{Pasquero2016Multiple} for the ideal case.

In the first part, two theoretical algorithms are presented for
the cases of ideal impact and Newtonian frictionless impact with
global dissipation index. The termination analysis of the
algorithms differentiates the two cases: in the ideal case, we
show that the algorithm always terminates and the disk exits from
the corner after a finite number of steps independently of the
initial impact velocity of the disk and the angle formed by the
walls; in the non--ideal case, although is not proved that the
disk exits from the corner in a finite number of steps, we show
that its velocity decreases to zero and the termination of the
algorithm can be fixed through an ``almost at rest'' condition.

In the second part, we present a numerical version of both the
theoretical algorithms that is more robust than the theoretical
ones with respect to noisy initial data and floating point
arithmetic computation. Moreover, we list and analyze the outputs
of the numerical algorithm in several cases.

\vskip 0.5truecm



\end{abstract}

\clearpage

\section*{Introduction}

The study of the behavior of a rigid or multibody system subject
to multiple contact and/or impact is a very actual argument of
investigation, finding application in several branches of
Classical Mechanics, from the analysis of the motion of billiard
balls to that of rocking blocks or granular materials.

The argument can be dealt following several different approaches,
ranging from completely theoretical to specifically numerical (see
e.g. \cite{Brogliato} for a wide but not complete bibliography).

In a recent paper \cite{Pasquero2016Multiple}, a geometric
approach framed in the context of jet--bundle theory was used to
analyze the behavior of a general mechanical system with a finite
number of degrees of freedom subject to multiple unilateral ideal
constraints. The analysis, simply based on very general arguments
of preservation of the kinetic energy, led to the construction of
a theoretical algorithm that determines the right velocity of the
system for several significant mechanical systems once the left
velocity and the geometric properties of the system are known.

In this paper we present the application of this algorithm to the
paradigmatic case of the planar system formed by a rigid disk
simultaneously impacting with both sides of a corner in two
possible situations: the ideal case, with frictionless contacts
and conservation of the kinetic energy of the disk; the so called
Newtonian frictionless impact with global dissipation index, a
non--ideal case with frictionless contact and no requirement of
conservation of kinetic energy.

In the ideal case, the algorithm is directly built on the
theoretical results of \cite{Pasquero2016Multiple}, and the main
result, apart from the analysis of the physical meaning of the
output, pertains to the termination analysis. This is an important
aspect of the approach that was not discussed in
\cite{Pasquero2016Multiple}: although for several meaningful
systems the theoretical algorithm evidently terminates, it is
clear that the requirement of conservation of kinetic energy
suggests the possibility of an infinite number of iterations of
the algorithm, reflecting a possible infinite number of rebounds
of the disk. We prove that the algorithm always terminates and
determines a velocity such that the disk exits from the corner.
However, notwithstanding the manifest simplicity of the mechanical
system, the analysis lights up two important aspects of multiple
ideal impacts: the first is that the geometry of the system can be
easily arranged in order to obtain as many iterations of the
algorithm as one can decide; the second is that the case of
effective double impact with both the sides of the corner,
although leading to a non trivial non--linear rule for the
determination of the right velocity, is however such that that the
double impact can happen only once.

The algorithm for the non--ideal case consists in a generalization
of the ideal one supposing that the (frictionless) walls of the
corner are partially or totally inelastic. The non--ideality is
introduced by a Newtonian restitution coefficient $\veps$ relating
the orthogonal components of the velocity of the disk with respect
to the walls before and after each step of the algorithm. In this
case we do not prove that the algorithm always determines a
velocity such that the disk exits from the corner, but we show
that the norm of the velocity decreases to $0$ for increasing
numbers of steps. This gives a second termination criterion for
the algorithm, with a meaningful physical interpretation. However,
even in the non--ideal case, we prove that double impacts of the
disk with both the sides of the corner can happen only once.

It is however well known that the data derived from real-world
measurements can be perturbed by errors, so that theoretical
algorithms that process such data can produce unreliable results.
Moreover, when an algorithm is implemented, because of the
floating point arithmetic, the computed values can be
perturbed by algorithmic errors. %The data uncertainty and the
%algorithmic errors can cause disastrous effects on the result.
 In case of pedestrian implementation of the theoretical
algorithms, small perturbations of the data can cause structural
changes in response, consequent different choices in the iterative
method, and than invalidate the final results.

For this reason we present a numerical algorithm, based on the
theoretical ones, that, differently to these ones, is robust with
respect to the errors introduced by the measurements and by the
use of the floating point arithmetic. This algorithm is obtained
introducing suitable thresholds changing the tests for the choice
of the iterative step to do.

\smallskip

The paper is then divided into two main parts: in the first, after
some short preliminaries, we introduce the theoretical algorithm
for the ideal impact in three different but equivalent forms and
we show the corresponding results. In particular we prove that the
disk exits from the corner after a finite number of steps
independently of the initial impact velocity of the disk and the
angle formed by the walls. Then we introduce the theoretical
algorithm for the non--ideal impact in two different but
equivalent forms and we show the corresponding results. In
particular, although we does not prove that disk exits from the
corner after a finite number of steps, we prove that its velocity
decreases to zero and the termination of the algorithm can be
fixed through an ``almost at rest'' condition.

In the second part we introduce the single numerical algorithm
that groups both the ideal and the non--ideal theoretical cases.
The numerical algorithm differs from the theoretical ones in the
criteria about the velocity that select the behavior of the disk
after a rebound. In the theoretical versions such criteria are
based on the on the nullity of some suitable components of the
velocity. It is however well known that an exact comparison with
zero makes the algorithm unstable, and so, in order to obtain a
more robust algorithm, we introduce two thresholds: one
determining when a single component of the velocity of the disk is
almost zero, one determining when the norm of the whole velocity
is almost zero. We compare the theoretical and numerical versions
of the algorithms showing that, starting from the same input, they
compute the same output in the same number of steps, or they
compute slightly different outputs, even if one of the versions
performs more steps. Finally, we illustrate the behavior of the
numerical algorithm by listing the outputs computed processing
several different meaningful inputs.

\smallskip

Since the main aim of the paper is focused on the analysis of the
algorithms giving the velocity of the disk after the impact, we
relegate to the appendix a brief but exhaustive sketch of the
geometric method determining the constitutive characterizations of
the multiple constraint that are the bases for the construction of
the algorithms. The Reader interested in a wider description of
the geometry and the impulsive aspects of general systems subject
to ideal multiple constraints can refer to
{\cite{Pasquero2016Multiple} and the references therein. In the
bulk of the paper, we will limit the mathematical aspects to the
bare necessary to describe the algorithms for the case of the disk
in the corner.

\smallskip

The list of possible references about multiple impacts is very
huge, and a bibliography claiming to be exhaustive on the argument
should be excessively long compared to the length of the paper.
Moreover, only few works would be reasonably pertinent to the
specific algorithm presented in the paper. Therefore, the list of
references has been based on the minimality criterion of making
the paper self--consistent. However, for large but not recent or
exhaustive lists of general references, see for example
\cite{Brogliato,Brogliato2000impacts,Pfeiffer2000multibody,Liu2008frictionless,Johnson,Stronge}.

\vskip 1truecm

\centerline{\bf {\LARGE PART 1: THEORETICAL ASPECTS}}

\section{Preliminaries}

A rigid disk of unitary mass and radius $R$ moves in the part of a
horizontal plane delimited by two walls $\SO, \ST$ forming an
angle $2\alpha \in (0, \pi)$ (see Fig. 1). We can describe the
geometry of the system by introducing local coordinates
$(x,y,\vth)$ where $x,y$ are the coordinates of the center of the
disk and $\vth$ is the orientation of the disk. Choosing $k = \tan
\alpha$, then $k>0$ and the walls can be described by the
cartesian relations $\SO: kx-y=0, \, \ST: kx+y=0$. Since we adopt
the so called {\it event--driven} approach,  we assume that
$(x,y)=(-\frac{R}{\sin \alpha},0)$ so that the disk is in contact
with both the walls. We assume the contacts as frictionless.

\begin{figure}[h]\label{Fig1}
\centering
%\hskip-3truecm
\includegraphics[width=0.6\textwidth, angle=90]{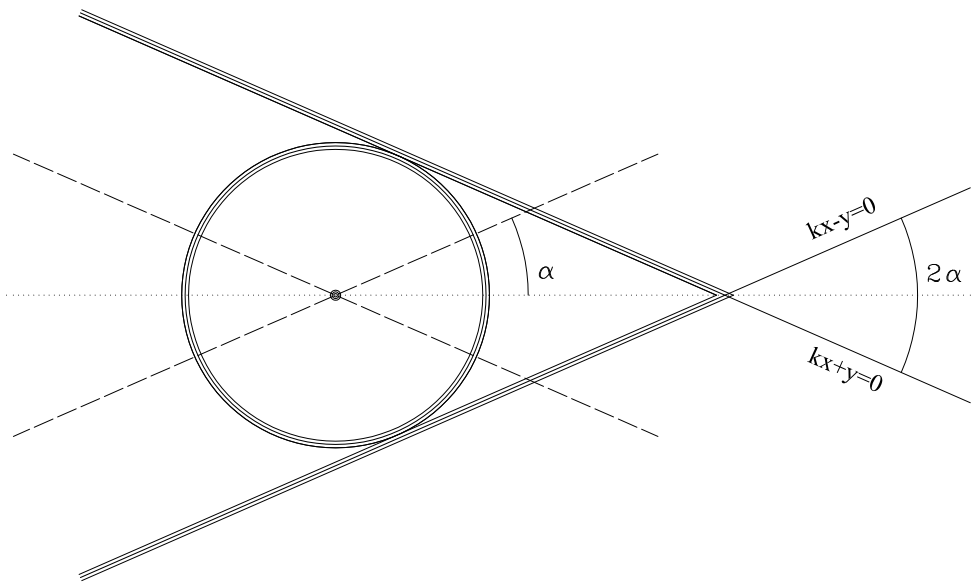}
\caption{Disk in contact with both sides of the corner}
\end{figure}

The kinetic state of the disk is known once the linear velocity of
the (center of the) disk is assigned through a pair ${\bf v} =
(\dot{x}, \dot{y})$, and the spin is assigned by $\dot{\vth}$. Of
course, the linear velocity can be assigned in the alternative
form $(v \cos \vph, v \sin \vph)$ with ($v>0$ and) $\vph \in
(-\pi, \pi]$. Without entering temporarily into mathematical
details, it is clear that the angle $\vph$ determines the nature
of the impact between disk and walls: if $\vph \in (-\alpha,
\alpha)$, the disk is subject to a multiple impact, if $\vph \in
[\alpha, \pi - \alpha)$ the disk is subject to an impact with
$\ST$ and not with $\SO$, if $\vph \in (- \pi +\alpha, - \alpha]$
the disk is subject to an impact with $\SO$ and not with $\ST$.
Otherwise, the velocity ${\bf v}$ is a so--called {\it exit
velocity}, the disk does not impact with the walls and its time
evolution separates it from one or both the walls.

We can divide the space of the linear velocity of the system in
four different zones $\Z_0, \Z_1, \Z_2, \Z_{12}$ with ${\bf v} \in
\Z_i$ if, with a slight abuse of notation, the linear velocity
determines an impact with $\S_i$.

\begin{Remark} If we restrict our attention to a physically
meaningful situation and taking into account Fig. \ref{Fig1}, the
initial velocity ${\bf v}_0=(\dot{x}_0, \dot{y}_0)$ is such that
${\bf v}_0 \in \Z_{12}$ or ${\bf v}_0$ is parallel to the walls,
that is $\vph_0\in[-\alpha, \alpha]$.
\end{Remark}

\medskip

We now resume the main consequences of the constitutive
characterization of multiple contact/impact presented in
\cite{Pasquero2016Multiple} and sketched in the Appendix for the
case of the disk in the corner:

\begin{itemize}

\item[0)] the constitutive characterization assigns a reactive
impulse ${\bf I}$ to any kinetic state of the disk supposed in
contact with the walls;

\item[1)] independently of the kinetic state of the disk, the
reactive impulse ${\bf I}$ does not involve the spin $\dot{\vth}$
of the disk, that remains unchanged in the impact. This is
coherent with the assumption that the contacts of the disk with
$\SO$ and $\ST$ are frictionless. Then we can focus our attention
only on the linear part ${\bf v}$ of the velocity of the disk;

\item[2)] for every linear velocity ${\bf v}$ we can determine the
orthogonal components $\velort_1({\bf v}), \velort_2({\bf v})$ of
${\bf v}$ with respect to $\SO, \ST$ respectively. The orthogonal
components indicate if ${\bf v}$ is an exit velocity or if ${\bf
v}$ gives rise to an impact with $\SO$ and/or $\ST$. In
particular, if ${\bf v} = (\dot{x}, \dot{y})$, we have that:

\begin{eqnarray}\label{CC_Calcolata}
\begin{array}{lcl}
{\bf v} \in \Z_0 & \Leftrightarrow  & \left\{
\begin{array}{l}
k\dot{x} + \dot{y}\le 0 \\
k\dot{x} - \dot{y}\le 0
\end{array}\right. \\ \\
{\bf v} \in \Z_1  & \Leftrightarrow &  \left\{
\begin{array}{l}
k\dot{x} + \dot{y}\le 0 \\
k\dot{x} - \dot{y} > 0
\end{array}\right. \\ \\
{\bf v} \in \Z_2  & \Leftrightarrow & \left\{
\begin{array}{l}
k\dot{x} + \dot{y} > 0 \\
k\dot{x} - \dot{y}\le 0
\end{array}\right. \\ \\
{\bf v} \in \Z_{12}  & \Leftrightarrow & \left\{
\begin{array}{l}
k\dot{x} + \dot{y} > 0 \\
k\dot{x} - \dot{y} > 0
\end{array}\right.
\end{array}
\end{eqnarray}

Note that the symmetry of the mechanical problem is reflected in
the symmetry of the zones $\Z_i$ with respect to the
$(\dot{x},\dot{y})$ components of the velocity. In fact we have
\begin{eqnarray*}
\begin{array}{lcl}
(\dot{x}_{n}, \dot{y}_{n}) \in \Z_0  \, \, & \Rightarrow &  \, \,
(\dot{x}_{n}, - \, \dot{y}_{n}) \in \Z_0 \\
(\dot{x}_{n}, \dot{y}_{n}) \in \Z_1  \, \, & \Rightarrow & \, \,
(\dot{x}_{n}, - \, \dot{y}_{n}) \in \Z_2 \\
(\dot{x}_{n}, \dot{y}_{n}) \in \Z_2  \, \, & \Rightarrow & \, \,
(\dot{x}_{n}, - \, \dot{y}_{n}) \in \Z_1 \\
(\dot{x}_{n}, \dot{y}_{n}) \in \Z_{12}  \, \, & \Rightarrow & \,
\, (\dot{x}_{n}, - \, \dot{y}_{n}) \in \Z_{12} \, .
\end{array}
\end{eqnarray*}

\item[3)] the constitutive characterization determines a rule
assigning a ``new'' velocity of the system once an ``old''
velocity is known. The rule is ${\bf v}_{new} = {\bf v}_{old} +
{\bf I}({\bf v}_{old})$, where ${\bf I}$ represents the reactive
impulse generated by the impact with the walls, and the function
${\bf I}({\bf v}_{old})$ depends on the ideal or non--ideal nature
of the impact.

Of course ${\bf I}({\bf v}_{old})\=0$ if ${\bf v} \in \Z_0$: in
this case, the reactive impulse acting on the disk is null and the
kinetic energy is obviously preserved, as well as the euclidean
norm $\| {\bf v}\|_{_{2}}=\sqrt{\dot{x}^2 + \dot{y}^2}$.

The rule ${\bf I}={\bf I}({\bf v}_{old})$ consists in a complete
or partial ``reflection'' of the orthogonal component
$\velort_i({\bf v}_{old})$ if ${\bf v}_{old} \in \Z_i, i=1,2$. In
the ideal case, the reactive impulse acting on the disk has the
form ${\bf I}({\bf v}_{old}) \= -2\velort_i({\bf v}_{old})$, so
that the reflection of the component $\velort_i({\bf v}_{old})$ is
complete. The kinetic energy and $\| {\bf v}\|_{_{2}}$ are once
again preserved. In the non--ideal case, the reactive impulse
acting on the disk has the form ${\bf I}({\bf v}_{old}) \=
-(1+\veps)\,\velort_i({\bf v}_{old})$ with $0\le \veps <1$, so
that the reflection of the component $\velort_i({\bf v}_{old})$ is
only partial. The kinetic energy of the system is not preserved
and  $\| {\bf v}_{new}\|_{_{2}} < \| {\bf v}_{old}\|_{_{2}}$ .

The rule ${\bf I}={\bf I}({\bf v}_{old})$ consists in a strongly
non--linear relation in the components of ${\bf v}_{old}$ in case
of multiple impact ${\bf v}_{old} \in \Z_{12}$  (see
(\ref{regola_x_y4}) and (\ref{regola_x_y_bis4}) below): in this
case, the reactive impulse has the form ${\bf I}({\bf v}_{old}) \=
\lambda (\velort_1({\bf v}_{old})+\velort_2({\bf v}_{old}))$ where
$\lambda$ is a suitable coefficient fit to obtain the conservation
of the kinetic energy in the ideal multiple impact (see
\cite{Pasquero2016Multiple}) or involving the Newtonian
restitution coefficient $\veps$ in the non--ideal one;

\item[4)] the iterative application of the rule determines an
algorithm fit to determine if and how the disk exits from the
corner, loosing contact with one or both the walls. The output of
the algorithms must become constant if the ``new'' velocity
belongs to the $\Z_0$ zone. Then, both in the ideal and the
non--ideal cases, the termination analysis of the algorithm is
necessarily related to the condition ${\bf v} \in \Z_0$ (or, of
course, on the number of steps). In the non--ideal case, the
termination analysis can be also based on the condition $\|{\bf
v}\|_{_{2}} \le constant$ upon proof that $\lim_{n\to +\infty}
\|{\bf v}\|_{_{2}} = 0$.

\end{itemize}

\section{Theoretical algorithm for the ideal case (TA${}_{id}$)}

In this section we present the iterative rule assigning the
``new'' velocity ${\bf v}_{n+1}$ of the disk as function of the
``old'' velocity ${\bf v}_{n}$ for the ideal case in three
different forms. Each one of the forms will be used to obtain
theoretical results about the algorithm.

\subsection{First expression of TA${}_{id}$: use of $(\dot{x},
\dot{y})$}

Given an initial velocity ${\bf v}_0 = (\dot{x}_0, \dot{y}_0)$,
the iterative rule determined by the constitutive characterization
of \cite{Pasquero2016Multiple} is such that:
\begin{subequations}
\begin{eqnarray}\label{regola_x_y1}
\hskip -3.2truecm \begin{array}{lllll} \mbox{If } (\dot{x}_{n},
\dot{y}_{n}) \in \Z_0 \, , & \mbox{that is if } & \left\{
\begin{array}{l}
k\dot{x}_n + \dot{y}_n\le 0 \\
k\dot{x}_n - \dot{y}_n\le 0
\end{array}\right. \, , & \mbox{then } &
\left\{
\begin{array}{l}
\dot{x}_{n+1} \= \dot{x}_{n} \\
\dot{y}_{n+1} \= \dot{y}_{n}
\end{array}\right.
\end{array}
\end{eqnarray}
\begin{eqnarray}\label{regola_x_y2}
\hskip -2truecm \begin{array}{lllll} \mbox{If } (\dot{x}_{n},
\dot{y}_{n}) \in \Z_1 \, , & \mbox{that is if } & \left\{
\begin{array}{l}
k\dot{x}_n + \dot{y}_n\le 0 \\
k\dot{x}_n - \dot{y}_n > 0
\end{array}\right. \, , & \mbox{then } &
\left\{
\begin{array}{l}
\dot{x}_{n+1} \= \dfrac{1-k^2}{1+k^2} \, \dot{x}_{n} +
\dfrac{2k}{1+k^2} \, \dot{y}_{n} \\ \\
\dot{y}_{n+1} \= \dfrac{2k}{1+k^2} \, \dot{x}_{n} -
\dfrac{1-k^2}{1+k^2} \, \dot{y}_{n}
\end{array}\right.
\end{array}
\end{eqnarray}
\begin{eqnarray}\label{regola_x_y3}
\hskip -2truecm \begin{array}{lllll} \mbox{If } (\dot{x}_{n},
\dot{y}_{n}) \in \Z_2 \, , & \mbox{that is if } & \left\{
\begin{array}{l}
k\dot{x}_n + \dot{y}_n > 0 \\
k\dot{x}_n - \dot{y}_n \le 0
\end{array}\right. \, , & \mbox{then } &
\left\{
\begin{array}{l}
\dot{x}_{n+1} \= \dfrac{1-k^2}{1+k^2} \, \dot{x}_{n} -
\dfrac{2k}{1+k^2} \, \dot{y}_{n} \\ \\
\dot{y}_{n+1} \= - \dfrac{2k}{1+k^2} \, \dot{x}_{n} -
\dfrac{1-k^2}{1+k^2} \, \dot{y}_{n}
\end{array}\right.
\end{array}
\end{eqnarray}
\begin{eqnarray}\label{regola_x_y4}
\hskip -2truecm \begin{array}{lllll} \mbox{If } (\dot{x}_{n},
\dot{y}_{n}) \in \Z_{12}  \, , & \mbox{that is if } & \left\{
\begin{array}{l}
k\dot{x}_n + \dot{y}_n > 0 \\
k\dot{x}_n - \dot{y}_n > 0
\end{array}\right. \, , & \mbox{then } &
\left\{
\begin{array}{l}
\dot{x}_{n+1} \= \dfrac{-k^4\dot{x}_n^2 +
(1-2k^2)\dot{y}_n^2}{k^4\dot{x}_n^2 + \dot{y}_n^2} \, \dot{x}_{n}
 \\ \\
\dot{y}_{n+1} \= \dfrac{k^2(k^2-2)\dot{x}_n^2 -
\dot{y}_n^2}{k^4\dot{x}_n^2 + \dot{y}_n^2} \, \dot{y}_{n}
\end{array}\right.
\end{array}
\end{eqnarray}
\end{subequations}

\begin{Remark} A straightforward calculation shows that independently of
the condition ${\bf v}_n \in \Z_{i}$ with $i=0,1,2,12$, we have
\begin{eqnarray}\label{rapporto_norme}
\dfrac{\left(\| {\bf v}_{n+1}\|_{{}_2}\right)^2}{\left(\| {\bf
v}_{n}\|_{{}_2}\right)^2} \= \dfrac{\dot{x}_{n+1}^2 \, + \,
\dot{y}_{n+1}^2}{\dot{x}_{n}^2 \, + \, \dot{y}_{n}^2} \=1 \, .
\end{eqnarray}
This is an easily predictable but not trivial consequence of the
preservation of the kinetic energy required in
\cite{Pasquero2016Multiple}. In fact, since the kinetic energy is
not an absolute quantity but it depends on the choice of a frame
of reference, the validity of (\ref{rapporto_norme}) follows from
the nature itself of the contact/impact, that does not affect the
angular coordinate of the disk, and the nature itself of the
constraint and its set of rest frames. Moreover, $\| {\bf
v}\|_{{}_2}$ is not the norm of the velocity vector of the disk
but only the Euclidean norm of the pair $(\dot{x},\dot{y})$ viewed
as an element of $\Reali^2$ (see Appendix for details).
\end{Remark}

\begin{Remark} The iterative rule (\ref{regola_x_y1}--\ref{regola_x_y4}) respects the symmetry
of the mechanical problem with respect to the $(\dot{x},\dot{y})$
components of the velocity. In fact an easy calculation shows
that, if $(\dot{x}_{n}, \dot{y}_{n}) \notin \Z_0$, then
\begin{eqnarray}\label{rispetto_simmetria}
\left\{
\begin{array}{lcr}
\dot{x}_{n+1}(\dot{x}_{n}, - \, \dot{y}_{n}) & = &
\dot{x}_{n+1}(\dot{x}_{n}, \dot{y}_{n}) \\
\dot{y}_{n+1}(\dot{x}_{n}, - \, \dot{y}_{n}) & = & -\,
\dot{y}_{n+1}(\dot{x}_{n}, \dot{y}_{n})
\end{array}\right.
\end{eqnarray}
\end{Remark}

\subsection{Second expression of TA${}_{id}$: use of $(\cos \vph, \sin
\vph)$}

The same algorithm can be expressed using ${\bf v}_0$ in the form
$(v_0 \cos \vph_0, v_0 \sin \vph_0)$ with $v_0>0, \, \vph \in
(\pi, \pi]$. Since $ k = \tan \alpha$ so that
\begin{eqnarray*}
\left\{
\begin{array}{lcl}
 \dfrac{1-k^2}{1+k^2} &=& \cos 2\alpha \\ \\
\dfrac{2k}{1+k^2} &=&  \sin 2\alpha \, ,
\end{array} \right.
\end{eqnarray*}
we immediately obtain that the matrices of the linear
transformations given by (\ref{regola_x_y2},\ref{regola_x_y3}) are
orthogonal but not special orthogonal. Moreover, thanks to
(\ref{rapporto_norme}), $v_0$ is factorized in every term. The
iterative rule becomes
\begin{subequations}
\begin{eqnarray}\label{regola_cosfi_sinfi1}
\hskip -3.5truecm \begin{array}{lllll} \mbox{If } \vph_{n} \in
\Z_0 \, , & \mbox{that is if } & \left\{
\begin{array}{l}
\cos \vph_n < 0 \\
| \tan \vph_n | \le \tan \alpha
\end{array}\right. \, ,  &
\mbox{then } & \vph_{n+1} \= \vph_{n}
\end{array}
\end{eqnarray}
\begin{eqnarray}\label{regola_cosfi_sinfi2}
\hskip -1.6truecm \begin{array}{lllll} \mbox{If } \vph_{n} \in
\Z_1 \, , & \mbox{that is if } & \left\{
\begin{array}{l}
\sin \vph_n < 0 \\
- \cot \alpha < \cot \vph_n  \le \cot \alpha
\end{array}\right. \, , &
\mbox{then }  & \vph_{n+1} \=  - \, \vph_{n} \, + \, 2\alpha
\end{array}
\end{eqnarray}
\begin{eqnarray}\label{regola_cosfi_sinfi3}
\hskip -1.6truecm \begin{array}{lllll} \mbox{If } \vph_{n} \in
\Z_2 \, , & \mbox{that is if } & \left\{
\begin{array}{l}
\sin \vph_n > 0 \\
- \cot \alpha < \cot \vph_n  \le \cot \alpha
\end{array}\right. \, ,  &
\mbox{then } & \vph_{n+1} \= - \, \vph_{n} \,  - \, 2\alpha
\end{array}
\end{eqnarray}
\begin{eqnarray}\label{regola_cosfi_sinfi4}
\begin{array}{l}
\hskip -1.8truecm \begin{array}{lllll} \mbox{If } \vph_{n} \in
\Z_{12}  , & \mbox{that is if } & \left\{
\begin{array}{l}
\cos \vph_n > 0 \\
| \tan \vph_n | < \tan \alpha
\end{array}\right. \, ,  &
\mbox{then } &
\end{array}
\\ \\
\hskip -1truecm \begin{array}{lllll} & \left\{\begin{array}{l}
\cos \vph_{n+1} \= - \, \left( \dfrac{\tan^4 \alpha \cos^2 \vph_n
- \sin^2 \vph_n}{\tan^4 \alpha \cos^2 \vph_n + \sin^2 \vph_n} \, +
\, 2 \, \dfrac{\tan^2 \alpha \sin^2 \vph_n}{\tan^4 \alpha\cos^2
\vph_n + \sin^2 \vph_n} \right) \, \cos \vph_{n}
 \\ \\
\sin \vph_{n+1} \= \phantom{- \, } \left( \dfrac{\tan^4
\alpha\cos^2 \vph_n - \sin^2 \vph_n}{\tan^4 \alpha\cos^2 \vph_n +
\sin^2 \vph_n} \, - \, 2 \, \dfrac{\tan^2 \alpha \cos^2
\vph_n}{\tan^4 \alpha\cos^2 \vph_n + \sin^2 \vph_n}  \right) \,
\sin \vph_{n}
\end{array}\right. &&&
\end{array}
\end{array}
\end{eqnarray}
\end{subequations}

\subsection{Third expression of TA${}_{id}$: use of $(\dot{\xi},
\dot{\eta})$}

A third version of the algorithm can be obtained by using a
standard change of coordinates $(\xi, \eta) \= (k x + y,k x - y)$
that identifies the velocity using its projections in the
directions of the walls. Then we have
\begin{eqnarray}\label{cambiamento_coordinate}
\left\{
\begin{array}{l}
\dot{\xi} \= k \dot{x} + \dot{y} \\
\dot{\eta} \= k \dot{x} - \dot{y}
\end{array} \right.
 \quad \Leftrightarrow \quad
\left\{
\begin{array}{l}
\dot{x} \= \dfrac{\dot{\xi} + \dot{\eta}}{2k} \\ \\
\dot{y} \= \dfrac{\dot{\xi} - \dot{\eta}}{2}
\end{array} \right.
 \, .
\end{eqnarray}
In this case we have:
\begin{subequations}
\begin{eqnarray}\label{regola_xi_eta1}
\hskip -6.2truecm \begin{array}{lllll} \mbox{If } (\dot{\xi}_{n},
\dot{\eta}_{n}) \in \Z_0 \, , & \mbox{that is if } & \left\{
\begin{array}{l}
\dot{\xi}_n\le 0 \\
\dot{\eta}_n \le 0
\end{array}\right. \, ,  & \mbox{then } &
\left\{
\begin{array}{l}
\dot{\xi}_{n+1} \= \dot{\xi}_{n} \\
\dot{\eta}_{n+1} \= \dot{\eta}_{n}
\end{array}\right.
\end{array}
\end{eqnarray}
\begin{eqnarray}\label{regola_xi_eta2}
\hskip -3.8truecm \begin{array}{lllll} \mbox{If } (\dot{\xi}_{n},
\dot{\eta}_{n}) \in \Z_1 \, , & \mbox{that is if } & \left\{
\begin{array}{l}
\dot{\xi}_n\le 0 \\
\dot{\eta}_n > 0
\end{array}\right. \, , & \mbox{then } &
\left\{
\begin{array}{l}
\dot{\xi}_{n+1} \= \dot{\xi}_{n} \, + \,
2 \, \dfrac{1-k^2}{1+k^2} \, \dot{\eta}_{n} \\ \\
\dot{\eta}_{n+1} \= - \, \dot{\eta}_{n}
\end{array}\right.
\end{array}
\end{eqnarray}
\begin{eqnarray}\label{regola_xi_eta3}
\hskip -3.8truecm \begin{array}{lllll} \mbox{If } (\dot{\xi}_{n},
\dot{\eta}_{n}) \in \Z_2\, ,  & \mbox{that is if } & \left\{
\begin{array}{l}
\dot{\xi}_n > 0 \\
\dot{\eta}_n \le 0
\end{array}\right.\, ,  & \mbox{then } &
\left\{
\begin{array}{l}
\dot{\xi}_{n+1} \= - \,  \dot{\xi}_{n} \\ \\
\dot{\eta}_{n+1} \= \dot{\eta}_{n} \, + \, 2 \,
\dfrac{1-k^2}{1+k^2} \, \dot{\xi}_{n}
\end{array}\right.
\end{array}
\end{eqnarray}
\begin{eqnarray}\label{regola_xi_eta4}
\hskip -2.8truecm \begin{array}{lllll} \mbox{If } (\dot{\xi}_{n},
\dot{\eta}_{n}) \in \Z_{12}\, ,  & \mbox{that is if } & \left\{
\begin{array}{l}
\dot{\xi}_n > 0 \\
\dot{\eta}_n > 0
\end{array}\right. \, ,  & \mbox{then } &
\left\{
\begin{array}{lcl}
\dot{\xi}_{n+1} &\=& - \, \dfrac{(1+k^2)(\dot{\xi}_n^2 +
\dot{\eta}_n^2) + 2
(1-k^2)\dot{\xi}_n\dot{\eta}_n}{(1+k^2)(\dot{\xi}_n^2 +
\dot{\eta}_n^2) - 2 (1-k^2)\dot{\xi}_n\dot{\eta}_n} \, \dot{\xi}_{n} \\
\\ && \, + \, 2 \, \dfrac{(1-k^2)(\dot{\xi}_n^2 + \dot{\eta}_n^2)
}{(1+k^2)(\dot{\xi}_n^2 + \dot{\eta}_n^2) - 2
(1-k^2)\dot{\xi}_n\dot{\eta}_n} \, \dot{\eta}_{n}
 \\ \\
\dot{\eta}_{n+1} &\=& 2 \, \dfrac{(1-k^2)(\dot{\xi}_n^2 +
\dot{\eta}_n^2) }{(1+k^2)(\dot{\xi}_n^2 + \dot{\eta}_n^2) - 2
(1-k^2)\dot{\xi}_n\dot{\eta}_n} \, \dot{\xi}_{n} \\ \\
&& \, - \, \dfrac{(1+k^2)(\dot{\xi}_n^2 + \dot{\eta}_n^2) + 2
(1-k^2)\dot{\xi}_n\dot{\eta}_n}{(1+k^2)(\dot{\xi}_n^2 +
\dot{\eta}_n^2) - 2 (1-k^2)\dot{\xi}_n\dot{\eta}_n} \,
\dot{\eta}_{n}
\end{array}\right.
\end{array}
\end{eqnarray}
\end{subequations}

\section{Theoretical results about the ideal impact}

Several results and some remarks can be listed about TA${}_{id}$.
Some of them can be straightforwardly obtained by one or more of
the expressions of the algorithm, some others requires a detailed
proof.

\begin{Result} If ${\bf v}_n \in \Z_1$  then ${\bf v}_{n+1} \in \Z_2$ or
${\bf v}_{n+1} \in \Z_0$. Analogously, if ${\bf v}_n \in \Z_2$
then ${\bf v}_{n+1} \in \Z_1$ or ${\bf v}_{n+1} \in \Z_0$.
\end{Result}

\begin{Proof} It follows immediately from
(\ref{regola_xi_eta2},\ref{regola_xi_eta3}). If ${\bf v}_n \in
\Z_1$  then $\eta_n
> 0$. Then $\eta_{n+1} = - \, \eta_n < 0$, so that ${\bf v}_{n+1}
\in \Z_2$ or ${\bf v}_{n+1} \in \Z_0$. The proof is analogous if
${\bf v}_n \in \Z_2$.
\end{Proof}

This shows that, if an iteration of TA${}_{id}$ gives a velocity
${\bf v} \notin \Z_{12}$, then all the following velocities do not
belong to $\Z_{12}$. In particular, if ${\bf v}_0 \notin \Z_{12}$,
than the evolution of the disk will be determined by a sequence of
single impacts, without multiple impacts.

\begin{Result} If ${\bf v}_n \in \Z_1$ and $k\ge 1$  then ${\bf v}_{n+1}
\in \Z_0$. Analogously, if ${\bf v}_n \in \Z_2$ and $k\ge 1$  then
${\bf v}_{n+1} \in \Z_0$. \end{Result}

\begin{Proof} It follows once again from
(\ref{regola_xi_eta2},\ref{regola_xi_eta3}). If ${\bf v}_n \in
\Z_1$  then $\xi_n \le 0$ and $\eta_n
> 0$. Therefore, if $k\ge 1$, we have $\xi_{n+1} = \xi_n + \, 2 \,
\dfrac{1-k^2}{1+k^2} \eta_n < 0$, so that ${\bf v}_{n+1} \in
\Z_0$. The proof is analogous if ${\bf v}_n \in \Z_2$.
\end{Proof}

This shows that if the angle $2\alpha \ge \frac{\pi}2$ and the
impact is not multiple we have only one iteration of TA${}_{id}$.
This is the case, for instance, when $2\alpha \ge \frac{\pi}2$,
the disk moves along one of the wall and impacts the other wall.

\begin{Result} If ${\bf v}_n \in \Z_1$  then there exists $\chi \in
\Naturali$ such that ${\bf v}_{n+\chi} \in \Z_0$. Analogously, if
${\bf v}_n \in \Z_2$ then there exists $\chi \in \Naturali$ such
that ${\bf v}_{n+\chi} \in \Z_0$. \end{Result}

\begin{Proof} This is a standard proof about reflections following from
(\ref{regola_cosfi_sinfi2},\ref{regola_cosfi_sinfi3}). If ${\bf
v}_n \in \Z_1$ or $\Z_2$ and $k\ge1$ the thesis follows from the
point 2) above.

If ${\bf v}_n \in \Z_1$ and $k<1$, then $\alpha \in
(0,\frac{\pi}{4})$ and $\vph_n \in (- \pi +\alpha, - \alpha]$. We
can construct the odd and even subsequences of the sequence
$\vph_{n+r}$ with $r\in \Naturali$. We have that:
\begin{eqnarray*}
\left\{
\begin{array}{lcl}
\vph_{n+2r} &\=& \vph_n - 2 \, (2r) \, \alpha\\
\vph_{n+2r+1} &\=& - \vph_n + 2 \, (2r+1) \, \alpha
\end{array} \right. .
\end{eqnarray*}
Then $\chi$ is the first natural number such that $\vph_n - 2 \,
(2\chi) \, \alpha \in (-\pi, -\pi + \alpha]\cup[\pi - \alpha,
\pi]$ or $- \vph_n + 2 \, (2\chi+1) \, \alpha \in (-\pi, -\pi +
\alpha]\cup[\pi, \pi - \alpha]$. The proof is analogous if ${\bf
v}_n \in \Z_2$.
\end{Proof}

This shows that, if an iteration of TA${}_{id}$ gives a velocity
${\bf v} \notin \Z_{12}$, then TA${}_{id}$ terminates, giving a
final exit velocity for the disk. Note moreover that the
reflection procedure of this situation is conceptually identical
to the well known one governing the (alternated) single impacts of
a disk with the walls of a corner in a sequence of configurations
of single (and not multiple) contacts between disk and walls.

\smallskip

The three results above pertain TA${}_{id}$ applied in the case of
single impact of the disk in the corner. However the most
significant results are those about multiple impacts. Note that
the condition ${\bf v}_n \in \Z_{12}$ implies that $\dot{x}_n >0$
and $\cos \vph_n
>0$. We have that:

\begin{Result} If ${\bf v}_n \in \Z_{12}$ has the direction of the angle
bisector, then ${\bf v}_{n+1} \= - \, {\bf v}_n \in \Z_0$.
\end{Result}

\begin{Proof} It follows immediately from (\ref{regola_x_y4}) requiring
$\dot{y}_n=0$ or alternatively from (\ref{regola_cosfi_sinfi4})
requiring $\cos \vph_n =1, \sin \vph_n=0$.
\end{Proof}

The main result about multiple impacts is however the following:
\begin{Theorem}\label{PrincipaleIdeale}
If ${\bf v}_n \in \Z_{12}$ then ${\bf v}_{n+1} \notin \Z_{12}$.
\end{Theorem}
\begin{Proof}
Let ${\bf v}_n \in \Z_{12}$ and let us consider
(\ref{regola_xi_eta4}). Then $\dot{\xi_n} > 0, \dot{\eta}_n >0$.
If $\dot{\xi_n} = \dot{\eta}_n$ then ${\bf v}_n$ has the direction
of the angle bisector and we have the thesis. Otherwise note that
\begin{eqnarray*}
(1+k^2)(\dot{\xi}_n^2 + \dot{\eta}_n^2) - 2
(1-k^2)\dot{\xi}_n\dot{\eta}_n \= (\dot{\xi}_n - \dot{\eta}_n)^2 +
k^2 \, (\dot{\xi}_n + \dot{\eta}_n)^2 >0 \,
\end{eqnarray*}
Defining for brevity $\rho_n \= \dfrac{2\xi_n\eta_n}{\dot{\xi}_n^2
+ \dot{\eta}_n^2}, \rho_n \in (0,1)$, we have then
\begin{eqnarray*}
(1+k^2) - \rho_n \, (1-k^2) >0
\end{eqnarray*}
and (13) can be rewritten as
\begin{eqnarray}
\left\{
\begin{array}{lcl}
\dot{\xi}_{n+1} &\=& - \, \dfrac{(1+k^2) + \rho_n \,
(1-k^2)}{(1+k^2) - \rho_n (1-k^2)} \, \dot{\xi}_{n}  \, + \, 2 \,
\dfrac{(1-k^2)}{(1+k^2) - \rho_n (1-k^2)} \, \dot{\eta}_{n}
 \\ \\
\dot{\eta}_{n+1} &\=& 2 \, \dfrac{(1-k^2)}{(1+k^2) - \rho_n
(1-k^2)} \, \dot{\xi}_{n}  \, - \, \dfrac{(1+k^2) + \rho_n \,
(1-k^2)}{(1+k^2) - \rho_n (1-k^2)} \, \dot{\eta}_{n}
\end{array}\right.
\end{eqnarray}
and so:
\begin{eqnarray}
\left\{
\begin{array}{lcl}
\dot{\xi}_{n+1} &\=& \dfrac{\dot{\xi}_n}{(1+k^2) - \rho_n (1-k^2)}
\left(- (1+k^2) \, - \, \rho_n (1-k^2) \, + \, 2\,
\dfrac{\dot{\eta}_n}{\dot{\xi}_n}(1-k^2)\right)
 \\ \\
\dot{\eta}_{n+1} &\=& \dfrac{\dot{\eta}_n}{(1+k^2) - \rho_n
(1-k^2)} \left(- (1+k^2) \, - \, \rho_n (1-k^2) \, + \, 2 \,
\dfrac{\dot{\xi}_n}{\dot{\eta}_n}(1-k^2)\right)
\end{array}\right.
\end{eqnarray}
Suppose by contradiction that ${\bf v}_{n+1} \in \Z_{12}$: then we
must have $\dot{\xi}_{n+1}>0, \dot{\eta}_{n+1}>0$, that is
\begin{eqnarray}
\left\{
\begin{array}{lcl}
 \left(2\, \dfrac{\dot{\eta}_n}{\dot{\xi}_n} \, - \, \rho_n\right)
(1-k^2) &>& (1+k^2)
 \\ \\
\left(2\, \dfrac{\dot{\xi}_n}{\dot{\eta}_n} \, - \, \rho_n\right)
(1-k^2) &>& (1+k^2)
\end{array}\right.
\end{eqnarray}
and recalling the expression of $\rho_n$,
\begin{eqnarray}
\left\{
\begin{array}{lcl}
2 \, (1-k^2)\dfrac{\dot{\eta}_n}{\dot{\xi}_n} \,
\dfrac{\dot{\eta}_n^2}{\dot{\xi}_n^2 + \dot{\eta}_n^2}&>& (1+k^2)
 \\ \\
2 \, (1-k^2)\dfrac{\dot{\xi}_n}{\dot{\eta}_n} \,
\dfrac{\dot{\xi}_n^2}{\dot{\xi}_n^2 + \dot{\eta}_n^2}&>& (1+k^2)
\end{array}\right.
\end{eqnarray}
This is possible only if $1-k^2 > 0$, so that let $k\in (0,1)$.
The function $f(k) = \dfrac{1+k^2}{2(1-k^2)} > \dfrac12$ if $k\in
(0,1)$. Then a necessary condition for ${\bf v}_{n+1} \in \Z_{12}$
is:
\begin{eqnarray}
\left\{
\begin{array}{lcl}
\dfrac{\dot{\eta}_n}{\dot{\xi}_n} \,
\dfrac{1}{\dfrac{\dot{\xi}_n^2 +
\dot{\eta}_n^2}{\dot{\eta}_n^2}}&>& \dfrac12
 \\ \\
\dfrac{\dot{\xi}_n}{\dot{\eta}_n} \,
\dfrac{1}{\dfrac{\dot{\xi}_n^2 +
\dot{\eta}_n^2}{\dot{\xi}_n^2}}&>& \dfrac12
\end{array}\right.
\end{eqnarray}
Let now for simplicity be $\frac{\dot{\eta}_n}{\dot{\xi}_n}= z>0$.
A straightforward calculation shows that (18) is equivalent to:
\begin{eqnarray}
\left\{
\begin{array}{lcl}
z \, \dfrac{1}{\dfrac{1}{z^2} +1}&>& \dfrac12
 \\ \\
\dfrac{1}{z} \, \dfrac{1}{1+z^2}&>& \dfrac12
\end{array}\right.
\Leftrightarrow \left\{
\begin{array}{lcl}
2 \, z^3 \, - \, z^2 \, - \, 1 &>& 0
\\
z^3 \, + \, z \, - \, 2 &<& 0
\end{array}\right.
\Leftrightarrow \left\{
\begin{array}{lcl}
z &>& 1
\\
 z  &<& 1
\end{array}\right. \, .
\end{eqnarray}
This is not possible, and then ${\bf v}_{n+1} \notin \Z_{12}$.
\end{Proof}
We have then the following:
\begin{Corollary}
For every initial velocity ${\bf v}_0$ of the disk, the algorithm
terminates after a finite number of steps.
\end{Corollary}
\begin{Proof}
If ${\bf v}_0 \in \Z_0$, there is nothing to proof. If ${\bf
v}_0\in \Z_1$ or ${\bf v}_0 \in \Z_2$, the thesis follows from
Result 3 of this section. If ${\bf v}_0 \in \Z_{12}$, then ${\bf
v}_1 \notin \Z_{12}$ and then we have the thesis.
\end{Proof}

\begin{Result}
Theorem \ref{PrincipaleIdeale} and Result 1 imply that the system
can have at most one multiple impact if and only if ${\bf v}_0 \in
\Z_{12}$.
\end{Result}

\section{Theoretical algorithm for the non--ideal case (TA${}_{nid}$)}

In this section we present  two different forms of the iterative
rule assigning the ``new'' velocity ${\bf v}_{n+1}$ of the disk as
function of the ``old'' velocity ${\bf v}_{n}$ in the non--ideal
case. The rule is derived by the theoretical characterization of
non--ideal impact presented in Appendix. Each one of the forms
will be used to obtain theoretical results about TA${}_{nid}$.

\subsection{First expression of TA${}_{nid}$: use of $(\dot{x},
\dot{y})$}

Given an initial velocity ${\bf v}_0 = (\dot{x}_0, \dot{y}_0)$,
the iterative rule is such that:
\begin{subequations}\label{AlgNonIdealXY}
\begin{eqnarray}\label{regola_x_y_bis1}
\hskip -5.5truecm \begin{array}{lllll} \mbox{If } (\dot{x}_{n},
\dot{y}_{n}) \in \Z_0 \, , & \mbox{that is if } & \left\{
\begin{array}{l}
k\dot{x}_n + \dot{y}_n\le 0 \\
k\dot{x}_n - \dot{y}_n\le 0
\end{array}\right. \, , & \mbox{then } &
\left\{
\begin{array}{l}
\dot{x}_{n+1} \= \dot{x}_{n} \\
\dot{y}_{n+1} \= \dot{y}_{n}
\end{array}\right.
\end{array}
\end{eqnarray}
\begin{eqnarray}\label{regola_x_y_bis2}
\hskip -3truecm \begin{array}{lllll} \mbox{If } (\dot{x}_{n},
\dot{y}_{n}) \in \Z_1 \, , & \mbox{that is if } & \left\{
\begin{array}{l}
k\dot{x}_n + \dot{y}_n\le 0 \\
k\dot{x}_n - \dot{y}_n > 0
\end{array}\right. \, , & \mbox{then } &
\left\{
\begin{array}{l}
\dot{x}_{n+1} \= \dfrac{1-\veps k^2}{1+k^2} \, \dot{x}_{n} +
\dfrac{(1+\veps)k}{1+k^2} \, \dot{y}_{n} \\ \\
\dot{y}_{n+1} \= \dfrac{(1+\veps)k}{1+k^2} \, \dot{x}_{n} -
\dfrac{\veps-k^2}{1+k^2} \, \dot{y}_{n}
\end{array}\right.
\end{array}
\end{eqnarray}
\begin{eqnarray}\label{regola_x_y_bis3}
\hskip -3truecm \begin{array}{lllll} \mbox{If } (\dot{x}_{n},
\dot{y}_{n}) \in \Z_2 \, , & \mbox{that is if } & \left\{
\begin{array}{l}
k\dot{x}_n + \dot{y}_n > 0 \\
k\dot{x}_n - \dot{y}_n \le 0
\end{array}\right. \, , & \mbox{then } &
\left\{
\begin{array}{l}
\dot{x}_{n+1} \= \dfrac{1-\veps k^2}{1+k^2} \, \dot{x}_{n} -
\dfrac{(1+\veps)k}{1+k^2} \, \dot{y}_{n} \\ \\
\dot{y}_{n+1} \= - \dfrac{(1+\veps)k}{1+k^2} \, \dot{x}_{n} -
\dfrac{\veps-k^2}{1+k^2} \, \dot{y}_{n}
\end{array}\right.
\end{array}
\end{eqnarray}
\begin{eqnarray}\label{regola_x_y_bis4}
\hskip -3truecm \begin{array}{lllll} \mbox{If } (\dot{x}_{n},
\dot{y}_{n}) \in \Z_{12}  \, , & \mbox{that is if } & \left\{
\begin{array}{l}
k\dot{x}_n + \dot{y}_n > 0 \\
k\dot{x}_n - \dot{y}_n > 0
\end{array}\right. \, , & \mbox{then } &
\left\{
\begin{array}{l}
\dot{x}_{n+1} \= \dfrac{-\veps k^4\dot{x}_n^2 +
(1-(1+\veps)k^2)\dot{y}_n^2}{k^4\dot{x}_n^2 + \dot{y}_n^2} \,
\dot{x}_{n}
 \\ \\
\dot{y}_{n+1} \= \dfrac{k^2(k^2-(1+\veps))\dot{x}_n^2 - \veps
\dot{y}_n^2}{k^4\dot{x}_n^2 + \dot{y}_n^2} \, \dot{y}_{n}
\end{array}\right.
\end{array}
\end{eqnarray}
\end{subequations}

\begin{Remark} The iterative rule (\ref{regola_x_y_bis1}--\ref{regola_x_y_bis4}) once
again respects the symmetry of the mechanical problem with respect
to the $(\dot{x},\dot{y})$ components of the velocity, since
relations (\ref{rispetto_simmetria}) hold for every ${\bf v}_n
\notin \Z_0$.
\end{Remark}

\subsection{Second expression of TA${}_{nid}$: use of $(\dot{\xi},
\dot{\eta})$}

The algorithm can be expressed once again by using the coordinates
$(\xi, \eta) \= (k x + y,k x - y)$. We obtain:
\begin{subequations}\label{AlgNonIdealXiEta}
\begin{eqnarray}\label{regola_xi_eta_bis1}
\hskip -6.2truecm \begin{array}{lllll} \mbox{If } (\dot{\xi}_{n},
\dot{\eta}_{n}) \in \Z_0  , & \mbox{that is if } & \left\{
\begin{array}{l}
\dot{\xi}_n\le 0 \\
\dot{\eta}_n \le 0
\end{array}\right.  ,  & \mbox{then } &
\left\{
\begin{array}{l}
\dot{\xi}_{n+1} \= \dot{\xi}_{n} \\
\dot{\eta}_{n+1} \= \dot{\eta}_{n}
\end{array}\right.
\end{array}
\end{eqnarray}
\begin{eqnarray}\label{regola_xi_eta_bis2}
\hskip -3truecm \begin{array}{lllll} \mbox{If } (\dot{\xi}_{n},
\dot{\eta}_{n}) \in \Z_1  , & \mbox{that is if } & \left\{
\begin{array}{l}
\dot{\xi}_n\le 0 \\
\dot{\eta}_n > 0
\end{array}\right.  , & \mbox{then } &
\left\{
\begin{array}{l}
\dot{\xi}_{n+1} \= \dot{\xi}_{n}  +
(1+\veps)  \dfrac{1-k^2}{1+k^2}  \dot{\eta}_{n} \\ \\
\dot{\eta}_{n+1} \= -  \veps  \dot{\eta}_{n}
\end{array}\right.
\end{array}
\end{eqnarray}
\begin{eqnarray}\label{regola_xi_eta_bis3}
\hskip -3truecm \begin{array}{lllll} \mbox{If } (\dot{\xi}_{n},
\dot{\eta}_{n}) \in \Z_2 ,  & \mbox{that is if } & \left\{
\begin{array}{l}
\dot{\xi}_n > 0 \\
\dot{\eta}_n \le 0
\end{array}\right. ,  & \mbox{then } &
\left\{
\begin{array}{l}
\dot{\xi}_{n+1} \= -  \veps  \dot{\xi}_{n} \\ \\
\dot{\eta}_{n+1} \= \dot{\eta}_{n}  +  (1+\veps)
\dfrac{1-k^2}{1+k^2}  \dot{\xi}_{n}
\end{array}\right.
\end{array}
\end{eqnarray}
\begin{eqnarray}\label{regola_xi_eta_bis4}
\hskip -3.5truecm \begin{array}{lllll} \mbox{If } (\dot{\xi}_{n},
\dot{\eta}_{n}) \in \Z_{12} ,  & \mbox{that is if } & \left\{
\begin{array}{l}
\dot{\xi}_n > 0 \\
\dot{\eta}_n > 0
\end{array}\right.  ,  & \mbox{then } &
\left\{
\begin{array}{lcl}
\dot{\xi}_{n+1} &=& -  \dfrac{\veps(1+k^2)(\dot{\xi}_n^2 +
\dot{\eta}_n^2) + 2
(1-k^2)\dot{\xi}_n\dot{\eta}_n}{(1+k^2)(\dot{\xi}_n^2 +
\dot{\eta}_n^2) - 2 (1-k^2)\dot{\xi}_n\dot{\eta}_n}  \dot{\xi}_{n} \\
\\ &&  +  (1+\veps)  \dfrac{(1-k^2)(\dot{\xi}_n^2 + \dot{\eta}_n^2)
}{(1+k^2)(\dot{\xi}_n^2 + \dot{\eta}_n^2) - 2
(1-k^2)\dot{\xi}_n\dot{\eta}_n}  \dot{\eta}_{n}
 \\ \\
\dot{\eta}_{n+1} &=& (1+\veps)  \dfrac{(1-k^2)(\dot{\xi}_n^2 +
\dot{\eta}_n^2) }{(1+k^2)(\dot{\xi}_n^2 + \dot{\eta}_n^2) - 2
(1-k^2)\dot{\xi}_n\dot{\eta}_n}  \dot{\xi}_{n} \\ \\
&&  -  \dfrac{\veps(1+k^2)(\dot{\xi}_n^2 + \dot{\eta}_n^2) + 2
(1-k^2)\dot{\xi}_n\dot{\eta}_n}{(1+k^2)(\dot{\xi}_n^2 +
\dot{\eta}_n^2) - 2 (1-k^2)\dot{\xi}_n\dot{\eta}_n} \dot{\eta}_{n}
\end{array}\right.
\end{array}
\end{eqnarray}
\end{subequations}

\begin{Remark} Note that, due to the change of coordinates
(\ref{cambiamento_coordinate}), we have:
\begin{eqnarray*}
\left(\| {\bf v}_{n}\|_{{}_2}\right)^2 \=
\dfrac{1+k^2}{4k^2}\left( \dot{\xi}_n^2 + \dot{\eta}_n^2 \right) +
\dfrac{1-k^2}{2k^2} \, \dot{\xi}_n \dot{\eta}_n
\end{eqnarray*} \end{Remark}

\section{Theoretical results about the non--ideal impact}

Several results and some remarks that can be listed about
TA${}_{nid}$ are strictly analogous to those about TA${}_{id}$.
For instance, Results 1, 2 and 4  can be immediately generalized
to the non--ideal case, with proofs and remarks analogous to those
presented in Sec.3. Theorem \ref{PrincipaleIdeale} too holds in
the non--ideal case, as we prove below in this section. Instead,
in the non--ideal case we cannot state the analogous of Result 3
of Sec.3, that in the ideal case is crucial to prove that
TA${}_{id}$ terminates. However, for TA${}_{nid}$, the termination
will be ensured on the basis of the criterion $\lim_{n\to +\infty}
\|{\bf v}_n\|_{_{2}} = 0$ of the following theorem
(\ref{ConvergenzaAZero}).

%%% INSERIRE EVENTUALMENTE QUI PARTI DA ScartiDaReinserire.tex %%%

\begin{Theorem}\label{PrincipaleNonIdeale}
If ${\bf v}_n \in \Z_{12}$ then ${\bf v}_{n+1} \notin \Z_{12}$.
\end{Theorem}

\begin{Proof}
Let ${\bf v}_n$ be in $\Z_{12}$, so that $\dot{\xi}_n >0$ and
$\dot{\eta}_n >0$. Recalling that $\veps \in [0,1)$ and $k\in
(0,+\infty)$, we set
\begin{eqnarray*}
\beta \= \dfrac{1-k^2}{1+k^2} \, \in (-1,1) , \, \quad z \=
\dfrac{\dot{\eta}_n}{\dot{\xi}_n} \, \in (0,+\infty) \, .
\end{eqnarray*}
We have $\dfrac{2\dot{\xi}_n\dot{\eta}_n}{\dot{\xi}_n^2 +
\dot{\eta}_n^2} \= \dfrac{2z}{1+z^2}$ and eqs.
(\ref{regola_xi_eta_bis4}) can be rewritten as
\begin{eqnarray*}
\left\{
\begin{array}{lcl}
\dot{\xi}_{n+1} &=& \dfrac{\dot{\xi}_n}{1-\beta \dfrac{2z}{1+z^2}}
\, \left( -\veps + \beta\left((1+\veps)z -
\dfrac{2z}{1+z^2}\right) \right) \\ \\
\dot{\eta}_{n+1} &=& \dfrac{\dot{\eta}_n}{1-\beta
\dfrac{2z}{1+z^2}} \, \left( -\veps +
\beta\left((1+\veps)\dfrac{1}{z} - \dfrac{2z}{1+z^2}\right)
\right)
\end{array}\right.
\end{eqnarray*}
where the two first factors of the RHSs are positive. Then ${\bf
v}_{n+1} \notin \Z_{12}$ if and only if the system of inequalities
\begin{eqnarray}\label{Sistema_NonIdeale}
\left\{
\begin{array}{l}
-\veps + \beta\left((1+\veps)z -
\dfrac{2z}{1+z^2}\right) \, > \, 0 \\ \\
-\veps + \beta\left((1+\veps)\dfrac{1}{z} -
\dfrac{2z}{1+z^2}\right) \, > \, 0
\end{array}\right.
\end{eqnarray}
does not admit solutions for $\veps \in [0,1), \beta \in (-1,1), z
\in (0,+\infty)$. Obviously (\ref{Sistema_NonIdeale}) does not
have solutions if $\beta = 0$ (that is when the amplitude
$2\alpha$ of the corner is $\frac{\pi}2$), if $z=1$ (that is when
${\bf v}_{n}$ is along the bisector of the corner), if $(1+\veps)z
- \frac{2z}{1+z^2}=0$, if $(1+\veps)\frac{1}{z} -
\frac{2z}{1+z^2}=0$.

\smallskip

If $\veps =0, \beta\in(0,1)$, then (\ref{Sistema_NonIdeale}) is
equivalent to
\begin{eqnarray*}
\left\{
\begin{array}{l}
z - \dfrac{2z}{1+z^2} \, > \, 0 \\ \\
\dfrac{1}{z} - \dfrac{2z}{1+z^2} \, > \, 0
\end{array}\right.
\qquad \Rightarrow \qquad \left\{
\begin{array}{l}
z^2 -1 \, > \, 0 \\ \\
1-z^2 \, > \, 0 \, \,
\end{array}\right.
\end{eqnarray*}
that does not have solutions $\forall \, z \in (0,+\infty)$. If
$\veps =0, \beta\in (-1,0)$, then (\ref{Sistema_NonIdeale}) is
equivalent to
\begin{eqnarray*}
\left\{
\begin{array}{l}
z - \dfrac{2z}{1+z^2} \, < \, 0 \\ \\
\dfrac{1}{z} - \dfrac{2z}{1+z^2} \, < \, 0
\end{array}\right.
\qquad \Rightarrow \qquad \left\{
\begin{array}{l}
z^2 -1 \, < \, 0 \\ \\
1-z^2 \, < \, 0 \, \,
\end{array}\right.
\end{eqnarray*}
that does not have solutions $\forall \, z \in (0,+\infty)$. If
$\veps \in (0,1), \beta \in (-1,0)$, then
(\ref{Sistema_NonIdeale}) implies that
\begin{eqnarray*}
\left\{
\begin{array}{l}
(1+\veps)z -
\dfrac{2z}{1+z^2}\, < \, 0 \\ \\
(1+\veps)\dfrac{1}{z} - \dfrac{2z}{1+z^2} \, < \, 0
\end{array}\right.
\qquad \Rightarrow \qquad \left\{
\begin{array}{l}
\veps < \dfrac{1- z^2}{1+z^2} \\ \\
\veps < \, - \, \dfrac{1- z^2}{1+z^2} \, \, ,
\end{array}\right.
\end{eqnarray*}
that does not have solutions $\forall \, z \in (0,+\infty)$.

\smallskip

If $\veps \in (0,1), \beta \in (0,1), z\in(0,1)$, then the first
inequality of (\ref{Sistema_NonIdeale}) can be verified only if
$(1+\veps)z - \frac{2z}{1+z^2}>0$. In this case we have:
\begin{eqnarray*}
\beta > \dfrac{\veps}{(1+\veps)z - \dfrac{2z}{1+z^2}} \=
\dfrac{\veps(1+z^2)}{(1+\veps)(z+z^3) - 2z}
\end{eqnarray*}
This can happen only if
\begin{eqnarray*}
 \dfrac{\veps(1+z^2)}{(1+\veps)(z+z^3) - 2z} <1 \qquad
 \Leftrightarrow \qquad \veps < \, - \,
 \dfrac{z(1+z)}{1+z^2} <0
\end{eqnarray*}
that is not admissible. If $\veps \in (0,1), \beta \in (0,1),
z\in(1,+\infty)$, then the second inequality of
(\ref{Sistema_NonIdeale}) can be verified only if
$(1+\veps)\frac{1}{z} - \frac{2z}{1+z^2}>0$. In this case we have:
\begin{eqnarray*}
\beta > \dfrac{\veps}{(1+\veps)\frac{1}{z} - \dfrac{2z}{1+z^2}} \=
\dfrac{\veps\, z (1+z^2)}{(1+\veps)(1+z^2) - 2z^2}
\end{eqnarray*}
This can happen only if
\begin{eqnarray*}
\dfrac{\veps\, z (1+z^2)}{(1+\veps)(1+z^2) - 2z^2} <1 \qquad
 \Leftrightarrow \qquad \veps < \, - \,
 \dfrac{1+z}{1+z^2} <0
\end{eqnarray*}
that is not admissible. It follows that (\ref{Sistema_NonIdeale})
cannot have solutions, and then ${\bf v}_{n+1} \notin \Z_{12}$.
\end{Proof}

\begin{Result}
Theorem \ref{PrincipaleNonIdeale} and Result 1 (that holds for
non--ideal impacts too) imply once again that the system can have
at most one multiple impact if and only if ${\bf v}_0 \in
\Z_{12}$.
\end{Result}

To proof the second important result about TA${}_{nid}$ we need to
introduce the convergent matrices and their properties. Let $A$ be
an  $N \times N$  matrix and let $\rho(A)$ its spectral radius,
that is the largest modulus of its eigenvalues. We recall that the
matrix $A$ is said to be convergent if $\lim_{k \rightarrow
+\infty} (A^k)_{ij} =0$ for each $i, j =1, \dots, N$, where
$(A^k)_{ij} $ is the $(i,j)$-th element of $A^k$. The following
theorem (see e.g. Theorem 4 in \cite[p. 14]{IsaacsonK}) shows some
well-known properties of  a convergent matrix.
\begin{Theorem}\label{ThIK94}
The following three statements are equivalent:
\begin{enumerate}
\item the matrix $A$ is convergent; \item  $\lim_{k \rightarrow
+\infty} \|A^k\| =0$ for some matrix norm induced by a vector
norm, that is defined by $\|A\|=\max_{\|x\|=1}\|Ax\|$; \item
$\rho(A) <1$.
\end{enumerate}
\end{Theorem}
\begin{Remark}
\label{RemIK94} Let $A$ be a convergent matrix and let $\| \ \cdot
\  \|$ be the induced matrix norm for which item 2 of Theorem
\ref{ThIK94} holds. Given a vector ${\bf w}$, we have, from a
property of the induced matrix norm,  that $ 0 \le \|A^k {\bf w}
\| \le \|A^k \|\|{\bf w}\| $ and so $\lim_{k \rightarrow +\infty}
\|A^k {\bf w}\|=0$. It follows that the vector $A^k {\bf w}$
converges to the zero vector.
\end{Remark}
%%%%%%%%

\begin{Theorem}\label{ConvergenzaAZero} If ${\bf v}_{n} \notin \Z_{0}$ for every $n$, then
$\lim_{n\to +\infty} \|{\bf v}_n\| = 0$.
\end{Theorem}
\begin{Proof}
Result 6 implies that ${\bf v}_{1} \notin \Z_{12}$. Let once again
be $\beta \= ({1-k^2})/({1+k^2})$. By hypothesis, due to Result 6,
we can take $k\in(0,1)$ and then $\beta \in (0,1)$.

Let us suppose that ${\bf v}_{1} \= (\dot{\xi}_{1},\dot{\eta}_{1})
\in \Z_1$. Result 5 and the hypothesis imply that ${\bf v}_{3} \=
(\dot{\xi}_{3},\dot{\eta}_{3}) \in \Z_1$. Applying
(\ref{regola_xi_eta_bis2},\ref{regola_xi_eta_bis3}) we have
\begin{eqnarray*}
\begin{array}{lcl}
\left(
\begin{array}{c}
\dot{\xi}_{3} \\
\dot{\eta}_{3}
\end{array}\right) & = &
\left(
\begin{array}{cc}
-\veps & -\beta\veps(1+\veps) \\
\beta(1+\veps)&\beta^2(1+\veps)^2 - \veps
\end{array}\right) \,
\left(
\begin{array}{c}
\dot{\xi}_{1} \\
\dot{\eta}_{1}
\end{array}\right)
\end{array}
\end{eqnarray*}
Therefore, for every $h\in \Naturali$, we have
\begin{eqnarray*}
\begin{array}{lcl}
\left(
\begin{array}{c}
\dot{\xi}_{2h+1} \\
\dot{\eta}_{2h+1}
\end{array}\right) & = &
\left(
\begin{array}{cc}
-\veps & -\beta\veps(1+\veps) \\
\beta(1+\veps)&\beta^2(1+\veps)^2 - \veps
\end{array}\right)^h \,
\left(
\begin{array}{c}
\dot{\xi}_{1} \\
\dot{\eta}_{1}
\end{array}\right)
\end{array}
\end{eqnarray*}
Moreover, if ${\bf v}_{1} \in \Z_1$, by Result 5 and the
hypothesis we have that ${\bf v}_{2} \in \Z_2$. A straigthforward
calculation shows that in this case, for every $h\in \Naturali,
h>0$, we have
\begin{eqnarray*}
\begin{array}{lcl}
\left(
\begin{array}{c}
\dot{\xi}_{2h} \\
\dot{\eta}_{2h}
\end{array}\right) & = &
\left(
\begin{array}{cc}
\beta^2(1+\veps)^2 -\veps & \beta(1+\veps) \\
-\beta\veps(1+\veps)& - \veps
\end{array}\right)^h \,
\left(
\begin{array}{c}
\dot{\xi}_{2} \\
\dot{\eta}_{2}
\end{array}\right)
\end{array}
\end{eqnarray*}
Since the two matrices
\begin{eqnarray*}
H_1 \= \left(
\begin{array}{cc}
-\veps & -\beta\veps(1+\veps) \\
\beta(1+\veps)&\beta^2(1+\veps)^2 - \veps
\end{array}\right) \qquad H_2 \= \left(
\begin{array}{cc}
\beta^2(1+\veps)^2-\veps & \beta(1+\veps) \\
-\beta\veps(1+\veps)& - \veps
\end{array}\right)
\end{eqnarray*}
have the same characteristic polynomial and eigenvalues, then
Theorem \ref{ThIK94} and Remark 6 imply  that $\lim_{n\to +\infty}
(\dot{\xi}_{n},\dot{\eta}_{n}) = (0,0)$  if the spectral radius
$\rho(H_1)=\rho(H_2)$ is such that $\rho(H_1) <1$. Therefore the
theorem follows upon proof that $\rho(H_1) <1$. Needless to say,
the proof is completely analogous if ${\bf v}_{1} \in \Z_2$.

\smallskip

The characteristic polynomial of $H_1$ is
\begin{eqnarray*}
p_{H_1}(\lambda) \= \lambda^2 - (\beta^2(1+\veps)^2 - 2\veps)
\lambda + \veps^2 \, ,
\end{eqnarray*}
whose corresponding eigenvalues are
\begin{eqnarray*}
\lambda_{12} \= \frac12\left(\beta^2(1+\veps)^2 - 2\veps \pm
\beta(1+\veps)\sqrt{\beta^2(1+\veps)^2 - 4\veps}\right)
\end{eqnarray*}
where we have $\lambda_1\lambda_2 = \veps^2$.

\smallskip

If $\beta^2(1+\veps)^2 - 4\veps<0$, the eigenvalues are complex
conjugates with the same module and then
$\rho(H_1)=|\lambda_1|=|\lambda_2| = \veps <1$.

\smallskip

If $\beta^2(1+\veps)^2 - 4\veps=0$, then $\lambda_1=\lambda_2 =
\veps <1$.

\smallskip

If $\beta^2(1+\veps)^2 - 4\veps>0$, the eigenvalues are both in
$\Reali$ and they have the same sign. In particular, since
$\beta^2(1+\veps)^2 - 2\veps>0$,  then
$\rho(H_1)=\max\{|\lambda_1|,|\lambda_2|\} \=
\frac12\left(\beta^2(1+\veps)^2 - 2\veps +
\beta(1+\veps)\sqrt{\beta^2(1+\veps)^2 - 4\veps}\right)$.

A standard study of
\begin{eqnarray*}
\rho(H_1)(\veps,\beta) \= \frac12\left(\beta^2(1+\veps)^2 - 2\veps
+ \beta(1+\veps)\sqrt{\beta^2(1+\veps)^2 - 4\veps}\right)
\end{eqnarray*}
in the compact set $\overline{\Theta} = \left\{(\veps,\beta) \big|
\veps \in [0,1], \beta \in \left[\dfrac{2\sqrt{\veps}}{1+\veps},
1\right]\right\}$ shows that, since $\dfrac{\partial
\rho(H_1)}{\partial \beta} >0$, the maximum is taken in the
segment $\{\beta = 1\}$ and $\underset{\overline{\Theta}}{\max}
(\rho(H_1))=1$. Then for every fixed $(\veps,\beta)\in {\Theta} =
\left\{(\veps,\beta) \big| \veps \in (0,1), \beta \in
\left(\dfrac{2\sqrt{\veps}}{1+\veps}, 1\right)\right\}$ we have
$\underset{{\Theta}}{\max}(\rho(H_1))<1$.

\smallskip

In conclusion, for every $\veps \in (0,1), k\in (0,1)$ we have
$\rho(H_1) \in [\veps,1)$. Since $\lim_{n\to +\infty}
\|(\dot{\xi}_{n},\dot{\eta}_{n})\| \= 0 \= \lim_{n\to +\infty}
\|(\dot{x}_{n},\dot{y}_{n})\|$ obviously implies that $\lim_{n\to
+\infty} \|{\bf v}_{n}\|_{_2} = 0$, we have the thesis.
\end{Proof}

\begin{Remark} Let be ${\bf v}_0 \in \Z_{12}$ and $\veps \in (0,1)$. A tedious but
straightforward calculation\footnote{The calculation was helped by
the use of the factorization command of CoCoa${}^{\copyright}$, a
freely available program for computing with multivariate
polynomials.} shows that
\begin{eqnarray*}
\| {\bf v}_{1}\|^2_{{}_2} - \| {\bf v}_{0}\|^2_{{}_2} \= (\veps^2
-1) \, \dfrac{(k^2 \dot{x}_0^2 + \dot{y}_0^2)^2}{k^4 \dot{x}_0^2 +
\dot{y}_0^2} \, < 0 \quad \Rightarrow \quad \| {\bf
v}_{1}\|_{{}_2} < \| {\bf v}_{0}\|_{{}_2}.
\end{eqnarray*}
Moreover, another straightforward calculation shows that
\begin{eqnarray*}\label{rapporto_norme_non_ideale}
\begin{array}{l}
{\bf v}_n \in \Z_1 \, \Rightarrow \| {\bf v}_{n+1}\|_{{}_2} < \|
{\bf v}_{n}\|_{{}_2}\, , \quad {\bf v}_n \in \Z_2 \, \Rightarrow
\| {\bf v}_{n+1}\|_{{}_2} < \| {\bf v}_{n}\|_{{}_2} \quad \forall
\, n\ge 0
\end{array}
\end{eqnarray*}
and then the whole sequence $(\| {\bf v}_{n}\|_{{}_2})_{n\ge0}$
decreases to $0$.
\end{Remark}

\begin{Remark} Theorem \ref{ConvergenzaAZero} states a physical property
of the mechanical system and not only a numerical property of
TA${}_{nid}$. For example, the same procedure of the proof applied
starting from the rule (\ref{AlgNonIdealXY}) instead of
(\ref{AlgNonIdealXiEta}) leads to the analysis of the spectral
radius of the matrices
\begin{eqnarray*}
K_1 \= \left(
\begin{array}{cc}
\dfrac{(1-\veps k^2)^2 - k^2(1+\veps)^2}{(1+k^2)^2} & \dfrac{k(1 - k^2)(1+\veps)^2}{(1+k^2)^2} \\
- \dfrac{k(1 - k^2)(1+\veps)^2}{(1+k^2)^2}& \dfrac{(\veps- k^2)^2
- k^2(1+\veps)^2}{(1+k^2)^2}
\end{array}\right) \\ \\
K_2 \= \left(
\begin{array}{cc}
\dfrac{(1-\veps k^2)^2 - k^2(1+\veps)^2}{(1+k^2)^2} & - \dfrac{k(1 - k^2)(1+\veps)^2}{(1+k^2)^2} \\
\dfrac{k(1 - k^2)(1+\veps)^2}{(1+k^2)^2}& \dfrac{(\veps- k^2)^2 -
k^2(1+\veps)^2}{(1+k^2)^2}
\end{array}\right) .
\end{eqnarray*}
It can be easily shown that the matrix $B \= \left(
\begin{array}{cc}
k & 1 \\
k & -1
\end{array}\right)
$ that expresses the change of coordinates
(\ref{cambiamento_coordinate}) is such that $K_1 \= B^{-1}H_1 B,
K_2 \= B^{-1}H_2 B$. The matrices $H_1, H_2$ and $K_1, K_2$ are
then respectively similar, they have the same eigenvalues and then
the same spectral radius. Similar arguments hold for every
admissible change of coordinates.
\end{Remark}

\begin{Remark}
For known results on matrices (the so called Gelfand's formula.
See e.g. Theorem 4 in~\cite[p. 28]{Lax}), the spectral radius $
\rho (H_1)$ of the matrix $H_1$ can be expressed as a limit of
matrix norms, that is
\begin{eqnarray*}
\rho (H_1)=\lim _{h \to \infty }\left\|H_1^{h}\right\|^{\frac
{1}{h}}.
\end{eqnarray*}
It follows that, for a large enough $h$, we have $\|H_1^h\|
\approx \rho(H_1)^h$, and so
\begin{eqnarray*}
\|H_1^h v_1 \| \le \|H_1^h\| \|v_1\| \approx \rho(H_1)^h \|v_1\|.
\end{eqnarray*}
Therefore the spectral radius $\rho(H_1)=\rho(K_1)$ gives also a
measure of the rate of convergence to $0$ of the velocity ${\bf
v}_n$. It follows then from the proof of Theorem
(\ref{ConvergenzaAZero}) that the bigger $\veps$ and $\beta$ are,
the slower the convergence is. This means that we can forecast
slow convergence to $0$ of the velocity for ``almost elastic''
walls and very small angles $\alpha$.
\end{Remark}

\vskip 1truecm

\centerline{\bf {\LARGE PART 2: NUMERICAL ASPECTS}}

\section{The numerical results}
We already said that data uncertainty due to possible errors in
the real--world measurements and algorithmic errors due to the use
of floating point arithmetic can cause disastrous effects on the
result. The theoretical algorithms TA can be easily implemented,
but they present some numerical drawbacks, because of their
sensitivity to the noise on the input data and of its instability
with respect to the floating point arithmetic: small perturbations
of $\dot \xi_n$ and $\dot \eta_n$  can cause structural changes in
response, e.g. it can happen that the exact $\dot \xi_n$ is a
small negative value while the computed $\dot \xi_n$ is a small
positive value, causing a different choice in the iterative
method.

For this reason we present a numerical algorithm obtained by
changing the tests for  the choice of the iterative step to do. A
threshold $S$ is introduced
 to consider as zero value the very small positive $\dot{
\xi }_n$ or $\dot{ \eta }_n$. Analogously, we introduce a test on
the norm of the final linear velocity, in order to consider as
(almost) at rest a disk whose computed velocity is less than a
very small threshold $S_v$.

Obviously, if $S=0$ and $S_v=0$ the  numerical algorithm coincides
with the theoretical ones.

\subsection{Numerical Algorithm (NA)}

\begin{itemize}
\item {\bf Input}: the coefficient $\varepsilon \in [0,1]$, the
angle $\alpha \in (0,\frac{\pi}4)$, the initial  velocity ${\bf
v}_0=(\dot x_0,\dot y_0)$, with $\| {\bf v}_0 \|_{_2}=1$, and the
threshold $S$ and $S_v$.

\item {\bf Output}: the final velocity ${\bf v}_f$.

\item {\bf First step}: $k=\tan(\alpha)$; $n=0$; $\dot \xi_0=k\dot
x_0+ \dot y_0$; $\dot \eta_0=k \dot x_0-\dot y_0$.

\item {\bf Core}: While ($\dot \xi_n >S \text{  or  } \dot \eta_n
>S$), $n<N_{max}$, and $\|(\dot x_n, \dot y_n)\|_{_2} > S_v$:
\begin{enumerate}
\item if ($\dot \xi_n\le S $) and ($\dot \eta_n  > S$)  then
\begin{eqnarray*}
\left\{
\begin{array}{l}
\dot x_{n+1}=\dfrac{1-\varepsilon k^2}{1+k^2}\dot
x_n+\dfrac{(1+\varepsilon)k }{1+k^2} \dot y_n \\ \\
\dot y_{n+1}=\dfrac{(1+\varepsilon) k }{1+k^2}\dot x_n-
\dfrac{\varepsilon-k^2}{1+k^2} \dot y_n   \ ;
\end{array}
\right.
\end{eqnarray*}
\item if ($\dot \xi_n > S $) and ($\dot \eta_n  \le  S$) then
\begin{eqnarray*}
\left\{
\begin{array}{l}
\dot x_{n+1} =\dfrac{1-\varepsilon k^2}{1+k^2}\dot x_n-
\dfrac{(1+\varepsilon)k}{1+k^2}\dot y_n \\ \\
\dot y_{n+1}=-\dfrac{(1+\varepsilon)k}{1+k^2}\dot x_n-
\dfrac{\varepsilon-k^2}{1+k^2}\dot y_n   \ ;
\end{array}
\right.
\end{eqnarray*}
\item  if ($\dot \xi_n > S $) and ($\dot \eta_n  > S$) then
\begin{eqnarray*}
\left\{
\begin{array}{lcl}
\dot x_{n+1}&=&\dfrac{-\varepsilon k^4 \dot
x_n^2+(1-(1+\varepsilon)k^2)\dot y_n ^2}{k^4\dot x_n^2+\dot y_n
^2} \dot x_n \\ \\
\dot y_{n+1} &=& \dfrac{ k^2(k^2-(1+\varepsilon))\dot
x_n^2-\varepsilon \dot y_n ^2}{k^4\dot x_n^2+\dot y_n ^2}\dot y_n
\ .
\end{array}
\right.
\end{eqnarray*}
\item  $\dot \xi_{n+1}=k\dot x_{n+1}+\dot y_{n+1} \,\, $; $\dot
\eta_{n+1}=k\dot x_{n+1}-\dot y_{n+1} \,\, $ ; $n=n+1$.
\end{enumerate}

\item ${\bf v}_f = (\dot x_n, \dot y_n)$.
\end{itemize}
The previous algorithm stops when both $\dot \xi_n$ and $\dot
\eta_n$ are less than $S$, or when the 2-norm of the computed
velocity is less than $S_v$, or when the  number of steps exceeds
the predefined number $N_{max}$ of cycles.

The following simple example shows the sensitivity of TA to the
noise on the input data and to the floating point arithmetic
computation even in a very simple case.

\begin{Example} We consider the  behaviors of
TA and NA when they process the initial velocity ${\bf v}_0=(\dot
x_0, \ \dot y_0)=(\frac 1{ \sqrt{2}},  \ \frac 1 {\sqrt{2}})$,
with $\alpha=\frac \pi 4$ and $\varepsilon =1$. From the
theoretical point of view, TA and NA process the input data in the
same way. Since ${\bf v}_0$ satisfies the condition of $\Z_2$, at
the first iteration both algorithms compute the new velocity ${\bf
v}_1=(\dot x_1, \ \dot y _1)= (- \dot y_0, \ -\dot x_0) = (-\frac
1 {\sqrt{2}},  \ -\frac 1 {\sqrt{2}})$. Since the coordinates of
${\bf v}_1$ satisfies the conditions in $\Z_0$ for both
algorithms, TA and NA stop and ${\bf v}_1$ is the final computed
velocity.

Nevertheless, when the algorithms are implemented, TA suffers from
the data error and the computational approximation, while NA has
the same behavior of the theoretical case, in the absence of
errors. In fact, since the computed  values of $\dot x_1$ and
$\dot y _1$ are perturbed by errors, we obtain $\dot x_1 = -
0.707106781186547$ and $\dot y_1= -0.707106781186548 $,  so that
$\dot \xi_1= -0.707106781186547$ and $\dot \eta_1=
4.440892098500626e-16$. If $S=0$, that is using an implementation
of TA, $\dot \xi_1$ and $\dot \eta_1$ satisfy the conditions of
$\Z_1$, and the algorithm compute a new iteration. Differently, NA
is more robust and, choosing $S=2\cdot eps $, the values $\dot
\xi_1$ and $\dot \eta_1$ satisfy the conditions of $\Z_0$, and the
algorithm stops, as in the theoretical case.
\end{Example}

\subsection{Comparison between TA and NA}

In the following we show that TA and NA compute the same final
velocity in the same number of steps or, even if one of the
algorithms executes more iterations, the final velocities are very
similar. We can conclude that NA is preferable when we deal with
real world measurements, since it produces analogous  final
velocities as TA, but it is  more robust with respect to the
errors on the input data.

\begin{Lemma}\label{small_changes}
Let $( \dot x_n, \dot y_n)$ be the linear velocity  at the current
step. Then:
\begin{itemize}
\item[i)] if $( \dot x_n, \dot y_n) \in \Z_1$ we have $ \quad |
\dot x_{n+1} - \dot x_n | \le \dot \eta_n \quad \text{and} \quad |
\dot y_{n+1} - \dot y_n  | \le 2 \dot \eta_n \ ; $ \item[ii)] if
$( \dot x_n, \dot y_n) \in \Z_2$ we have $ \quad | \dot x_{n+1} -
\dot x_n | \le \dot \xi_n \quad \text{and} \quad  | \dot y_{n+1} -
\dot y_n | \le 2\dot \xi_n \ . $
\end{itemize}
\end{Lemma}
\begin{Proof}
By direct computation, the relations in $\Z_1$ give
\begin{eqnarray*}
\dot x_{n+1} - \dot x_n   =  -\frac{(1+\varepsilon)k}{1+k^2} \dot
\eta_n \quad \text{and} \quad \dot y_{n+1} - \dot y_n   =
\frac{(1+\varepsilon)}{1+k^2} \dot \eta_n
\end{eqnarray*}
so that the thesis follows, since  $\varepsilon, k <1$ and $\dot
\eta_n>0$. Analogously for $\Z_2$, changing the role of $\dot
\xi_n$ and $\dot \eta_n$.
\end{Proof}

\begin{Lemma}\label{diff}
Let  ${\bf v}= (\dot x, \dot y)$  and  ${\bf v_p}= ( \dot x
+\delta_x,\dot y +\delta_y)$ be two velocity vectors such that
both ${\bf v}, {\bf v_p} \in \Z_{1}$ or both ${\bf v}, {\bf v_p}
\in\Z_{2}$. Let ${\bf w}=(\dot t, \dot z)$ and ${\bf w_p}=(\dot
t_p, \dot z_p)$ be the new computed velocity vectors starting from
${\bf v}$ and ${\bf v_p}$, respectively. Then
\begin{itemize}
\item[i)] $ \left | \dot t_p - \dot t \right |\le
|\delta_x|+|\delta_y| \quad \text{and} \quad  \left | \dot z_p -
\dot z \right | \le |\delta_x|+|\delta_y|    \, ; $

\item[ii)] $\left | (k\dot t_p + \dot z_p)- (k\dot t + \dot z)
\right |\le \dfrac 3 2 |\delta_x|+|\delta_y| \quad \text{and}
\quad \\ \phantom{XXX}  \left | (k\dot t_p - \dot z_p)- (k\dot t -
\dot z) \right |\le \dfrac 3 2 |\delta_x|+|\delta_y| \, .$
\end{itemize}
\end{Lemma}
\begin{Proof}
If ${\bf v}, {\bf v_p} \in \Z_{1}$ then, from
(\ref{regola_x_y_bis2}) we have
\begin{equation*}
\dot t_p - \dot t  \= \frac{(1-\varepsilon k^2)}{1+k^2}\delta_x
+\frac{(1+\varepsilon) k}{1+k^2}\delta_y \qquad \textrm{and}
\qquad \dot z_p - \dot z \= \frac{(1+\varepsilon)
k}{1+k^2}\delta_x -\frac{\varepsilon - k^2}{1+k^2}\delta_y \ .
\end{equation*}
so that $i)$ follows from $\varepsilon, k < 1$. Moreover, we have
that
\begin{eqnarray*}
\begin{array}{lcl}
(k\dot t_p + \dot z_p)- (k\dot t + \dot z) &=& \dfrac{k
(1-\varepsilon k^2)}{1+k^2}\delta_x +\dfrac{k^2
(1+\varepsilon)}{1+k^2}\delta_y   +\dfrac{k
(1+\varepsilon)}{1+k^2}\delta_x -\dfrac{\varepsilon -
k^2}{1+k^2}\delta_y
\\    \\
&=& \dfrac{ k  (2 -\varepsilon k^2 +\varepsilon )} {1+k^2}
\delta_x  +   \dfrac{(2k^2 + \varepsilon k^2  -\varepsilon
)}{1+k^2}       \delta_y
\ , \\ \\
(k\dot t_p - \dot z_p)- (k\dot t - \dot z) &=& \dfrac{k
(1-\varepsilon k^2)}{1+k^2}\delta_x +  \dfrac{k^2
(1+\varepsilon)}{1+k^2}\delta_y - \dfrac{k (1+\varepsilon)}{1+k^2}
\delta_x +\dfrac{\varepsilon - k^2}{1+k^2} \delta_y
\\       \\ &=& -k \varepsilon \delta_x +\varepsilon  \delta_y \ .
\end{array}
\end{eqnarray*}
Since  $\varepsilon, k \le 1 $, we have
\begin{eqnarray*}
\begin{array}{lcl}
&& 0 \le  \dfrac{ k  (2 -\varepsilon k^2 +\varepsilon )} {1+k^2} =
\dfrac{ k } {1+k^2}  (2 +\varepsilon( 1-k^2  ) ) \le 3\dfrac{ k }
{1+k^2} \le \dfrac 3 2 \qquad \textrm{and} \\ \\
&&\left |\dfrac{(2k^2 + \varepsilon k^2 -\varepsilon  )}{1+k^2}
\right | = \begin{cases}
1 & \text{if } \   k=1  \ ,  \\ \\
\dfrac {2k^2 - \varepsilon(1-k^2)}{1+k^2} \le  \dfrac {2k^2
}{1+k^2} \le1 &\text{if } \
k \neq 1, \varepsilon \le \dfrac{2 k^2}{1-k^2}\ ,  \\ \\
\dfrac {\varepsilon(1-k^2)-2k^2}{1+k^2}   \le \dfrac
{\varepsilon(1-k^2)}{1+k^2} \le 1  &\text{if } \ k \neq 1,
\varepsilon  >  \dfrac{2 k^2}{1-k^2}  \ .
\end{cases}
\end{array}
\end{eqnarray*}
Then $ii)$ follows. Analogous computation holds when  ${\bf v}$
and  ${\bf v_p}$ are in $\Z_2$.
\end{Proof}

\begin{Theorem}
Starting from the same input, that is the same angle $\alpha$, the
same parameter $\varepsilon$ and the same  initial velocity vector
$(\dot x_0, \ \dot y_0)$,  TA and NA either compute the same
output in the same number of steps, or they compute slightly
different outputs, even if one of the algorithms  performs more
steps.
\end{Theorem}
\begin{Proof}
First of all, we observe that, if NA stops because $\|{\bf
v}_f\|_{_2} < S_v$, then the disk is almost at rest and so, even
if TA computes further iterations, its output is similar to  the
final velocity computed by NA. Let $(\dot x_n, \dot y_n)$ be the
velocity vector, with corresponding $\dot \xi_n=k \dot x_n+\dot
y_n$ and $\dot \eta_n =k \dot x_n -\dot y_n$, processed at the
$n$-th step by both TA and NA. This is certainly verified at the
first step, when $n=0$.

\medskip

In the following cases both TA and NA have the same behavior, that
is they compute the same velocity vector $(\dot x_{n+1}, \ \dot
y_{n+1})$.
\begin{itemize}
\item If $\dot \xi_n \le 0$ and $\dot \eta_n \le 0$ both
algorithms stop. \item If $\dot \xi_n  \le 0$ and $\dot \eta_n \ge
S$, or $\dot \xi_n  \ge S $ and $\dot \eta_n \le 0$, or $\dot
\xi_n  \ge S$ and $\dot \eta_n \ge S$, both algorithms compute the
same new velocity vector.
\end{itemize}

In the other cases, we show that the TA and NA have different
behavior, but they compute similar outputs.
\begin{enumerate}
\item If $0< \dot \xi_n < S$ and  $0< \dot \eta_n< S$, then NA
stops and the disk can be considered almost at rest, since $k \gg
S $ and
$$\| (\dot x_n, \ \dot y_n) \|^2_2=\frac{1+k^2}{4k^2}
(\dot \xi_n^2+\dot \eta_n^2) +\frac{1-k^2}{2k^2} \dot \xi_n \dot
\eta_n \le \frac {S^2}{ k^2}  \ . $$ The TA consider $(\dot x_n, \
\dot y_n) \in \Z_{12}$ and computes a new iteration. Since the
$2$-norm of the velocity vector decreases at each step also the
output of the TA corresponds to an almost at rest disk. \item If
$0<\dot \xi_n<S$ and $\dot \eta_n \le 0$, then NA stops, while TA
computes a new velocity vector $(\dot x_{n+1}, \ \dot y_{n+1})$
using the relations in $\Z_2$. From Lemma~\ref{small_changes},
$(\dot x_{n+1}, \ \dot y_{n+1})$ differs from $(\dot x_n, \ \dot
y_n)$, component-wise, for less than $2S$, since $0<\dot \xi_n<S$.
Moreover,  $\dot \xi_{n+1} =-\varepsilon \dot \xi_n$, that is $ -S
<\dot \xi_{n+1}  <0 $ and so, if $\dot \eta_{n+1} $ is negative,
then TA stops and its output is similar to the one of NA.
Otherwise, if $ \dot \eta_{n+1}$ is positive,  since $\dot \eta_n
\le 0$, we have
$$0<  \dot \eta_{n+1} = \dot \eta_n +(1+\varepsilon)\frac{1-k^2}{1+k^2} \dot \xi_n < (1+\varepsilon)\frac{1-k^2}{1+k^2} \dot \xi_n < 2S
$$
and so  the disk is almost at rest, we conclude, as in item 1,
that TA and NA produce similar outputs.

Analogously if $0<\dot \eta_n<S$ and $\dot \xi_n \le 0$.

\item If $0<\dot \xi_n\le S$ and $\dot \eta_n >S $, the NA and TA
process  the velocity vector $(\dot x_n , \ \dot y_n) $ in
different ways. We have analogous behaviors if $\dot \xi_n > S  $
and $0< \dot \eta_n\le S$,  changing the role of $\dot \xi_n$ and
$\dot \eta_n$.

Let $0< \dot \xi_n \le  S$ and  $\dot \eta_n >S$.

The NA  computes $(\dot x_{n+1}^{(1)}, \ \dot x_{n+1}^{(1)})$,
$\dot \xi_{n+1}^{(1)}=k \dot x_{n+1}^{(1)} + \dot y_{n+1}^{(1)} $
and $\dot \eta_{n+1}^{(1)}=k \dot x_{n+1}^{(1)} - \dot
y_{n+1}^{(1)} $, where
\begin{eqnarray*}
 \dot x_{n+1}^{(1)} &=&  \dfrac{\dot
\xi_{n+1}^{(1)} +\dot\eta_{n+1}^{(1)}}{2k}=\dfrac{1}{2k} \dot
\xi_n
+\dfrac{1 -k^2 -2\varepsilon k^2}{2k(1+k^2)} \dot \eta_n\\ \\
\dot y_{n+1}^{(1)} &=& \dfrac{\dot \xi_{n+1}^{(1)} - \dot
\eta_{n+1}^{(1)}}{2}= \dfrac{1}{2} \dot \xi_n +\dfrac{1 -k^2
+2\varepsilon }{2(1+k^2)} \dot \eta_n \ .
\end{eqnarray*}
Furthermore,
$$  \dot \eta_{n+1}^{(1)}=-\varepsilon \dot \eta_n <0
\qquad \text{and} \qquad 0<\dot \xi_{n+1}^{(1)}=\dot \xi_n
+\frac{(1+\varepsilon)(1-k^2)}{1+k^2} \dot \eta_n <S+2 \dot \eta_n
\ . $$

The TA  computes $(\dot x_{n+1}^{(12)}, \ \dot x_{n+1}^{(12)})$,
$\dot \xi_{n+1}^{(12)}=k \dot x_{n+1}^{(12)} + \dot y_{n+1}^{(12)}
$ and $\dot \eta_{n+1}^{(12)}=k \dot x_{n+1}^{(12)} - \dot
y_{n+1}^{(12)} $, where
\begin{eqnarray*}
\dot x_{n+1}^{(12)} &=&   \frac{\dot \xi_{n+1}^{(12)} +\dot \eta_{n+1}^{(12)}}{2k}= (\dot \xi_n+\dot \eta_n)\frac {(1-k^2-2\varepsilon k^2 )(\dot \xi_n^2+\dot \eta_n^2) -2(1-k^2 )\dot \xi_n\dot \eta_n}{8kD} \\
 \\
\dot y_{n+1}^{(12)} &=&   \frac{\dot \xi_{n+1}^{(12)} -\dot
\eta_{n+1}^{(12)}}{2}=
 (\dot \eta_n- \dot \xi_n)\frac {(1-k^2+2\varepsilon)(\dot \xi_n^2+\dot \eta_n^2) +2 (1-k^2 ) \dot \xi_n\dot \eta_n}{8D}\ ,
\end{eqnarray*}
where $D=\dfrac{(\dot \xi_n^2+\dot \eta_n^2) (1+k^2) +2\dot \xi_n
\dot \eta_n(k^2-1)}{4}$. Furthermore,
$$
\dot \eta_{n+1}^{(12)}= \frac{ -\varepsilon (\dot \xi_n^2+\dot
\eta_n^2)(\dot \eta_n(1+k^2)-\dot \xi_n(1-k^2)) -\dot \xi_n(1-
k^2)(\dot \xi_n^2+\dot \eta_n^2) }{4D}
$$
and, since $\dot \eta_n>\dot \xi_n$, then also $\dot
\eta_{n+1}^{(12)} <0 $.

The linear velocities $(\dot x_{n+1}^{(1)}, \dot y_{n+1}^{(1)})$
and $(\dot x_{n+1}^{(12)}, \dot y_{n+1}^{(12)})$ are very similar.
In fact, since  $4D=(\dot \xi_n^2+\dot \eta_n^2) (1+k^2) +2\dot
\xi_n \dot \eta_n(k^2-1) \approx \dot \eta_n^2 (1+k^2)$,
\begin{eqnarray*}
\left |\dot x_{n+1}^{(1)} -\dot x_{n+1}^{(12)}  \right | &=&
\frac{\dot \xi_n k(1+\varepsilon) \left ( \dot \eta_n
^2(3-k^2)+\dot \xi_n ^2(1+k^2)\right )}{4D(1+k^2)}\\ \\
&\approx&\frac{\dot \xi_n k(1+\varepsilon)(3-k^2)}{(1+k^2)^2}
<2\dot \xi_n < 2S\\ \\
\left |\dot y_{n+1}^{(1)} -\dot y_{n+1}^{(12)}  \right | &=&
\frac{\dot \xi_n(1+\varepsilon)\left| \dot \eta_n ^2(3k^2-1)+ \dot
\xi_n ^2(1+k^2)\right|}{4D(1+k^2)}\\ \\
&\approx& \frac{\dot \xi_n (1+\varepsilon)|3k^2-1|}{(1+k^2)^2} < 2
\dot \xi_n<2S \ .
\end{eqnarray*}

It follows that
\begin{eqnarray*}
&&\left | \dot \xi_{n+1}^{(1)} - \dot \xi_{n+1}^{(12)} \right | <
|k| \left |\dot x_{n+1}^{(1)}  -\dot x_{n+1}^{(12)}  \right | +
\left | \dot y_{n+1}^{(1)}  - \dot y_{n+1}^{(12)} \right | <4 S \\
\\
&&\left | \dot \eta_{n+1}^{(1)} - \dot \eta_{n+1}^{(12)} \right |
< |k| \left |\dot x_{n+1}^{(1)}  -\dot x_{n+1}^{(12)}  \right | +
\left | \dot y_{n+1}^{(1)}  - \dot y_{n+1}^{(12)} \right | <4 S \
.
\end{eqnarray*}

\vskip0.5truecm

\noindent Summing up, we have
\begin{eqnarray*}
&\dot \xi^{(1)}_{n+1} > 0 \qquad \text{and} \qquad  \dot
\eta_{n+1}^{(1)} <0 \ , &  \\ \\
&\dot \xi^{(1)}_{n+1} -4S <\dot \xi^{(12)}_{n+1} < \dot
\xi^{(1)}_{n+1} +4S \qquad \text{and} \qquad \dot \eta_{n+1}^{(1)}
-4S < \dot \eta_{n+1}^{(12)} <0 \ . &
\end{eqnarray*}
Since $\dot \eta_{n+1}^{(1)} , \dot \eta_{n+1}^{(12)}  <0$, if one
of TA and NA does not stop, then it computes a new velocity vector
using the formul\ae~of $\Z_2$.

In general, there are the following  cases.
\begin{enumerate}[label=\alph*)]
\item \underline{Let $0 < \dot \xi_{n+1}^{(1)}  \le S$ and $\dot
\xi_{n+1}^{(12)}  \le 0$}.  NA and TA stop, since $\dot
\eta_{n+1}^{(1)} , \dot \eta _{n+1}^{(12)}   <0$,  and they return
the similar outputs $(\dot x_{n+1}^{(1)} , \ \dot y_{n+1}^{(1)} )$
and  $(\dot x_{n+1}^{(12)} , \ \dot y_{n+1}^{(12)}  )$.
\item \underline{Let $0 < \dot \xi_{n+1}^{(1)}  \le S$ and $\dot
\xi_{n+1}^{(12)} > 0$}.
  NA stops, since $\dot \eta_{n+1}^{(1)}  <0$. \\
TA uses the formul\ae~in $\Z_2$, since $\dot \eta_{n+1}^{(12)}
<0$. We have   $0< \dot \xi_{n+1}^{(12)}  \le \dot \xi_{n+1}^{(1)}
+ 4 S<5 S$ and $\eta_{n+1}^{(12)}  <0$,  and  so
Lemma~\ref{small_changes} implies that TA performs small changes
to  $(\dot x_{n+1}^{(12)}, \dot y_{n+1}^{(12)})$.

Moreover, $\dot \xi_{n+2}= - \varepsilon \dot \xi_{n+1}^{(12)} $
and so $-5S < \dot \xi_{n+2}<0$. If $\dot \eta_{n+2} \le 0$, then
TA stops with a similar output as NA. Otherwise, if $\dot
\eta_{n+2}  >0$, then
$$ 0 < \dot \eta_{n+2}  = \dot \eta_{n+1}^{(12)} + (1+\varepsilon)\frac{1-k^2}{1+k^2} \dot \xi_{n+1} ^{(12)}$$
and so, since $\dot \eta_{n+1}^{(12)}<0$,
$$ 0 < - \dot \eta_{n+1}^{(12)} < (1+\varepsilon)\frac{1-k^2}{1+k^2} \dot \xi_{n+1} ^{(12})< 10 S \ . $$
In this case, since $\dot \xi_{n+1}^{(12)}$ and $\dot
\eta_{n+1}^{(12)}$ are very small, the disk, at the $(n+1)$-th
step,  is almost at rest and the output of TA and NA is very
similar, independently of the number of steps performed by TA,
after NA has stopped.
\item \underline{Let $\xi_{n+1}^{(1)}  > S$ and $\dot
\xi_{n+1}^{(12)}  \le 0$}. TA stops, since $\eta_{n+1}^{(12)} <0$.
NA uses the formul\ae~in $\Z_2$. Since $0>\xi_{n+1}^{(12)}  >
\xi_{n+1}^{(1)}  - 4S$, then $S<\xi_{n+1}^{(1)} < 4S$ and, from
Lemma~\ref{small_changes}, NA performs small changes to  $(\dot
x_{n+1}^{(1)}, \dot y_{n+1}^{(1)})$.

Since  $\dot \xi_{n+2}= - \varepsilon \dot \xi_{n+1}^{(1)} $, we
have $-4 S< \dot \xi_{n+2}<0$ and, if $\dot \eta_{n+2} \le 0$,
then NA stops and it returns a final velocity similar to the
output of TA. Otherwise, if $\dot \eta_{n+2} > 0$, since $\dot
\eta_{n+1}^{(1)} < 0$, we have that
$$ 0 < - \dot \eta_{n+1}^{(1)}  < (1+\varepsilon)\frac{1-k^2}{1+k^2} \dot \xi_{n+1} ^{(1)}< 8 S \ .$$
In this case, since $\dot \xi_{n+1}^{(11)}$ and $\dot
\eta_{n+1}^{(12)}$ are very small, the disk, at the $(n+1)$-th
step,  is almost at rest and the output of TA and NA is very
similar, independently of the number of steps performed by NA,
after TA has stopped.
\item \underline{Let $\dot \xi_{n+1}^{(1)}  > S$ and $\dot
\xi_{n+1}^{(12)}
> 0$}.  Both TA and NA compute, using the same formul\ae~in
$\Z_2$,  a new iteration starting from two similar velocities.
Lemma~\ref{diff} implies that the new computed velocities slightly
differ from each other. Moreover, the new  values of $\dot \xi$,
equal to $-\varepsilon \dot \xi_{n+1}^{(1)} $ and $-\varepsilon
\dot \xi_{n+1}^{(12)} $ respectively,  are negative, so that we
can repeat an analysis of the behavior of TA and NA analogous to
the one presented in items  a -- d, changing the role of $\dot
\xi$ and $\dot \eta$.
\end{enumerate}
\end{enumerate}

Table~\ref{cases} illustrates the possible different cases after
the $n$-th step.
\begin{table}[h]
\begin{tabular}{| c | c | c | c |}
\cline{2-4}
\multicolumn{1}{ c| }{ }&  \multicolumn{1}{ c| }{ $ \phantom{\dfrac12} \dot \xi_{n+1} \le 0  $}  &  \multicolumn{1}{ c| }{ $ 0 < \dot \xi _{n+1}<S  $}&  \multicolumn{1}{ c |}{ $\dot \xi _{n+1}\ge S$}\\
\hline
$\dot \eta_{n+1} \le 0$ &TA and NA  &  TA small changes & TA and NA \\
&stop  & NA stops & in $\Z_{2}$\\
\hline
$0 < \dot \eta_{n+1} < S$ & TA small changes &  TA disk almost at rest & TA in $\Z_{12}$\\
&NA stops  & NA stops &  NA in $\Z_2$ \\
\hline
$\dot \eta_{n+1} \ge S$ & TA and NA & TA in $\Z_{12}$ & TA and NA \\
& in $\Z_{1}$ &NA in $\Z_{1}$ & in $\Z_{12}$\\
\hline
\end{tabular}
\caption{Behavior of TA and NA, starting from the same $(\dot
\xi_n, \dot \eta_n)$. }\label{cases}
\end{table}

In conclusion, the following cases happen:
\begin{enumerate}
\item both algorithms have the same behavior (same iterations and
same steps number); \item NA stops and TA does not stop and it
makes small changes computing the new velocities or vice versa;
\item NA stops and TA does not stop, but the disk is almost at
rest; \item both algorithms computes velocities whose difference
is very small.
\end{enumerate}
\end{Proof}

\subsection{Numerical examples}

We consider several examples, obtained by varying  the coefficient
$\varepsilon$ and the angle $\alpha$.  Furthermore, for each pair
$(\varepsilon, \ \alpha)$,  we consider several initial linear
velocity ${\bf v}_0 = (\dot{x}_0, \ \dot{y}_0 )$, with $\|{\bf
v}_0\|_{_2}=1$. The following Table~\ref{coeff} shows some
possible values of $\varepsilon $ and $\alpha$ and, denoting by $k
= \tan \alpha$, some initial velocities $\widehat{\bf v}_0$ such
that ${\bf v}_0=\dfrac{\widehat {\bf v}_0 }{\|\widehat {\bf v}_0
\|_{_2}}$.
\begin{table}[ht]
 \begin{center}
 \begin{tabular}{|c|| c | c | c | c | c| c| c|}
 \hline
  \multicolumn{8}{|c|}{Values of $ \varepsilon $}  \\
\hline
 \#& $1$ & $2 $ & $3 $ & $ 4 $  & $5$ & $6$ & $7$ \\
\hline
$\varepsilon $ &   $1$ &  $0.95$ &$0.75$ &$0.5$ & $0.25$  &$0.05$ & $0$  \\
 \hline
 \hline
  \multicolumn{8}{|c|}{Values of $ \alpha $}  \\
\hline
 \#& $1$ & $2 $ & $3 $ & $ 4 $  & $5$ & $6$ & $7$ \\
\hline
$\alpha $ &   $\pi /4$ &  $\pi /6$&   $\pi /8$&   $\pi /12$&   $\pi /16$&   $\pi /32$&   $\pi /64$  \\
 \hline
\hline
  \multicolumn{8}{|c|}{Values of $\widehat{\bf v}_0$}  \\
\hline
 \#& $1$ & $2 $ & $3 $ & $ 4 $  & $5$ & $6$ & $7$ \\
\hline
$\widehat{\bf v}_0$ &  $(1, 0)$&     $(1,k/3)$&     $(1, 2k/3)$&   $(1,k)$&   $(1,1/k)$  & $(0,1)$&       $(-1, k)$   \\
 \hline
\end{tabular}
\end{center}
\caption{Coefficients, angles, initial velocity}   \label{coeff}
\end{table}

\vskip3truecm

The following tables collect the results of NA corresponding to
pairs $(\varepsilon, \alpha)$. The entries of each line of a table
show, respectively, the number of the example, the initial
velocity vector ${\bf v}_0$, the final velocity vector ${\bf v}_f$
and its norm,  the case to which ${\bf v}_0$ belongs, the number
$N$ of steps to obtain ${\bf v}_f$ and if the algorithm stops
because ${\bf v}_f \in \Z_0$ or because $\|{\bf v}_f \|_{_2} <
S_v$.

The results are obtained processing the previous data by the  NA,
implemented  in MatLab, using $S=2\cdot 2^{-52}=4.44 \cdot
10^{-16} $, which corresponds to twice  the machine precision, $
S_v= 10^{-12}$, and the maximum number of step $N_{max}= 10^4$.
Note that a similar threshold $S_v$ ensure that, with an input
velocity of $1$ kilometer/second, the rest condition is fixed for
an output velocity of less than $1$ nanometer/second.

Later on we denote with $i.j.k$ the example where  $\varepsilon $
assumes  the $i$-th  value, $\alpha$ assumes the $j$-th value and
$v_0$ the $k$-th value of the Table~\ref{coeff}, e.g. the case
$3.2$ is obtained using $\varepsilon = 0.75$ and $\alpha= \pi /6$
and the case $6.4.2$ is obtained using $\varepsilon = 0.05$,
$\alpha= \pi /12$ and $\widehat{\bf v}_0= (1,k/3)$.

The cases *.*.1 and *.*.7 are test situations: in fact all the
cases *.*.1 are such that ${\bf v}_0$ has the direction of the
angle bisector and the behavior of the algorithm is known from
theoretical results; all the cases *.*.7 are such that ${\bf
v}_0\in \Z_0$ and then ${\bf v}_0$ is immediately the output
velocity. Moreover, the cases *.*.2 and *.*.3 are the only ones
with multiple impact. In the case *.*.4 the disk moves along in
contact with one wall and impacts with the other. In the case
*.*.5 the initial velocity is orthogonal to one wall. In the case
*.*.6 the initial velocity is orthogonal to the bisector of the
angle.

The cases *.1.* too are test situations, since when $\alpha =
\pi/4$ or, that is the same, $k=1$, once again the behavior of the
algorithm is known from theoretical results.

The case 7.*.*, such that $\varepsilon=0$, does not imply that the
disk stops after the first impact, since only that the orthogonal
component of the velocity with respect to the impacted wall is
annihilated.

\begin{Remark}
The case 3.7.*, where $\alpha = \pi/64$, is the first case for
which, in perfect agreement with the theoretical considerations
about the algorithm, the algorithm stops because of the presence
of the threshold $S_v$ on the norm of ${\bf v}_n$. Moreover, the
case 3.*.*,  where $\varepsilon=0.75$, is the first where it
becomes more evident that, once again in perfect agreement with
the theoretical considerations about the algorithm, the number of
steps of the algorithm increases when the angle $\alpha$ between
the walls becomes smaller. Both this events are more highlighted
by the subsequent examples.
\end{Remark}

%\clearpage
\subsection*{ Ideal case:  $\varepsilon=1$. }
\begin{table}[h]
\begin{tabular}{|c||c||c||c||c||c||c|}
 \hline
  \multicolumn{7}{|c|}{ $ \alpha = \pi / 4 =     0.78540 \qquad \Rightarrow \qquad k = 1 $}  \\
\hline
  $\#$ & $(\dot{x}_0, \dot{y}_0 )$ & $ (\dot{x}_f ,  \dot{y}_f) $  &$ \|{\bf v}_f \|$& $\Z({\bf v}_0)$ & $N$ & Stop \\
\hline
 1.1.1  & (      1.000  ,      0.000 )& (-1.00e+00 , -0.00e+00) &    1   & $\Z_{12}$  & 1  & $\Z_0$  \\
  1.1.2  & (      0.949  ,      0.316 )& (-9.49e-01 , -3.16e-01) &    1   & $\Z_{12}$  & 1  & $\Z_0$  \\
  1.1.3  & (      0.832  ,      0.555 )& (-8.32e-01 , -5.55e-01) &    1   & $\Z_{12}$  & 1  & $\Z_0$  \\
  1.1.4  & (      0.707  ,      0.707 )& (-7.07e-01 , -7.07e-01) &    1   & $\Z_{2}$  & 1  & $\Z_0$  \\
  1.1.5  & (      0.707  ,      0.707 )& (-7.07e-01 , -7.07e-01) &    1   & $\Z_{2}$  & 1  & $\Z_0$  \\
  1.1.6  & (      0.000  ,      1.000 )& (-1.00e+00 , -1.11e-16) &  1  & $\Z_{2}$  & 1  & $\Z_0$  \\
  1.1.7  & (     -0.707  ,      0.707 )& (-7.07e-01 , 7.07e-01) & 1  & $\Z_{0}$  & 0  & $\Z_0$  \\
 \hline
\end{tabular}
 \end{table}
 \begin{table}[h]
 \begin{tabular}{|c||c||c||c||c||c||c|}
 \hline
  \multicolumn{7}{|c|}{ $ \alpha = \pi / 6 =     0.52360 \qquad \Rightarrow \qquad k=    0.57735 $}  \\
\hline
  $\#$ & $(\dot{x}_0, \dot{y}_0 )$ & $ (\dot{x}_f ,  \dot{y}_f) $  &$ \|{\bf v}_f \|$& $\Z({\bf v}_0)$ & $N$ & Stop \\
\hline
 1.2.1  & (      1.000  ,      0.000 )& (-1.00e+00 , -0.00e+00) &  1& $\Z_{12}$  & 1  & $\Z_0$  \\
  1.2.2  & (      0.982  ,      0.189 )& (-9.82e-01 , -1.89e-01) &  1& $\Z_{12}$  & 2  & $\Z_0$  \\
  1.2.3  & (      0.933  ,      0.359 )& (-9.55e-01 , 2.95e-01) &  1& $\Z_{12}$  & 2  & $\Z_0$  \\
  1.2.4  & (      0.866  ,      0.500 )& (-8.66e-01 , 5.00e-01) &  1& $\Z_{2}$  & 2  & $\Z_0$  \\
  1.2.5  & (      0.500  ,      0.866 )& (-1.00e+00 , 8.33e-17) &  1& $\Z_{2}$  & 2  & $\Z_0$  \\
  1.2.6  & (      0.000  ,      1.000 )& (-8.66e-01 , -5.00e-01) &  1& $\Z_{2}$  & 1  & $\Z_0$  \\
  1.2.7  & (     -0.866  ,      0.500 )& (-8.66e-01 , 5.00e-01) &  1& $\Z_{0}$  & 0  & $\Z_0$  \\
 \hline
\end{tabular}
 \end{table}
 \begin{table}[h]
 \begin{tabular}{|c||c||c||c||c||c||c|}
 \hline
  \multicolumn{7}{|c|}{ $ \alpha = \pi / 8 =     0.39270 \qquad \Rightarrow \qquad k=    0.41421 $}  \\
\hline
  $\#$ & $(\dot{x}_0, \dot{y}_0 )$ & $ (\dot{x}_f ,  \dot{y}_f) $  &$ \|{\bf v}_f \|$& $\Z({\bf v}_0)$ & $N$ & Stop \\
\hline
 1.3.1  & (      1.000  ,      0.000 )& (-1.00e+00 , -0.00e+00) &  1& $\Z_{12}$  & 1  & $\Z_0$  \\
  1.3.2  & (      0.991  ,      0.137 )& (-9.38e-01 , 3.45e-01) &  1& $\Z_{12}$  & 3  & $\Z_0$  \\
  1.3.3  & (      0.964  ,      0.266 )& (-9.82e-01 , -1.88e-01) &  1& $\Z_{12}$  & 3  & $\Z_0$  \\
  1.3.4  & (      0.924  ,      0.383 )& (-9.24e-01 , -3.83e-01) &  1& $\Z_{2}$  & 3  & $\Z_0$  \\
  1.3.5  & (      0.383  ,      0.924 )& (-9.24e-01 , 3.83e-01) &  1& $\Z_{2}$  & 2  & $\Z_0$  \\
  1.3.6  & (      0.000  ,      1.000 )& (-1.00e+00 , 0.00e+00) &  1& $\Z_{2}$  & 2  & $\Z_0$  \\
  1.3.7  & (     -0.924  ,      0.383 )& (-9.24e-01 , 3.83e-01) &  1& $\Z_{0}$  & 0  & $\Z_0$  \\
 \hline
\end{tabular}
 \end{table}
 \begin{table}[h]
 \begin{tabular}{|c||c||c||c||c||c||c|}
 \hline
  \multicolumn{7}{|c|}{ $ \alpha = \pi / 12 =     0.26180 \qquad \Rightarrow \qquad k=    0.26795 $}  \\
\hline
  $\#$ & $(\dot{x}_0, \dot{y}_0 )$ & $ (\dot{x}_f ,  \dot{y}_f) $  &$ \|{\bf v}_f \|$& $\Z({\bf v}_0)$ & $N$ & Stop \\
\hline
 1.4.1  & (      1.000  ,      0.000 )& (-1.00e+00 , -0.00e+00) &  1& $\Z_{12}$  & 1  & $\Z_0$  \\
  1.4.2  & (      0.996  ,      0.089 )& (-9.92e-01 , 1.27e-01) &  1& $\Z_{12}$  & 4  & $\Z_0$  \\
  1.4.3  & (      0.984  ,      0.176 )& (-9.94e-01 , -1.06e-01) &  1& $\Z_{12}$  & 5  & $\Z_0$  \\
  1.4.4  & (      0.966  ,      0.259 )& (-9.66e-01 , -2.59e-01) &  1& $\Z_{2}$  & 5  & $\Z_0$  \\
  1.4.5  & (      0.259  ,      0.966 )& (-9.66e-01 , -2.59e-01) &  1& $\Z_{2}$  & 3  & $\Z_0$  \\
  1.4.6  & (      0.000  ,      1.000 )& (-1.00e+00 , -1.04e-16) &  1& $\Z_{2}$  & 3  & $\Z_0$  \\
  1.4.7  & (     -0.966  ,      0.259 )& (-9.66e-01 , 2.59e-01) &  1& $\Z_{0}$  & 0  & $\Z_0$  \\
 \hline
\end{tabular}
 \end{table}
 \begin{table}[h]
 \begin{tabular}{|c||c||c||c||c||c||c|}
 \hline
  \multicolumn{7}{|c|}{ $ \alpha = \pi / 16 =     0.19635 \qquad \Rightarrow \qquad k=    0.19891 $}  \\
\hline
  $\#$ & $(\dot{x}_0, \dot{y}_0 )$ & $ (\dot{x}_f ,  \dot{y}_f) $  &$ \|{\bf v}_f \|$& $\Z({\bf v}_0)$ & $N$ & Stop \\
\hline
 1.5.1  & (      1.000  ,      0.000 )& (-1.00e+00 , -0.00e+00) &  1& $\Z_{12}$  & 1  & $\Z_0$  \\
  1.5.2  & (      0.998  ,      0.066 )& (-9.99e-01 , 3.58e-02) &  1& $\Z_{12}$  & 6  & $\Z_0$  \\
  1.5.3  & (      0.991  ,      0.131 )& (-9.97e-01 , -7.36e-02) &  1& $\Z_{12}$  & 7  & $\Z_0$  \\
  1.5.4  & (      0.981  ,      0.195 )& (-9.81e-01 , -1.95e-01) &  1& $\Z_{2}$  & 7  & $\Z_0$  \\
  1.5.5  & (      0.195  ,      0.981 )& (-9.81e-01 , 1.95e-01) &  1& $\Z_{2}$  & 4  & $\Z_0$  \\
  1.5.6  & (      0.000  ,      1.000 )& (-1.00e+00 , 0.00e+00) &  1& $\Z_{2}$  & 4  & $\Z_0$  \\
  1.5.7  & (     -0.981  ,      0.195 )& (-9.81e-01 , 1.95e-01) &  1& $\Z_{0}$  & 0  & $\Z_0$  \\
 \hline
\end{tabular}
 \end{table}
 \begin{table}[h]
 \begin{tabular}{|c||c||c||c||c||c||c|}
 \hline
  \multicolumn{7}{|c|}{ $ \alpha = \pi / 32 =     0.09817 \qquad \Rightarrow \qquad k=    0.09849 $}  \\
\hline
  $\#$ & $(\dot{x}_0, \dot{y}_0 )$ & $ (\dot{x}_f ,  \dot{y}_f) $  &$ \|{\bf v}_f \|$& $\Z({\bf v}_0)$ & $N$ & Stop \\
\hline
 1.6.1  & (      1.000  ,      0.000 )& (-1.00e+00 , -0.00e+00) &  1& $\Z_{12}$  & 1  & $\Z_0$  \\
  1.6.2  & (      0.999  ,      0.033 )& (-1.00e+00 , -1.84e-02) &  1& $\Z_{12}$  & 14  & $\Z_0$  \\
  1.6.3  & (      0.998  ,      0.066 )& (-9.99e-01 , -3.38e-02) &  1& $\Z_{12}$  & 15  & $\Z_0$  \\
  1.6.4  & (      0.995  ,      0.098 )& (-9.95e-01 , -9.80e-02) &  1& $\Z_{2}$  & 15  & $\Z_0$  \\
  1.6.5  & (      0.098  ,      0.995 )& (-9.95e-01 , 9.80e-02) &  1& $\Z_{2}$  & 8  & $\Z_0$  \\
  1.6.6  & (      0.000  ,      1.000 )& (-1.00e+00 , 8.25e-17) &  1& $\Z_{2}$  & 8  & $\Z_0$  \\
  1.6.7  & (     -0.995  ,      0.098 )& (-9.95e-01 , 9.80e-02) &  1& $\Z_{0}$  & 0  & $\Z_0$  \\
 \hline
\end{tabular}
 \end{table}
 \begin{table}[h]
 \begin{tabular}{|c||c||c||c||c||c||c|}
 \hline
  \multicolumn{7}{|c|}{ $ \alpha = \pi / 64 =     0.04909 \qquad \Rightarrow \qquad k=    0.04913 $}  \\
\hline
  $\#$ & $(\dot{x}_0, \dot{y}_0 )$ & $ (\dot{x}_f ,  \dot{y}_f) $  &$ \|{\bf v}_f \|$& $\Z({\bf v}_0)$ & $N$ & Stop \\
\hline
 1.7.1  & (      1.000  ,      0.000 )& (-1.00e+00 , -0.00e+00) &  1& $\Z_{12}$  & 1  & $\Z_0$  \\
  1.7.2  & (      1.000  ,      0.016 )& (-1.00e+00 , -1.45e-02) &  1& $\Z_{12}$  & 30  & $\Z_0$  \\
  1.7.3  & (      0.999  ,      0.033 )& (-1.00e+00 , -1.65e-02) &  1& $\Z_{12}$  & 31  & $\Z_0$  \\
  1.7.4  & (      0.999  ,      0.049 )& (-9.99e-01 , -4.91e-02) &  1& $\Z_{2}$  & 31  & $\Z_0$  \\
  1.7.5  & (      0.049  ,      0.999 )& (-9.99e-01 , 4.91e-02) &  1& $\Z_{2}$  & 16  & $\Z_0$  \\
  1.7.6  & (      0.000  ,      1.000 )& (-1.00e+00 , 4.15e-17) &  1& $\Z_{2}$  & 16  & $\Z_0$  \\
  1.7.7  & (     -0.999  ,      0.049 )& (-9.99e-01 , 4.91e-02) & 1& $\Z_{0}$  & 0  & $\Z_0$  \\
 \hline
\end{tabular}
 \end{table}

 \clearpage
\subsection*{ Non ideal case: $\varepsilon=0.95$. }
\begin{table}[h]
 \begin{tabular}{|c||c||c||c||c||c||c|}
 \hline
  \multicolumn{7}{|c|}{ $ \alpha = \pi / 4 =     0.78540 \qquad \Rightarrow \qquad k=   1  $}  \\
\hline
  $\#$ & $(\dot{x}_0, \dot{y}_0 )$ & $ (\dot{x}_f ,  \dot{y}_f) $  &$ \|{\bf v}_f \|$& $\Z({\bf v}_0)$ & $N$ & Stop \\
\hline
 2.1.1  & (      1.000  ,      0.000 )& (-9.50e-01 , -0.00e+00) & 9.50e-01& $\Z_{12}$  & 1  & $\Z_0$  \\
  2.1.2  & (      0.949  ,      0.316 )& (-9.01e-01 , -3.00e-01) & 9.50e-01& $\Z_{12}$  & 1  & $\Z_0$  \\
  2.1.3  & (      0.832  ,      0.555 )& (-7.90e-01 , -5.27e-01) & 9.50e-01& $\Z_{12}$  & 1  & $\Z_0$  \\
  2.1.4  & (      0.707  ,      0.707 )& (-6.72e-01 , -6.72e-01) & 9.50e-01& $\Z_{2}$  & 1  & $\Z_0$  \\
  2.1.5  & (      0.707  ,      0.707 )& (-6.72e-01 , -6.72e-01) & 9.50e-01& $\Z_{2}$  & 1  & $\Z_0$  \\
  2.1.6  & (      0.000  ,      1.000 )& (-9.75e-01 , 2.50e-02) & 9.75e-01& $\Z_{2}$  & 1  & $\Z_0$  \\
  2.1.7  & (     -0.707  ,      0.707 )& (-7.07e-01 , 7.07e-01) & 1.00e+00& $\Z_{0}$  & 0  & $\Z_0$  \\
 \hline
\end{tabular}
 \end{table}
 \begin{table}[h]
 \begin{tabular}{|c||c||c||c||c||c||c|}
 \hline
  \multicolumn{7}{|c|}{ $ \alpha = \pi / 6 =     0.52360 \qquad \Rightarrow \qquad k=    0.57735 $}  \\
\hline
  $\#$ & $(\dot{x}_0, \dot{y}_0 )$ & $ (\dot{x}_f ,  \dot{y}_f) $  &$ \|{\bf v}_f \|$& $\Z({\bf v}_0)$ & $N$ & Stop \\
\hline
 2.2.1  & (      1.000  ,      0.000 )& (-9.50e-01 , -0.00e+00) & 9.50e-01& $\Z_{12}$  & 1  & $\Z_0$  \\
  2.2.2  & (      0.982  ,      0.189 )& (-9.33e-01 , -1.80e-01) & 9.50e-01& $\Z_{12}$  & 2  & $\Z_0$  \\
  2.2.3  & (      0.933  ,      0.359 )& (-8.94e-01 , 2.72e-01) & 9.35e-01& $\Z_{12}$  & 2  & $\Z_0$  \\
  2.2.4  & (      0.866  ,      0.500 )& (-8.02e-01 , 4.63e-01) & 9.26e-01& $\Z_{2}$  & 3  & $\Z_0$  \\
  2.2.5  & (      0.500  ,      0.866 )& (-9.38e-01 , -2.06e-02) & 9.38e-01& $\Z_{2}$  & 2  & $\Z_0$  \\
  2.2.6  & (      0.000  ,      1.000 )& (-8.44e-01 , -4.63e-01) & 9.63e-01& $\Z_{2}$  & 1  & $\Z_0$  \\
  2.2.7  & (     -0.866  ,      0.500 )& (-8.66e-01 , 5.00e-01) & 1.00e+00& $\Z_{0}$  & 0  & $\Z_0$  \\
 \hline
\end{tabular}
 \end{table}
 \begin{table}[h]
 \begin{tabular}{|c||c||c||c||c||c||c|}
 \hline
  \multicolumn{7}{|c|}{ $ \alpha = \pi / 8 =     0.39270 \qquad \Rightarrow \qquad k=    0.41421 $}  \\
\hline
  $\#$ & $(\dot{x}_0, \dot{y}_0 )$ & $ (\dot{x}_f ,  \dot{y}_f) $  &$ \|{\bf v}_f \|$& $\Z({\bf v}_0)$ & $N$ & Stop \\
\hline
 2.3.1  & (      1.000  ,      0.000 )& (-9.50e-01 , -0.00e+00) & 9.50e-01& $\Z_{12}$  & 1  & $\Z_0$  \\
  2.3.2  & (      0.991  ,      0.137 )& (-8.77e-01 , 3.28e-01) & 9.37e-01& $\Z_{12}$  & 3  & $\Z_0$  \\
  2.3.3  & (      0.964  ,      0.266 )& (-8.97e-01 , -1.65e-01) & 9.12e-01& $\Z_{12}$  & 3  & $\Z_0$  \\
  2.3.4  & (      0.924  ,      0.383 )& (-8.34e-01 , -3.44e-01) & 9.02e-01& $\Z_{2}$  & 4  & $\Z_0$  \\
  2.3.5  & (      0.383  ,      0.924 )& (-8.65e-01 , 3.33e-01) & 9.27e-01& $\Z_{2}$  & 2  & $\Z_0$  \\
  2.3.6  & (      0.000  ,      1.000 )& (-9.51e-01 , -3.38e-02) & 9.51e-01& $\Z_{2}$  & 2  & $\Z_0$  \\
  2.3.7  & (     -0.924  ,      0.383 )& (-9.24e-01 , 3.83e-01) & 1.00e+00& $\Z_{0}$  & 0  & $\Z_0$  \\
 \hline
\end{tabular}
 \end{table}
 \begin{table}[h]
 \begin{tabular}{|c||c||c||c||c||c||c|}
 \hline
  \multicolumn{7}{|c|}{ $ \alpha = \pi / 12 =     0.26180 \qquad \Rightarrow \qquad k=    0.26795 $}  \\
\hline
  $\#$ & $(\dot{x}_0, \dot{y}_0 )$ & $ (\dot{x}_f ,  \dot{y}_f) $  &$ \|{\bf v}_f \|$& $\Z({\bf v}_0)$ & $N$ & Stop \\
\hline
 2.4.1  & (      1.000  ,      0.000 )& (-9.50e-01 , -0.00e+00) & 9.50e-01& $\Z_{12}$  & 1  & $\Z_0$  \\
  2.4.2  & (      0.996  ,      0.089 )& (-8.87e-01 , 9.37e-02) & 8.92e-01& $\Z_{12}$  & 4  & $\Z_0$  \\
  2.4.3  & (      0.984  ,      0.176 )& (-8.59e-01 , -8.55e-02) & 8.63e-01& $\Z_{12}$  & 5  & $\Z_0$  \\
  2.4.4  & (      0.966  ,      0.259 )& (-8.29e-01 , -2.19e-01) & 8.57e-01& $\Z_{2}$  & 6  & $\Z_0$  \\
  2.4.5  & (      0.259  ,      0.966 )& (-8.83e-01 , -1.96e-01) & 9.04e-01& $\Z_{2}$  & 3  & $\Z_0$  \\
  2.4.6  & (      0.000  ,      1.000 )& (-9.27e-01 , 4.60e-02) & 9.28e-01& $\Z_{2}$  & 3  & $\Z_0$  \\
  2.4.7  & (     -0.966  ,      0.259 )& (-9.66e-01 , 2.59e-01) & 1.00e+00& $\Z_{0}$  & 0  & $\Z_0$  \\
 \hline
\end{tabular}
 \end{table}
 \begin{table}[h]
 \begin{tabular}{|c||c||c||c||c||c||c|}
 \hline
  \multicolumn{7}{|c|}{ $ \alpha = \pi / 16 =     0.19635 \qquad \Rightarrow \qquad k=    0.19891 $}  \\
\hline
  $\#$ & $(\dot{x}_0, \dot{y}_0 )$ & $ (\dot{x}_f ,  \dot{y}_f) $  &$ \|{\bf v}_f \|$& $\Z({\bf v}_0)$ & $N$ & Stop \\
\hline
 2.5.1  & (      1.000  ,      0.000 )& (-9.50e-01 , -0.00e+00) & 9.50e-01& $\Z_{12}$  & 1  & $\Z_0$  \\
  2.5.2  & (      0.998  ,      0.066 )& (-8.42e-01 , 7.41e-03) & 8.42e-01& $\Z_{12}$  & 6  & $\Z_0$  \\
  2.5.3  & (      0.991  ,      0.131 )& (-8.16e-01 , -5.77e-02) & 8.18e-01& $\Z_{12}$  & 7  & $\Z_0$  \\
  2.5.4  & (      0.981  ,      0.195 )& (-8.00e-01 , -1.54e-01) & 8.14e-01& $\Z_{2}$  & 8  & $\Z_0$  \\
  2.5.5  & (      0.195  ,      0.981 )& (-8.75e-01 , 1.20e-01) & 8.83e-01& $\Z_{2}$  & 4  & $\Z_0$  \\
  2.5.6  & (      0.000  ,      1.000 )& (-9.05e-01 , -5.77e-02) & 9.06e-01& $\Z_{2}$  & 4  & $\Z_0$  \\
  2.5.7  & (     -0.981  ,      0.195 )& (-9.81e-01 , 1.95e-01) & 1.00e+00& $\Z_{0}$  & 0  & $\Z_0$  \\
 \hline
\end{tabular}
 \end{table}
 \begin{table}[h]
 \begin{tabular}{|c||c||c||c||c||c||c|}
 \hline
  \multicolumn{7}{|c|}{ $ \alpha = \pi / 32 =     0.09817 \qquad \Rightarrow \qquad k=    0.09849 $}  \\
\hline
  $\#$ & $(\dot{x}_0, \dot{y}_0 )$ & $ (\dot{x}_f ,  \dot{y}_f) $  &$ \|{\bf v}_f \|$& $\Z({\bf v}_0)$ & $N$ & Stop \\
\hline
 2.6.1  & (      1.000  ,      0.000 )& (-9.50e-01 , -0.00e+00) & 9.50e-01& $\Z_{12}$  & 1  & $\Z_0$  \\
  2.6.2  & (      0.999  ,      0.033 )& (-6.70e-01 , -1.43e-02) & 6.70e-01& $\Z_{12}$  & 14  & $\Z_0$  \\
  2.6.3  & (      0.998  ,      0.066 )& (-6.61e-01 , -3.65e-02) & 6.62e-01& $\Z_{12}$  & 15  & $\Z_0$  \\
  2.6.4  & (      0.995  ,      0.098 )& (-6.59e-01 , -4.75e-02) & 6.61e-01& $\Z_{2}$  & 16  & $\Z_0$  \\
  2.6.5  & (      0.098  ,      0.995 )& (-8.06e-01 , -1.62e-02) & 8.06e-01& $\Z_{2}$  & 8  & $\Z_0$  \\
  2.6.6  & (      0.000  ,      1.000 )& (-8.25e-01 , -6.41e-02) & 8.27e-01& $\Z_{2}$  & 7  & $\Z_0$  \\
  2.6.7  & (     -0.995  ,      0.098 )& (-9.95e-01 , 9.80e-02) & 1.00e+00& $\Z_{0}$  & 0  & $\Z_0$  \\
 \hline
\end{tabular}
 \end{table}
 \begin{table}[h]
 \begin{tabular}{|c||c||c||c||c||c||c|}
 \hline
  \multicolumn{7}{|c|}{ $ \alpha = \pi / 64 =     0.04909 \qquad \Rightarrow \qquad k=    0.04913 $}  \\
\hline
  $\#$ & $(\dot{x}_0, \dot{y}_0 )$ & $ (\dot{x}_f ,  \dot{y}_f) $  &$ \|{\bf v}_f \|$& $\Z({\bf v}_0)$ & $N$ & Stop \\
\hline
 2.7.1  & (      1.000  ,      0.000 )& (-9.50e-01 , -0.00e+00) & 9.50e-01& $\Z_{12}$  & 1  & $\Z_0$  \\
  2.7.2  & (      1.000  ,      0.016 )& (-4.29e-01 , 7.51e-03) & 4.29e-01& $\Z_{12}$  & 31  & $\Z_0$  \\
  2.7.3  & (      0.999  ,      0.033 )& (-4.27e-01 , 1.20e-02) & 4.28e-01& $\Z_{12}$  & 32  & $\Z_0$  \\
  2.7.4  & (      0.999  ,      0.049 )& (-4.27e-01 , 1.49e-02) & 4.27e-01& $\Z_{2}$  & 33  & $\Z_0$  \\
  2.7.5  & (      0.049  ,      0.999 )& (-6.84e-01 , 1.95e-02) & 6.84e-01& $\Z_{2}$  & 14  & $\Z_0$  \\
  2.7.6  & (      0.000  ,      1.000 )& (-7.02e-01 , -1.50e-02) & 7.02e-01& $\Z_{2}$  & 14  & $\Z_0$  \\
  2.7.7  & (     -0.999  ,      0.049 )& (-9.99e-01 , 4.91e-02) & 1.00e+00& $\Z_{0}$  & 0  & $\Z_0$  \\
 \hline
\end{tabular}
 \end{table}
 \clearpage
 \subsection*{ Non ideal case: $\varepsilon=0.75$. }
\begin{table}[h]
 \begin{tabular}{|c||c||c||c||c||c||c|}
 \hline
  \multicolumn{7}{|c|}{ $ \alpha = \pi / 4 =     0.78540 \qquad \Rightarrow \qquad k = 1 $}  \\
\hline
  $\#$ & $(\dot{x}_0, \dot{y}_0 )$ & $ (\dot{x}_f ,  \dot{y}_f) $  &$ \|{\bf v}_f \|$& $\Z({\bf v}_0)$ & $N$ & Stop \\
\hline
 3.1.1  & (      1.000  ,      0.000 )& (-7.50e-01 , -0.00e+00) & 7.50e-01& $\Z_{12}$  & 1  & $\Z_0$  \\
  3.1.2  & (      0.949  ,      0.316 )& (-7.12e-01 , -2.37e-01) & 7.50e-01& $\Z_{12}$  & 1  & $\Z_0$  \\
  3.1.3  & (      0.832  ,      0.555 )& (-6.24e-01 , -4.16e-01) & 7.50e-01& $\Z_{12}$  & 1  & $\Z_0$  \\
  3.1.4  & (      0.707  ,      0.707 )& (-5.30e-01 , -5.30e-01) & 7.50e-01& $\Z_{2}$  & 1  & $\Z_0$  \\
  3.1.5  & (      0.707  ,      0.707 )& (-5.30e-01 , -5.30e-01) & 7.50e-01& $\Z_{2}$  & 1  & $\Z_0$  \\
  3.1.6  & (      0.000  ,      1.000 )& (-8.75e-01 , 1.25e-01) & 8.84e-01& $\Z_{2}$  & 1  & $\Z_0$  \\
  3.1.7  & (     -0.707  ,      0.707 )& (-7.07e-01 , 7.07e-01) & 1.00e+00& $\Z_{0}$  & 0  & $\Z_0$  \\
 \hline
\end{tabular}
 \end{table}
 \begin{table}[h]
 \begin{tabular}{|c||c||c||c||c||c||c|}
 \hline
  \multicolumn{7}{|c|}{ $ \alpha = \pi / 6 =     0.52360 \qquad \Rightarrow \qquad k=    0.57735 $}  \\
\hline
  $\#$ & $(\dot{x}_0, \dot{y}_0 )$ & $ (\dot{x}_f ,  \dot{y}_f) $  &$ \|{\bf v}_f \|$& $\Z({\bf v}_0)$ & $N$ & Stop \\
\hline
 3.2.1  & (      1.000  ,      0.000 )& (-7.50e-01 , -0.00e+00) & 7.50e-01& $\Z_{12}$  & 1  & $\Z_0$  \\
  3.2.2  & (      0.982  ,      0.189 )& (-7.36e-01 , -1.42e-01) & 7.50e-01& $\Z_{12}$  & 2  & $\Z_0$  \\
  3.2.3  & (      0.933  ,      0.359 )& (-6.56e-01 , 1.94e-01) & 6.84e-01& $\Z_{12}$  & 2  & $\Z_0$  \\
  3.2.4  & (      0.866  ,      0.500 )& (-5.67e-01 , 3.15e-01) & 6.49e-01& $\Z_{2}$  & 3  & $\Z_0$  \\
  3.2.5  & (      0.500  ,      0.866 )& (-7.03e-01 , -8.12e-02) & 7.08e-01& $\Z_{2}$  & 2  & $\Z_0$  \\
  3.2.6  & (      0.000  ,      1.000 )& (-7.58e-01 , -3.13e-01) & 8.20e-01& $\Z_{2}$  & 1  & $\Z_0$  \\
  3.2.7  & (     -0.866  ,      0.500 )& (-8.66e-01 , 5.00e-01) & 1.00e+00& $\Z_{0}$  & 0  & $\Z_0$  \\
 \hline
\end{tabular}
 \end{table}
 \begin{table}[h]
 \begin{tabular}{|c||c||c||c||c||c||c|}
 \hline
  \multicolumn{7}{|c|}{ $ \alpha = \pi / 8 =     0.39270 \qquad \Rightarrow \qquad k=    0.41421 $}  \\
\hline
  $\#$ & $(\dot{x}_0, \dot{y}_0 )$ & $ (\dot{x}_f ,  \dot{y}_f) $  &$ \|{\bf v}_f \|$& $\Z({\bf v}_0)$ & $N$ & Stop \\
\hline
 3.3.1  & (      1.000  ,      0.000 )& (-7.50e-01 , -0.00e+00) & 7.50e-01& $\Z_{12}$  & 1  & $\Z_0$  \\
  3.3.2  & (      0.991  ,      0.137 )& (-6.46e-01 , 2.48e-01) & 6.92e-01& $\Z_{12}$  & 3  & $\Z_0$  \\
  3.3.3  & (      0.964  ,      0.266 )& (-5.84e-01 , -1.04e-01) & 5.93e-01& $\Z_{12}$  & 3  & $\Z_0$  \\
  3.3.4  & (      0.924  ,      0.383 )& (-5.24e-01 , -1.95e-01) & 5.59e-01& $\Z_{2}$  & 4  & $\Z_0$  \\
  3.3.5  & (      0.383  ,      0.924 )& (-6.42e-01 , 1.65e-01) & 6.63e-01& $\Z_{2}$  & 2  & $\Z_0$  \\
  3.3.6  & (      0.000  ,      1.000 )& (-7.66e-01 , -1.39e-01) & 7.78e-01& $\Z_{2}$  & 2  & $\Z_0$  \\
  3.3.7  & (     -0.924  ,      0.383 )& (-9.24e-01 , 3.83e-01) & 1.00e+00& $\Z_{0}$  & 0  & $\Z_0$  \\
 \hline
\end{tabular}
 \end{table}
 \begin{table}[h]
 \begin{tabular}{|c||c||c||c||c||c||c|}
 \hline
  \multicolumn{7}{|c|}{ $ \alpha = \pi / 12 =     0.26180 \qquad \Rightarrow \qquad k=    0.26795 $}  \\
\hline
  $\#$ & $(\dot{x}_0, \dot{y}_0 )$ & $ (\dot{x}_f ,  \dot{y}_f) $  &$ \|{\bf v}_f \|$& $\Z({\bf v}_0)$ & $N$ & Stop \\
\hline
 3.4.1  & (      1.000  ,      0.000 )& (-7.50e-01 , -0.00e+00) & 7.50e-01& $\Z_{12}$  & 1  & $\Z_0$  \\
  3.4.2  & (      0.996  ,      0.089 )& (-5.26e-01 , 2.44e-02) & 5.26e-01& $\Z_{12}$  & 4  & $\Z_0$  \\
  3.4.3  & (      0.984  ,      0.176 )& (-4.22e-01 , -6.60e-02) & 4.27e-01& $\Z_{12}$  & 5  & $\Z_0$  \\
  3.4.4  & (      0.966  ,      0.259 )& (-4.04e-01 , -6.57e-02) & 4.10e-01& $\Z_{2}$  & 6  & $\Z_0$  \\
  3.4.5  & (      0.259  ,      0.966 )& (-5.94e-01 , -2.27e-02) & 5.95e-01& $\Z_{2}$  & 3  & $\Z_0$  \\
  3.4.6  & (      0.000  ,      1.000 )& (-6.77e-01 , 1.58e-01) & 6.95e-01& $\Z_{2}$  & 3  & $\Z_0$  \\
  3.4.7  & (     -0.966  ,      0.259 )& (-9.66e-01 , 2.59e-01) & 1.00e+00& $\Z_{0}$  & 0  & $\Z_0$  \\
 \hline
\end{tabular}
 \end{table}
 \begin{table}[h]
 \begin{tabular}{|c||c||c||c||c||c||c|}
 \hline
  \multicolumn{7}{|c|}{ $ \alpha = \pi / 16 =     0.19635 \qquad \Rightarrow \qquad k=    0.19891 $}  \\
\hline
  $\#$ & $(\dot{x}_0, \dot{y}_0 )$ & $ (\dot{x}_f ,  \dot{y}_f) $  &$ \|{\bf v}_f \|$& $\Z({\bf v}_0)$ & $N$ & Stop \\
\hline
 3.5.1  & (      1.000  ,      0.000 )& (-7.50e-01 , -0.00e+00) & 7.50e-01& $\Z_{12}$  & 1  & $\Z_0$  \\
  3.5.2  & (      0.998  ,      0.066 )& (-3.69e-01 , -1.44e-03) & 3.69e-01& $\Z_{12}$  & 6  & $\Z_0$  \\
  3.5.3  & (      0.991  ,      0.131 )& (-2.97e-01 , -4.66e-02) & 3.01e-01& $\Z_{12}$  & 8  & $\Z_0$  \\
  3.5.4  & (      0.981  ,      0.195 )& (-2.92e-01 , 3.11e-03) & 2.92e-01& $\Z_{2}$  & 8  & $\Z_0$  \\
  3.5.5  & (      0.195  ,      0.981 )& (-5.36e-01 , -4.63e-02) & 5.38e-01& $\Z_{2}$  & 4  & $\Z_0$  \\
  3.5.6  & (      0.000  ,      1.000 )& (-6.24e-01 , -6.13e-02) & 6.27e-01& $\Z_{2}$  & 3  & $\Z_0$  \\
  3.5.7  & (     -0.981  ,      0.195 )& (-9.81e-01 , 1.95e-01) & 1.00e+00& $\Z_{0}$  & 0  & $\Z_0$  \\
 \hline
\end{tabular}
 \end{table}
 \begin{table}[h]
 \begin{tabular}{|c||c||c||c||c||c||c|}
 \hline
  \multicolumn{7}{|c|}{ $ \alpha = \pi / 32 =     0.09817 \qquad \Rightarrow \qquad k=    0.09849 $}  \\
\hline
  $\#$ & $(\dot{x}_0, \dot{y}_0 )$ & $ (\dot{x}_f ,  \dot{y}_f) $  &$ \|{\bf v}_f \|$& $\Z({\bf v}_0)$ & $N$ & Stop \\
\hline
 3.6.1  & (      1.000  ,      0.000 )& (-7.50e-01 , -0.00e+00) & 7.50e-01& $\Z_{12}$  & 1  & $\Z_0$  \\
  3.6.2  & (      0.999  ,      0.033 )& (-3.93e-02 , -1.60e-03) & 3.94e-02& $\Z_{12}$  & 20  & $\Z_0$  \\
  3.6.3  & (      0.998  ,      0.066 )& (-3.51e-02 , 2.16e-03) & 3.52e-02& $\Z_{12}$  & 22  & $\Z_0$  \\
  3.6.4  & (      0.995  ,      0.098 )& (-3.48e-02 , 1.36e-03) & 3.49e-02& $\Z_{2}$  & 23  & $\Z_0$  \\
  3.6.5  & (      0.098  ,      0.995 )& (-3.85e-01 , 7.72e-03) & 3.85e-01& $\Z_{2}$  & 6  & $\Z_0$  \\
  3.6.6  & (      0.000  ,      1.000 )& (-4.48e-01 , -3.63e-02) & 4.49e-01& $\Z_{2}$  & 6  & $\Z_0$  \\
  3.6.7  & (     -0.995  ,      0.098 )& (-9.95e-01 , 9.80e-02) & 1.00e+00& $\Z_{0}$  & 0  & $\Z_0$  \\
 \hline
\end{tabular}
 \end{table}
 \begin{table}[h]
 \begin{tabular}{|c||c||c||c||c||c||c|}
 \hline
  \multicolumn{7}{|c|}{ $ \alpha = \pi / 64 =     0.04909 \qquad \Rightarrow \qquad k=    0.04913 $}  \\
\hline
  $\#$ & $(\dot{x}_0, \dot{y}_0 )$ & $ (\dot{x}_f ,  \dot{y}_f) $  &$ \|{\bf v}_f \|$& $\Z({\bf v}_0)$ & $N$ & Stop \\
\hline
 3.7.1  & (      1.000  ,      0.000 )& (-7.50e-01 , -0.00e+00) & 7.50e-01& $\Z_{12}$  & 1  & $\Z_0$  \\
  3.7.2  & (      1.000  ,      0.016 )& (9.00e-13 , 3.57e-13) & 9.68e-13& $\Z_{12}$  & 712  & $\|  \cdot \| $  \\
  3.7.3  & (      0.999  ,      0.033 )& (9.02e-13 , 3.58e-13) & 9.71e-13& $\Z_{12}$  & 714  & $\|  \cdot \| $  \\
  3.7.4  & (      0.999  ,      0.049 )& (9.18e-13 , 3.65e-13) & 9.88e-13& $\Z_{2}$  & 714  & $\|  \cdot \| $  \\
  3.7.5  & (      0.049  ,      0.999 )& (-2.41e-01 , -8.67e-03) & 2.41e-01& $\Z_{2}$  & 9  & $\Z_0$  \\
  3.7.6  & (      0.000  ,      1.000 )& (-2.82e-01 , 5.15e-03) & 2.82e-01& $\Z_{2}$  & 9  & $\Z_0$  \\
  3.7.7  & (     -0.999  ,      0.049 )& (-9.99e-01 , 4.91e-02) & 1.00e+00& $\Z_{0}$  & 0  & $\Z_0$  \\
 \hline
\end{tabular}
 \end{table}

 \clearpage
\subsection*{ Non ideal case: $\varepsilon=0.5$. }
\begin{table}[h]
 \begin{tabular}{|c||c||c||c||c||c||c|}
 \hline
  \multicolumn{7}{|c|}{ $ \alpha = \pi / 4 =     0.78540 \qquad \Rightarrow \qquad k = 1 $}  \\
\hline
  $\#$ & $(\dot{x}_0, \dot{y}_0 )$ & $ (\dot{x}_f ,  \dot{y}_f) $  &$ \|{\bf v}_f \|$& $\Z({\bf v}_0)$ & $N$ & Stop \\
\hline
 4.1.1  & (      1.000  ,      0.000 )& (-5.00e-01 , -0.00e+00) & 5.00e-01& $\Z_{12}$  & 1  & $\Z_0$  \\
  4.1.2  & (      0.949  ,      0.316 )& (-4.74e-01 , -1.58e-01) & 5.00e-01& $\Z_{12}$  & 1  & $\Z_0$  \\
  4.1.3  & (      0.832  ,      0.555 )& (-4.16e-01 , -2.77e-01) & 5.00e-01& $\Z_{12}$  & 1  & $\Z_0$  \\
  4.1.4  & (      0.707  ,      0.707 )& (-3.54e-01 , -3.54e-01) & 5.00e-01& $\Z_{2}$  & 1  & $\Z_0$  \\
  4.1.5  & (      0.707  ,      0.707 )& (-3.54e-01 , -3.54e-01) & 5.00e-01& $\Z_{2}$  & 1  & $\Z_0$  \\
  4.1.6  & (      0.000  ,      1.000 )& (-7.50e-01 , 2.50e-01) & 7.91e-01& $\Z_{2}$  & 1  & $\Z_0$  \\
  4.1.7  & (     -0.707  ,      0.707 )& (-7.07e-01 , 7.07e-01) & 1.00e+00& $\Z_{0}$  & 0  & $\Z_0$  \\
 \hline
\end{tabular}
 \end{table}
 \begin{table}[h]
 \begin{tabular}{|c||c||c||c||c||c||c|}
 \hline
  \multicolumn{7}{|c|}{ $ \alpha = \pi / 6 =     0.52360 \qquad \Rightarrow \qquad k=    0.57735 $}  \\
\hline
  $\#$ & $(\dot{x}_0, \dot{y}_0 )$ & $ (\dot{x}_f ,  \dot{y}_f) $  &$ \|{\bf v}_f \|$& $\Z({\bf v}_0)$ & $N$ & Stop \\
\hline
 4.2.1  & (      1.000  ,      0.000 )& (-5.00e-01 , -0.00e+00) & 5.00e-01& $\Z_{12}$  & 1  & $\Z_0$  \\
  4.2.2  & (      0.982  ,      0.189 )& (-4.91e-01 , -9.45e-02) & 5.00e-01& $\Z_{12}$  & 2  & $\Z_0$  \\
  4.2.3  & (      0.933  ,      0.359 )& (-3.75e-01 , 1.23e-01) & 3.95e-01& $\Z_{12}$  & 2  & $\Z_0$  \\
  4.2.4  & (      0.866  ,      0.500 )& (-3.11e-01 , 1.48e-01) & 3.45e-01& $\Z_{2}$  & 3  & $\Z_0$  \\
  4.2.5  & (      0.500  ,      0.866 )& (-4.37e-01 , -1.08e-01) & 4.51e-01& $\Z_{2}$  & 2  & $\Z_0$  \\
  4.2.6  & (      0.000  ,      1.000 )& (-6.50e-01 , -1.25e-01) & 6.61e-01& $\Z_{2}$  & 1  & $\Z_0$  \\
  4.2.7  & (     -0.866  ,      0.500 )& (-8.66e-01 , 5.00e-01) & 1.00e+00& $\Z_{0}$  & 0  & $\Z_0$  \\
 \hline
\end{tabular}
 \end{table}
 \begin{table}[h]
 \begin{tabular}{|c||c||c||c||c||c||c|}
 \hline
  \multicolumn{7}{|c|}{ $ \alpha = \pi / 8 =     0.39270 \qquad \Rightarrow \qquad k=    0.41421 $}  \\
\hline
  $\#$ & $(\dot{x}_0, \dot{y}_0 )$ & $ (\dot{x}_f ,  \dot{y}_f) $  &$ \|{\bf v}_f \|$& $\Z({\bf v}_0)$ & $N$ & Stop \\
\hline
 4.3.1  & (      1.000  ,      0.000 )& (-5.00e-01 , -0.00e+00) & 5.00e-01& $\Z_{12}$  & 1  & $\Z_0$  \\
  4.3.2  & (      0.991  ,      0.137 )& (-3.80e-01 , 1.42e-01) & 4.06e-01& $\Z_{12}$  & 3  & $\Z_0$  \\
  4.3.3  & (      0.964  ,      0.266 )& (-2.52e-01 , -7.92e-02) & 2.64e-01& $\Z_{12}$  & 3  & $\Z_0$  \\
  4.3.4  & (      0.924  ,      0.383 )& (-2.20e-01 , -4.04e-02) & 2.24e-01& $\Z_{2}$  & 4  & $\Z_0$  \\
  4.3.5  & (      0.383  ,      0.924 )& (-3.94e-01 , 2.80e-02) & 3.95e-01& $\Z_{2}$  & 2  & $\Z_0$  \\
  4.3.6  & (      0.000  ,      1.000 )& (-5.62e-01 , -2.03e-01) & 5.98e-01& $\Z_{2}$  & 2  & $\Z_0$  \\
  4.3.7  & (     -0.924  ,      0.383 )& (-9.24e-01 , 3.83e-01) & 1.00e+00& $\Z_{0}$  & 0  & $\Z_0$  \\
 \hline
\end{tabular}
 \end{table}
 \begin{table}[h]
 \begin{tabular}{|c||c||c||c||c||c||c|}
 \hline
  \multicolumn{7}{|c|}{ $ \alpha = \pi / 12 =     0.26180 \qquad \Rightarrow \qquad k=    0.26795 $}  \\
\hline
  $\#$ & $(\dot{x}_0, \dot{y}_0 )$ & $ (\dot{x}_f ,  \dot{y}_f) $  &$ \|{\bf v}_f \|$& $\Z({\bf v}_0)$ & $N$ & Stop \\
\hline
 4.4.1  & (      1.000  ,      0.000 )& (-5.00e-01 , -0.00e+00) & 5.00e-01& $\Z_{12}$  & 1  & $\Z_0$  \\
  4.4.2  & (      0.996  ,      0.089 )& (-1.77e-01 , 2.70e-02) & 1.79e-01& $\Z_{12}$  & 4  & $\Z_0$  \\
  4.4.3  & (      0.984  ,      0.176 )& (-7.68e-02 , 1.96e-02) & 7.93e-02& $\Z_{12}$  & 7  & $\Z_0$  \\
  4.4.4  & (      0.966  ,      0.259 )& (-6.83e-02 , -5.65e-03) & 6.85e-02& $\Z_{2}$  & 7  & $\Z_0$  \\
  4.4.5  & (      0.259  ,      0.966 )& (-3.22e-01 , 5.39e-02) & 3.26e-01& $\Z_{2}$  & 3  & $\Z_0$  \\
  4.4.6  & (      0.000  ,      1.000 )& (-4.87e-01 , 1.90e-02) & 4.88e-01& $\Z_{2}$  & 2  & $\Z_0$  \\
  4.4.7  & (     -0.966  ,      0.259 )& (-9.66e-01 , 2.59e-01) & 1.00e+00& $\Z_{0}$  & 0  & $\Z_0$  \\
 \hline
\end{tabular}
 \end{table}
 \begin{table}[h]
 \begin{tabular}{|c||c||c||c||c||c||c|}
 \hline
  \multicolumn{7}{|c|}{ $ \alpha = \pi / 16 =     0.19635 \qquad \Rightarrow \qquad k=    0.19891 $}  \\
\hline
  $\#$ & $(\dot{x}_0, \dot{y}_0 )$ & $ (\dot{x}_f ,  \dot{y}_f) $  &$ \|{\bf v}_f \|$& $\Z({\bf v}_0)$ & $N$ & Stop \\
\hline
 4.5.1  & (      1.000  ,      0.000 )& (-5.00e-01 , -0.00e+00) & 5.00e-01& $\Z_{12}$  & 1  & $\Z_0$  \\
  4.5.2  & (      0.998  ,      0.066 )& (-2.37e-02 , 3.83e-03) & 2.40e-02& $\Z_{12}$  & 9  & $\Z_0$  \\
  4.5.3  & (      0.991  ,      0.131 )& (-4.87e-03 , 4.16e-04) & 4.89e-03& $\Z_{12}$  & 14  & $\Z_0$  \\
  4.5.4  & (      0.981  ,      0.195 )& (-4.41e-03 , -1.20e-04) & 4.41e-03& $\Z_{2}$  & 15  & $\Z_0$  \\
  4.5.5  & (      0.195  ,      0.981 )& (-2.74e-01 , -1.70e-02) & 2.74e-01& $\Z_{2}$  & 3  & $\Z_0$  \\
  4.5.6  & (      0.000  ,      1.000 )& (-4.08e-01 , 6.38e-02) & 4.13e-01& $\Z_{2}$  & 3  & $\Z_0$  \\
  4.5.7  & (     -0.981  ,      0.195 )& (-9.81e-01 , 1.95e-01) & 1.00e+00& $\Z_{0}$  & 0  & $\Z_0$  \\
 \hline
\end{tabular}
 \end{table}
 \begin{table}[h]
 \begin{tabular}{|c||c||c||c||c||c||c|}
 \hline
  \multicolumn{7}{|c|}{ $ \alpha = \pi / 32 =     0.09817 \qquad \Rightarrow \qquad k=    0.09849 $}  \\
\hline
  $\#$ & $(\dot{x}_0, \dot{y}_0 )$ & $ (\dot{x}_f ,  \dot{y}_f) $  &$ \|{\bf v}_f \|$& $\Z({\bf v}_0)$ & $N$ & Stop \\
\hline
 4.6.1  & (      1.000  ,      0.000 )& (-5.00e-01 , -0.00e+00) & 5.00e-01& $\Z_{12}$  & 1  & $\Z_0$  \\
  4.6.2  & (      0.999  ,      0.033 )& (9.17e-13 , -2.96e-13) & 9.64e-13& $\Z_{12}$  & 433  & $\|  \cdot \| $  \\
  4.6.3  & (      0.998  ,      0.066 )& (9.05e-13 , 2.93e-13) & 9.51e-13& $\Z_{12}$  & 436  & $\|  \cdot \| $  \\
  4.6.4  & (      0.995  ,      0.098 )& (9.39e-13 , 3.03e-13) & 9.86e-13& $\Z_{2}$  & 436  & $\|  \cdot \| $  \\
  4.6.5  & (      0.098  ,      0.995 )& (-1.67e-01 , 1.25e-02) & 1.67e-01& $\Z_{2}$  & 5  & $\Z_0$  \\
  4.6.6  & (      0.000  ,      1.000 )& (-2.51e-01 , -3.75e-04) & 2.51e-01& $\Z_{2}$  & 4  & $\Z_0$  \\
  4.6.7  & (     -0.995  ,      0.098 )& (-9.95e-01 , 9.80e-02) & 1.00e+00& $\Z_{0}$  & 0  & $\Z_0$  \\
 \hline
\end{tabular}
 \end{table}
 \begin{table}[h]
 \begin{tabular}{|c||c||c||c||c||c||c|}
 \hline
  \multicolumn{7}{|c|}{ $ \alpha = \pi / 64 =     0.04909 \qquad \Rightarrow \qquad k=    0.04913 $}  \\
\hline
  $\#$ & $(\dot{x}_0, \dot{y}_0 )$ & $ (\dot{x}_f ,  \dot{y}_f) $  &$ \|{\bf v}_f \|$& $\Z({\bf v}_0)$ & $N$ & Stop \\
\hline
 4.7.1  & (      1.000  ,      0.000 )& (-5.00e-01 , -0.00e+00) & 5.00e-01& $\Z_{12}$  & 1  & $\Z_0$  \\
  4.7.2  & (      1.000  ,      0.016 )& (9.87e-13 , -1.48e-13) & 9.98e-13& $\Z_{12}$  & 1869  & $\|  \cdot \| $  \\
  4.7.3  & (      0.999  ,      0.033 )& (9.83e-13 , 1.48e-13) & 9.94e-13& $\Z_{12}$  & 1872  & $\|  \cdot \| $  \\
  4.7.4  & (      0.999  ,      0.049 )& (9.76e-13 , -1.47e-13) & 9.87e-13& $\Z_{2}$  & 1873  & $\|  \cdot \| $  \\
  4.7.5  & (      0.049  ,      0.999 )& (-9.21e-02 , 3.09e-03) & 9.21e-02& $\Z_{2}$  & 6  & $\Z_0$  \\
  4.7.6  & (      0.000  ,      1.000 )& (-1.38e-01 , -3.70e-03) & 1.38e-01& $\Z_{2}$  & 6  & $\Z_0$  \\
  4.7.7  & (     -0.999  ,      0.049 )& (-9.99e-01 , 4.91e-02) & 1.00e+00& $\Z_{0}$  & 0  & $\Z_0$  \\
 \hline
\end{tabular}
 \end{table}
 \clearpage
\subsection*{ Non ideal case: $\varepsilon=0.25$. }
\begin{table}[h]
 \begin{tabular}{|c||c||c||c||c||c||c|}
 \hline
  \multicolumn{7}{|c|}{ $ \alpha = \pi / 4 =     0.78540 \qquad \Rightarrow \qquad k = 1 $}  \\
\hline
  $\#$ & $(\dot{x}_0, \dot{y}_0 )$ & $ (\dot{x}_f ,  \dot{y}_f) $  &$ \|{\bf v}_f \|$& $\Z({\bf v}_0)$ & $N$ & Stop \\
\hline
 5.1.1  & (      1.000  ,      0.000 )& (-2.50e-01 , -0.00e+00) & 2.50e-01& $\Z_{12}$  & 1  & $\Z_0$  \\
  5.1.2  & (      0.949  ,      0.316 )& (-2.37e-01 , -7.91e-02) & 2.50e-01& $\Z_{12}$  & 1  & $\Z_0$  \\
  5.1.3  & (      0.832  ,      0.555 )& (-2.08e-01 , -1.39e-01) & 2.50e-01& $\Z_{12}$  & 1  & $\Z_0$  \\
  5.1.4  & (      0.707  ,      0.707 )& (-1.77e-01 , -1.77e-01) & 2.50e-01& $\Z_{2}$  & 1  & $\Z_0$  \\
  5.1.5  & (      0.707  ,      0.707 )& (-1.77e-01 , -1.77e-01) & 2.50e-01& $\Z_{2}$  & 1  & $\Z_0$  \\
  5.1.6  & (      0.000  ,      1.000 )& (-6.25e-01 , 3.75e-01) & 7.29e-01& $\Z_{2}$  & 1  & $\Z_0$  \\
  5.1.7  & (     -0.707  ,      0.707 )& (-7.07e-01 , 7.07e-01) & 1.00e+00& $\Z_{0}$  & 0  & $\Z_0$  \\
 \hline
\end{tabular}
 \end{table}
 \begin{table}[h]
 \begin{tabular}{|c||c||c||c||c||c||c|}
 \hline
  \multicolumn{7}{|c|}{ $ \alpha = \pi / 6 =     0.52360 \qquad \Rightarrow \qquad k=    0.57735 $}  \\
\hline
  $\#$ & $(\dot{x}_0, \dot{y}_0 )$ & $ (\dot{x}_f ,  \dot{y}_f) $  &$ \|{\bf v}_f \|$& $\Z({\bf v}_0)$ & $N$ & Stop \\
\hline
 5.2.1  & (      1.000  ,      0.000 )& (-2.50e-01 , -0.00e+00) & 2.50e-01& $\Z_{12}$  & 1  & $\Z_0$  \\
  5.2.2  & (      0.982  ,      0.189 )& (-2.45e-01 , -4.72e-02) & 2.50e-01& $\Z_{12}$  & 2  & $\Z_0$  \\
  5.2.3  & (      0.933  ,      0.359 )& (-1.23e-01 , 6.58e-02) & 1.39e-01& $\Z_{12}$  & 3  & $\Z_0$  \\
  5.2.4  & (      0.866  ,      0.500 )& (-8.96e-02 , 1.66e-02) & 9.12e-02& $\Z_{2}$  & 3  & $\Z_0$  \\
  5.2.5  & (      0.500  ,      0.866 )& (-2.03e-01 , -8.12e-02) & 2.19e-01& $\Z_{2}$  & 2  & $\Z_0$  \\
  5.2.6  & (      0.000  ,      1.000 )& (-5.41e-01 , 6.25e-02) & 5.45e-01& $\Z_{2}$  & 1  & $\Z_0$  \\
  5.2.7  & (     -0.866  ,      0.500 )& (-8.66e-01 , 5.00e-01) & 1.00e+00& $\Z_{0}$  & 0  & $\Z_0$  \\
 \hline
\end{tabular}
 \end{table}
 \begin{table}[h]
 \begin{tabular}{|c||c||c||c||c||c||c|}
 \hline
  \multicolumn{7}{|c|}{ $ \alpha = \pi / 8 =     0.39270 \qquad \Rightarrow \qquad k=    0.41421 $}  \\
\hline
  $\#$ & $(\dot{x}_0, \dot{y}_0 )$ & $ (\dot{x}_f ,  \dot{y}_f) $  &$ \|{\bf v}_f \|$& $\Z({\bf v}_0)$ & $N$ & Stop \\
\hline
 5.3.1  & (      1.000  ,      0.000 )& (-2.50e-01 , -0.00e+00) & 2.50e-01& $\Z_{12}$  & 1  & $\Z_0$  \\
  5.3.2  & (      0.991  ,      0.137 )& (-1.36e-01 , 4.14e-02) & 1.42e-01& $\Z_{12}$  & 3  & $\Z_0$  \\
  5.3.3  & (      0.964  ,      0.266 )& (-2.28e-02 , 7.07e-03) & 2.39e-02& $\Z_{12}$  & 5  & $\Z_0$  \\
  5.3.4  & (      0.924  ,      0.383 )& (-1.15e-02 , -1.95e-03) & 1.16e-02& $\Z_{2}$  & 6  & $\Z_0$  \\
  5.3.5  & (      0.383  ,      0.924 )& (-1.80e-01 , -2.68e-02) & 1.82e-01& $\Z_{2}$  & 2  & $\Z_0$  \\
  5.3.6  & (      0.000  ,      1.000 )& (-4.42e-01 , -6.69e-02) & 4.47e-01& $\Z_{2}$  & 1  & $\Z_0$  \\
  5.3.7  & (     -0.924  ,      0.383 )& (-9.24e-01 , 3.83e-01) & 1.00e+00& $\Z_{0}$  & 0  & $\Z_0$  \\
 \hline
\end{tabular}
 \end{table}
 \begin{table}[h]
 \begin{tabular}{|c||c||c||c||c||c||c|}
 \hline
  \multicolumn{7}{|c|}{ $ \alpha = \pi / 12 =     0.26180 \qquad \Rightarrow \qquad k=    0.26795 $}  \\
\hline
  $\#$ & $(\dot{x}_0, \dot{y}_0 )$ & $ (\dot{x}_f ,  \dot{y}_f) $  &$ \|{\bf v}_f \|$& $\Z({\bf v}_0)$ & $N$ & Stop \\
\hline
 5.4.1  & (      1.000  ,      0.000 )& (-2.50e-01 , -0.00e+00) & 2.50e-01& $\Z_{12}$  & 1  & $\Z_0$  \\
  5.4.2  & (      0.996  ,      0.089 )& (6.96e-13 , 3.73e-13) & 7.90e-13& $\Z_{12}$  & 92  & $\|  \cdot \| $  \\
  5.4.3  & (      0.984  ,      0.176 )& (7.93e-13 , 4.25e-13) & 9.00e-13& $\Z_{12}$  & 96  & $\|  \cdot \| $  \\
  5.4.4  & (      0.966  ,      0.259 )& (7.32e-13 , -3.93e-13) & 8.31e-13& $\Z_{2}$  & 97  & $\|  \cdot \| $  \\
  5.4.5  & (      0.259  ,      0.966 )& (-1.35e-01 , 1.99e-02) & 1.36e-01& $\Z_{2}$  & 2  & $\Z_0$  \\
  5.4.6  & (      0.000  ,      1.000 )& (-3.38e-01 , -7.00e-02) & 3.45e-01& $\Z_{2}$  & 2  & $\Z_0$  \\
  5.4.7  & (     -0.966  ,      0.259 )& (-9.66e-01 , 2.59e-01) & 1.00e+00& $\Z_{0}$  & 0  & $\Z_0$  \\
 \hline
\end{tabular}
 \end{table}
 \begin{table}[h]
 \begin{tabular}{|c||c||c||c||c||c||c|}
 \hline
  \multicolumn{7}{|c|}{ $ \alpha = \pi / 16 =     0.19635 \qquad \Rightarrow \qquad k=    0.19891 $}  \\
\hline
  $\#$ & $(\dot{x}_0, \dot{y}_0 )$ & $ (\dot{x}_f ,  \dot{y}_f) $  &$ \|{\bf v}_f \|$& $\Z({\bf v}_0)$ & $N$ & Stop \\
\hline
 5.5.1  & (      1.000  ,      0.000 )& (-2.50e-01 , -0.00e+00) & 2.50e-01& $\Z_{12}$  & 1  & $\Z_0$  \\
  5.5.2  & (      0.998  ,      0.066 )& (8.25e-13 , 2.97e-13) & 8.77e-13& $\Z_{12}$  & 190  & $\|  \cdot \| $  \\
  5.5.3  & (      0.991  ,      0.131 )& (8.76e-13 , -3.16e-13) & 9.31e-13& $\Z_{12}$  & 193  & $\|  \cdot \| $  \\
  5.5.4  & (      0.981  ,      0.195 )& (8.36e-13 , 3.01e-13) & 8.89e-13& $\Z_{2}$  & 194  & $\|  \cdot \| $  \\
  5.5.5  & (      0.195  ,      0.981 )& (-1.09e-01 , 1.75e-02) & 1.11e-01& $\Z_{2}$  & 3  & $\Z_0$  \\
  5.5.6  & (      0.000  ,      1.000 )& (-2.76e-01 , -1.62e-02) & 2.77e-01& $\Z_{2}$  & 2  & $\Z_0$  \\
  5.5.7  & (     -0.981  ,      0.195 )& (-9.81e-01 , 1.95e-01) & 1.00e+00& $\Z_{0}$  & 0  & $\Z_0$  \\
 \hline
\end{tabular}
 \end{table}
 \begin{table}[h]
 \begin{tabular}{|c||c||c||c||c||c||c|}
 \hline
  \multicolumn{7}{|c|}{ $ \alpha = \pi / 32 =     0.09817 \qquad \Rightarrow \qquad k=    0.09849 $}  \\
\hline
  $\#$ & $(\dot{x}_0, \dot{y}_0 )$ & $ (\dot{x}_f ,  \dot{y}_f) $  &$ \|{\bf v}_f \|$& $\Z({\bf v}_0)$ & $N$ & Stop \\
\hline
 5.6.1  & (      1.000  ,      0.000 )& (-2.50e-01 , -0.00e+00) & 2.50e-01& $\Z_{12}$  & 1  & $\Z_0$  \\
  5.6.2  & (      0.999  ,      0.033 )& (9.84e-13 , 1.64e-13) & 9.98e-13& $\Z_{12}$  & 836  & $\|  \cdot \| $  \\
  5.6.3  & (      0.998  ,      0.066 )& (9.63e-13 , 1.61e-13) & 9.76e-13& $\Z_{12}$  & 840  & $\|  \cdot \| $  \\
  5.6.4  & (      0.995  ,      0.098 )& (9.84e-13 , 1.64e-13) & 9.97e-13& $\Z_{2}$  & 840  & $\|  \cdot \| $  \\
  5.6.5  & (      0.098  ,      0.995 )& (-6.08e-02 , -5.81e-03) & 6.11e-02& $\Z_{2}$  & 4  & $\Z_0$  \\
  5.6.6  & (      0.000  ,      1.000 )& (-1.53e-01 , 8.29e-03) & 1.53e-01& $\Z_{2}$  & 3  & $\Z_0$  \\
  5.6.7  & (     -0.995  ,      0.098 )& (-9.95e-01 , 9.80e-02) & 1.00e+00& $\Z_{0}$  & 0  & $\Z_0$  \\
 \hline
\end{tabular}
 \end{table}
 \begin{table}[h]
 \begin{tabular}{|c||c||c||c||c||c||c|}
 \hline
  \multicolumn{7}{|c|}{ $ \alpha = \pi / 64 =     0.04909 \qquad \Rightarrow \qquad k=    0.04913 $}  \\
\hline
  $\#$ & $(\dot{x}_0, \dot{y}_0 )$ & $ (\dot{x}_f ,  \dot{y}_f) $  &$ \|{\bf v}_f \|$& $\Z({\bf v}_0)$ & $N$ & Stop \\
\hline
 5.7.1  & (      1.000  ,      0.000 )& (-2.50e-01 , -0.00e+00) & 2.50e-01& $\Z_{12}$  & 1  & $\Z_0$  \\
  5.7.2  & (      1.000  ,      0.016 )& (9.92e-13 , -8.15e-14) & 9.95e-13& $\Z_{12}$  & 3417  & $\|  \cdot \| $  \\
  5.7.3  & (      0.999  ,      0.033 )& (9.94e-13 , 8.18e-14) & 9.98e-13& $\Z_{12}$  & 3420  & $\|  \cdot \| $  \\
  5.7.4  & (      0.999  ,      0.049 )& (9.92e-13 , -8.16e-14) & 9.95e-13& $\Z_{2}$  & 3421  & $\|  \cdot \| $  \\
  5.7.5  & (      0.049  ,      0.999 )& (-3.19e-02 , 1.39e-03) & 3.19e-02& $\Z_{2}$  & 4  & $\Z_0$  \\
  5.7.6  & (      0.000  ,      1.000 )& (-7.98e-02 , -2.53e-03) & 7.99e-02& $\Z_{2}$  & 4  & $\Z_0$  \\
  5.7.7  & (     -0.999  ,      0.049 )& (-9.99e-01 , 4.91e-02) & 1.00e+00& $\Z_{0}$  & 0  & $\Z_0$  \\
 \hline
\end{tabular}
 \end{table}
 \clearpage
\subsection*{ Non ideal case: $\varepsilon=0.05$. }
\begin{table}[h]
 \begin{tabular}{|c||c||c||c||c||c||c|}
 \hline
  \multicolumn{7}{|c|}{ $ \alpha = \pi / 4 =     0.78540 \qquad \Rightarrow \qquad k = 1 $}  \\
\hline
  $\#$ & $(\dot{x}_0, \dot{y}_0 )$ & $ (\dot{x}_f ,  \dot{y}_f) $  &$ \|{\bf v}_f \|$& $\Z({\bf v}_0)$ & $N$ & Stop \\
\hline
 6.1.1  & (      1.000  ,      0.000 )& (-5.00e-02 , -0.00e+00) & 5.00e-02& $\Z_{12}$  & 1  & $\Z_0$  \\
  6.1.2  & (      0.949  ,      0.316 )& (-4.74e-02 , -1.58e-02) & 5.00e-02& $\Z_{12}$  & 1  & $\Z_0$  \\
  6.1.3  & (      0.832  ,      0.555 )& (-4.16e-02 , -2.77e-02) & 5.00e-02& $\Z_{12}$  & 1  & $\Z_0$  \\
  6.1.4  & (      0.707  ,      0.707 )& (-3.54e-02 , -3.54e-02) & 5.00e-02& $\Z_{2}$  & 1  & $\Z_0$  \\
  6.1.5  & (      0.707  ,      0.707 )& (-3.54e-02 , -3.54e-02) & 5.00e-02& $\Z_{2}$  & 1  & $\Z_0$  \\
  6.1.6  & (      0.000  ,      1.000 )& (-5.25e-01 , 4.75e-01) & 7.08e-01& $\Z_{2}$  & 1  & $\Z_0$  \\
  6.1.7  & (     -0.707  ,      0.707 )& (-7.07e-01 , 7.07e-01) & 1.00e+00& $\Z_{0}$  & 0  & $\Z_0$  \\
 \hline
\end{tabular}
 \end{table}
 \begin{table}[h]
 \begin{tabular}{|c||c||c||c||c||c||c|}
 \hline
  \multicolumn{7}{|c|}{ $ \alpha = \pi / 6 =     0.52360 \qquad \Rightarrow \qquad k=    0.57735 $}  \\
\hline
  $\#$ & $(\dot{x}_0, \dot{y}_0 )$ & $ (\dot{x}_f ,  \dot{y}_f) $  &$ \|{\bf v}_f \|$& $\Z({\bf v}_0)$ & $N$ & Stop \\
\hline
 6.2.1  & (      1.000  ,      0.000 )& (-5.00e-02 , -0.00e+00) & 5.00e-02& $\Z_{12}$  & 1  & $\Z_0$  \\
  6.2.2  & (      0.982  ,      0.189 )& (-4.91e-02 , -9.45e-03) & 5.00e-02& $\Z_{12}$  & 2  & $\Z_0$  \\
  6.2.3  & (      0.933  ,      0.359 )& (4.71e-13 , 3.49e-13) & 5.86e-13& $\Z_{12}$  & 30  & $\|  \cdot \| $  \\
  6.2.4  & (      0.866  ,      0.500 )& (5.08e-13 , -3.77e-13) & 6.33e-13& $\Z_{2}$  & 31  & $\|  \cdot \| $  \\
  6.2.5  & (      0.500  ,      0.866 )& (-3.81e-02 , -2.06e-02) & 4.33e-02& $\Z_{2}$  & 2  & $\Z_0$  \\
  6.2.6  & (      0.000  ,      1.000 )& (-4.55e-01 , 2.12e-01) & 5.02e-01& $\Z_{2}$  & 1  & $\Z_0$  \\
  6.2.7  & (     -0.866  ,      0.500 )& (-8.66e-01 , 5.00e-01) & 1.00e+00& $\Z_{0}$  & 0  & $\Z_0$  \\
 \hline
\end{tabular}
 \end{table}
 \begin{table}[h]
 \begin{tabular}{|c||c||c||c||c||c||c|}
 \hline
  \multicolumn{7}{|c|}{ $ \alpha = \pi / 8 =     0.39270 \qquad \Rightarrow \qquad k=    0.41421 $}  \\
\hline
  $\#$ & $(\dot{x}_0, \dot{y}_0 )$ & $ (\dot{x}_f ,  \dot{y}_f) $  &$ \|{\bf v}_f \|$& $\Z({\bf v}_0)$ & $N$ & Stop \\
\hline
 6.3.1  & (      1.000  ,      0.000 )& (-5.00e-02 , -0.00e+00) & 5.00e-02& $\Z_{12}$  & 1  & $\Z_0$  \\
  6.3.2  & (      0.991  ,      0.137 )& (6.35e-13 , -3.06e-13) & 7.04e-13& $\Z_{12}$  & 65  & $\|  \cdot \| $  \\
  6.3.3  & (      0.964  ,      0.266 )& (8.37e-13 , 4.03e-13) & 9.29e-13& $\Z_{12}$  & 68  & $\|  \cdot \| $  \\
  6.3.4  & (      0.924  ,      0.383 )& (7.47e-13 , -3.60e-13) & 8.29e-13& $\Z_{2}$  & 69  & $\|  \cdot \| $  \\
  6.3.5  & (      0.383  ,      0.924 )& (-3.33e-02 , -1.19e-02) & 3.54e-02& $\Z_{2}$  & 2  & $\Z_0$  \\
  6.3.6  & (      0.000  ,      1.000 )& (-3.71e-01 , 1.04e-01) & 3.85e-01& $\Z_{2}$  & 1  & $\Z_0$  \\
  6.3.7  & (     -0.924  ,      0.383 )& (-9.24e-01 , 3.83e-01) & 1.00e+00& $\Z_{0}$  & 0  & $\Z_0$  \\
 \hline
\end{tabular}
 \end{table}
 \begin{table}[h]
 \begin{tabular}{|c||c||c||c||c||c||c|}
 \hline
  \multicolumn{7}{|c|}{ $ \alpha = \pi / 12 =     0.26180 \qquad \Rightarrow \qquad k=    0.26795 $}  \\
\hline
  $\#$ & $(\dot{x}_0, \dot{y}_0 )$ & $ (\dot{x}_f ,  \dot{y}_f) $  &$ \|{\bf v}_f \|$& $\Z({\bf v}_0)$ & $N$ & Stop \\
\hline
 6.4.1  & (      1.000  ,      0.000 )& (-5.00e-02 , -0.00e+00) & 5.00e-02& $\Z_{12}$  & 1  & $\Z_0$  \\
  6.4.2  & (      0.996  ,      0.089 )& (9.29e-13 , -2.80e-13) & 9.70e-13& $\Z_{12}$  & 167  & $\|  \cdot \| $  \\
  6.4.3  & (      0.984  ,      0.176 )& (8.28e-13 , -2.50e-13) & 8.65e-13& $\Z_{12}$  & 171  & $\|  \cdot \| $  \\
  6.4.4  & (      0.966  ,      0.259 )& (9.28e-13 , -2.80e-13) & 9.69e-13& $\Z_{2}$  & 171  & $\|  \cdot \| $  \\
  6.4.5  & (      0.259  ,      0.966 )& (-2.47e-02 , -4.38e-03) & 2.51e-02& $\Z_{2}$  & 2  & $\Z_0$  \\
  6.4.6  & (      0.000  ,      1.000 )& (-2.63e-01 , 2.03e-02) & 2.63e-01& $\Z_{2}$  & 1  & $\Z_0$  \\
  6.4.7  & (     -0.966  ,      0.259 )& (-9.66e-01 , 2.59e-01) & 1.00e+00& $\Z_{0}$  & 0  & $\Z_0$  \\
 \hline
\end{tabular}
 \end{table}
 \begin{table}[h]
 \begin{tabular}{|c||c||c||c||c||c||c|}
 \hline
  \multicolumn{7}{|c|}{ $ \alpha = \pi / 16 =     0.19635 \qquad \Rightarrow \qquad k=    0.19891 $}  \\
\hline
  $\#$ & $(\dot{x}_0, \dot{y}_0 )$ & $ (\dot{x}_f ,  \dot{y}_f) $  &$ \|{\bf v}_f \|$& $\Z({\bf v}_0)$ & $N$ & Stop \\
\hline
 6.5.1  & (      1.000  ,      0.000 )& (-5.00e-02 , -0.00e+00) & 5.00e-02& $\Z_{12}$  & 1  & $\Z_0$  \\
  6.5.2  & (      0.998  ,      0.066 )& (9.54e-13 , -2.12e-13) & 9.77e-13& $\Z_{12}$  & 309  & $\|  \cdot \| $  \\
  6.5.3  & (      0.991  ,      0.131 )& (9.10e-13 , -2.02e-13) & 9.33e-13& $\Z_{12}$  & 313  & $\|  \cdot \| $  \\
  6.5.4  & (      0.981  ,      0.195 )& (9.69e-13 , -2.15e-13) & 9.93e-13& $\Z_{2}$  & 313  & $\|  \cdot \| $  \\
  6.5.5  & (      0.195  ,      0.981 )& (-1.92e-02 , -1.47e-03) & 1.93e-02& $\Z_{2}$  & 2  & $\Z_0$  \\
  6.5.6  & (      0.000  ,      1.000 )& (-2.01e-01 , -1.00e-02) & 2.01e-01& $\Z_{2}$  & 1  & $\Z_0$  \\
  6.5.7  & (     -0.981  ,      0.195 )& (-9.81e-01 , 1.95e-01) & 1.00e+00& $\Z_{0}$  & 0  & $\Z_0$  \\
 \hline
\end{tabular}
 \end{table}
 \begin{table}[h]
 \begin{tabular}{|c||c||c||c||c||c||c|}
 \hline
  \multicolumn{7}{|c|}{ $ \alpha = \pi / 32 =     0.09817 \qquad \Rightarrow \qquad k=    0.09849 $}  \\
\hline
  $\#$ & $(\dot{x}_0, \dot{y}_0 )$ & $ (\dot{x}_f ,  \dot{y}_f) $  &$ \|{\bf v}_f \|$& $\Z({\bf v}_0)$ & $N$ & Stop \\
\hline
 6.6.1  & (      1.000  ,      0.000 )& (-5.00e-02 , -0.00e+00) & 5.00e-02& $\Z_{12}$  & 1  & $\Z_0$  \\
  6.6.2  & (      0.999  ,      0.033 )& (9.80e-13 , 1.07e-13) & 9.85e-13& $\Z_{12}$  & 1282  & $\|  \cdot \| $  \\
  6.6.3  & (      0.998  ,      0.066 )& (9.75e-13 , 1.06e-13) & 9.81e-13& $\Z_{12}$  & 1286  & $\|  \cdot \| $  \\
  6.6.4  & (      0.995  ,      0.098 )& (9.90e-13 , 1.08e-13) & 9.96e-13& $\Z_{2}$  & 1286  & $\|  \cdot \| $  \\
  6.6.5  & (      0.098  ,      0.995 )& (-1.00e-02 , 9.60e-04) & 1.00e-02& $\Z_{2}$  & 3  & $\Z_0$  \\
  6.6.6  & (      0.000  ,      1.000 )& (-1.05e-01 , -8.90e-03) & 1.06e-01& $\Z_{2}$  & 2  & $\Z_0$  \\
  6.6.7  & (     -0.995  ,      0.098 )& (-9.95e-01 , 9.80e-02) & 1.00e+00& $\Z_{0}$  & 0  & $\Z_0$  \\
 \hline
\end{tabular}
 \end{table}
 \begin{table}[h]
 \begin{tabular}{|c||c||c||c||c||c||c|}
 \hline
  \multicolumn{7}{|c|}{ $ \alpha = \pi / 64 =     0.04909 \qquad \Rightarrow \qquad k=    0.04913 $}  \\
\hline
  $\#$ & $(\dot{x}_0, \dot{y}_0 )$ & $ (\dot{x}_f ,  \dot{y}_f) $  &$ \|{\bf v}_f \|$& $\Z({\bf v}_0)$ & $N$ & Stop \\
\hline
 6.7.1  & (      1.000  ,      0.000 )& (-5.00e-02 , -0.00e+00) & 5.00e-02& $\Z_{12}$  & 1  & $\Z_0$  \\
  6.7.2  & (      1.000  ,      0.016 )& (9.98e-13 , 5.42e-14) & 1.00e-12& $\Z_{12}$  & 5172  & $\|  \cdot \| $  \\
  6.7.3  & (      0.999  ,      0.033 )& (9.97e-13 , 5.42e-14) & 9.99e-13& $\Z_{12}$  & 5176  & $\|  \cdot \| $  \\
  6.7.4  & (      0.999  ,      0.049 )& (9.96e-13 , -5.41e-14) & 9.97e-13& $\Z_{2}$  & 5177  & $\|  \cdot \| $  \\
  6.7.5  & (      0.049  ,      0.999 )& (-5.12e-03 , 1.52e-04) & 5.12e-03& $\Z_{2}$  & 3  & $\Z_0$  \\
  6.7.6  & (      0.000  ,      1.000 )& (-5.38e-02 , -3.94e-04) & 5.38e-02& $\Z_{2}$  & 2  & $\Z_0$  \\
  6.7.7  & (     -0.999  ,      0.049 )& (-9.99e-01 , 4.91e-02) & 1.00e+00& $\Z_{0}$  & 0  & $\Z_0$  \\
 \hline
\end{tabular}
 \end{table}
 \clearpage
\subsection*{ Non ideal case: $\varepsilon=0$. }

\begin{table}[h]
 \begin{tabular}{|c||c||c||c||c||c||c|}
 \hline
  \multicolumn{7}{|c|}{ $ \alpha = \pi / 4 =     0.78540 \qquad \Rightarrow \qquad k = 1 $}  \\
\hline
  $\#$ & $(\dot{x}_0, \dot{y}_0 )$ & $ (\dot{x}_f ,  \dot{y}_f) $  &$ \|{\bf v}_f \|$& $\Z({\bf v}_0)$ & $N$ & Stop \\
\hline
 7.1.1  & (      1.000  ,      0.000 )& (0.00e+00 , -0.00e+00) & 0.00e+00& $\Z_{12}$  & 1  & $\Z_0$  \\
  7.1.2  & (      0.949  ,      0.316 )& (2.11e-17 , -6.32e-17) & 6.66e-17& $\Z_{12}$  & 1  & $\Z_0$  \\
  7.1.3  & (      0.832  ,      0.555 )& (5.68e-17 , -8.53e-17) & 1.02e-16& $\Z_{12}$  & 1  & $\Z_0$  \\
  7.1.4  & (      0.707  ,      0.707 )& (1.11e-16 , -5.55e-17) & 1.24e-16& $\Z_{2}$  & 1  & $\Z_0$  \\
  7.1.5  & (      0.707  ,      0.707 )& (-5.55e-17 , 1.11e-16) & 1.24e-16& $\Z_{2}$  & 1  & $\Z_0$  \\
  7.1.6  & (      0.000  ,      1.000 )& (-5.00e-01 , 5.00e-01) & 7.07e-01& $\Z_{2}$  & 1  & $\Z_0$  \\
  7.1.7  & (     -0.707  ,      0.707 )& (-7.07e-01 , 7.07e-01) & 1.00e+00& $\Z_{0}$  & 0  & $\Z_0$  \\
 \hline
\end{tabular}
 \end{table}
 \begin{table}[h]
 \begin{tabular}{|c||c||c||c||c||c||c|}
 \hline
  \multicolumn{7}{|c|}{ $ \alpha = \pi / 6 =     0.52360 \qquad \Rightarrow \qquad k=    0.57735 $}  \\
\hline
  $\#$ & $(\dot{x}_0, \dot{y}_0 )$ & $ (\dot{x}_f ,  \dot{y}_f) $  &$ \|{\bf v}_f \|$& $\Z({\bf v}_0)$ & $N$ & Stop \\
\hline
 7.2.1  & (      1.000  ,      0.000 )& (0.00e+00 , -0.00e+00) & 0.00e+00& $\Z_{12}$  & 1  & $\Z_0$  \\
  7.2.2  & (      0.982  ,      0.189 )& (0.00e+00 , 0.00e+00) & 0.00e+00& $\Z_{12}$  & 2  & $\Z_0$  \\
  7.2.3  & (      0.933  ,      0.359 )& (4.85e-13 , 2.80e-13) & 5.60e-13& $\Z_{12}$  & 40  & $\|  \cdot \| $  \\
  7.2.4  & (      0.866  ,      0.500 )& (7.88e-13 , 4.55e-13) & 9.09e-13& $\Z_{2}$  & 40  & $\|  \cdot \| $  \\
  7.2.5  & (      0.500  ,      0.866 )& (0.00e+00 , 0.00e+00) & 0.00e+00& $\Z_{2}$  & 1  & $\Z_0$  \\
  7.2.6  & (      0.000  ,      1.000 )& (-4.33e-01 , 2.50e-01) & 5.00e-01& $\Z_{2}$  & 1  & $\Z_0$  \\
  7.2.7  & (     -0.866  ,      0.500 )& (-8.66e-01 , 5.00e-01) & 1.00e+00& $\Z_{0}$  & 0  & $\Z_0$  \\
 \hline
\end{tabular}
 \end{table}
 \begin{table}[h]
 \begin{tabular}{|c||c||c||c||c||c||c|}
 \hline
  \multicolumn{7}{|c|}{ $ \alpha = \pi / 8 =     0.39270 \qquad \Rightarrow \qquad k=    0.41421 $}  \\
\hline
  $\#$ & $(\dot{x}_0, \dot{y}_0 )$ & $ (\dot{x}_f ,  \dot{y}_f) $  &$ \|{\bf v}_f \|$& $\Z({\bf v}_0)$ & $N$ & Stop \\
\hline
 7.3.1  & (      1.000  ,      0.000 )& (0.00e+00 , -0.00e+00) & 0.00e+00& $\Z_{12}$  & 1  & $\Z_0$  \\
  7.3.2  & (      0.991  ,      0.137 )& (6.87e-13 , -2.85e-13) & 7.44e-13& $\Z_{12}$  & 77  & $\|  \cdot \| $  \\
  7.3.3  & (      0.964  ,      0.266 )& (6.64e-13 , 2.75e-13) & 7.19e-13& $\Z_{12}$  & 80  & $\|  \cdot \| $  \\
  7.3.4  & (      0.924  ,      0.383 )& (8.40e-13 , 3.48e-13) & 9.09e-13& $\Z_{2}$  & 80  & $\|  \cdot \| $  \\
  7.3.5  & (      0.383  ,      0.924 )& (4.74e-17 , 0.00e+00) & 4.74e-17& $\Z_{2}$  & 1  & $\Z_0$  \\
  7.3.6  & (      0.000  ,      1.000 )& (-3.54e-01 , 1.46e-01) & 3.83e-01& $\Z_{2}$  & 1  & $\Z_0$  \\
  7.3.7  & (     -0.924  ,      0.383 )& (-9.24e-01 , 3.83e-01) & 1.00e+00& $\Z_{0}$  & 0  & $\Z_0$  \\
 \hline
\end{tabular}
 \end{table}
 \begin{table}[h]
 \begin{tabular}{|c||c||c||c||c||c||c|}
 \hline
  \multicolumn{7}{|c|}{ $ \alpha = \pi / 12 =     0.26180 \qquad \Rightarrow \qquad k=    0.26795 $}  \\
\hline
  $\#$ & $(\dot{x}_0, \dot{y}_0 )$ & $ (\dot{x}_f ,  \dot{y}_f) $  &$ \|{\bf v}_f \|$& $\Z({\bf v}_0)$ & $N$ & Stop \\
\hline
 7.4.1  & (      1.000  ,      0.000 )& (0.00e+00 , -0.00e+00) & 0.00e+00& $\Z_{12}$  & 1  & $\Z_0$  \\
  7.4.2  & (      0.996  ,      0.089 )& (8.55e-13 , -2.29e-13) & 8.86e-13& $\Z_{12}$  & 189  & $\|  \cdot \| $  \\
  7.4.3  & (      0.984  ,      0.176 )& (8.85e-13 , 2.37e-13) & 9.16e-13& $\Z_{12}$  & 192  & $\|  \cdot \| $  \\
  7.4.4  & (      0.966  ,      0.259 )& (8.48e-13 , -2.27e-13) & 8.78e-13& $\Z_{2}$  & 193  & $\|  \cdot \| $  \\
  7.4.5  & (      0.259  ,      0.966 )& (0.00e+00 , 0.00e+00) & 0.00e+00& $\Z_{2}$  & 1  & $\Z_0$  \\
  7.4.6  & (      0.000  ,      1.000 )& (-2.50e-01 , 6.70e-02) & 2.59e-01& $\Z_{2}$  & 1  & $\Z_0$  \\
  7.4.7  & (     -0.966  ,      0.259 )& (-9.66e-01 , 2.59e-01) & 1.00e+00& $\Z_{0}$  & 0  & $\Z_0$  \\
 \hline
\end{tabular}%
 \end{table}

\clearpage

\section{Appendix. Theoretical aspects of the multiple impact}

In this section we recall the geometric structures and properties
involved in the study of multiple impacts and the constitutive
characterization of the multiple contact/impact, as presented in
\cite{Pasquero2016Multiple}, for the particular system given by
the disk in the corner. The arguments are presented in a very
synthetic way, just to made the paper self consistent. For a more
general and exhaustive description we refer to
\cite{Pasquero2016Multiple} and the references therein.

\subsection{Geometry of the system}

The geometric setup suitable to study in a frame independent and
time dependent way the mechanical system formed by a rigid disk of
mass $m$ having multiple contact/impact with two walls forming a
corner consists in:

\begin{itemize}

\item a bundle $\pi_t : \M \to \Euclideo$, being $\M$  a
$(3+1)$--dimensional differentiable manifold and $\Euclideo$ the
affine time line. The elements of $\M$ are called space--time
configurations of the system and the simplest coordinates
describing $\M$ are the fibred coordinates $(t, x, y, \vth)$,
where $t$ is the time coordinate, $(x,y)$ are the coordinate of
the center of the disk and $\vth$ is the orientation of the disk;

\item the first jet--extension $\pi: \J \to \M$ of the bundle
$\M$, representing the space of absolute velocities of the system.
It is a $(6+1)$--dimensional affine subbundle of the tangent
bundle $T(\M)$ of $\M$ that can be referred to jet--coordinates
$(t, x, y, \vth, \dot{x}, \dot{y}, \dot{\vth})$. Using these local
coordinates, the elements of $\J$ have the form ${\bf p} =
\frac{\partial}{\partial t} + \dot{x} \frac{\partial}{\partial x}
+ \dot{y} \frac{\partial}{\partial y}+ \dot{\vth}
\frac{\partial}{\partial \vth}$;

\item the vertical vector bundle $\pi: V(\M) \to \M$ of the
vectors of $T(\M)$ that are vertical with respect to $\pi_t$, that
is, that are tangent to the fibers of $\M$. The bundle $V(\M)$ is
the vector bundle modelling the affine bundle $\J$ and it
represents both the space of the relative velocities of the system
(once a frame of reference is assigned) and the space of possible
impulses acting on the system. It can be referred to the same
jet--coordinates $(t, x, y, \vth, \dot{x}, \dot{y}, \dot{\vth})$
of $\J$ and, using these local coordinates, the elements of
$V(\M)$ have the form ${\bf v} = v_{x} \frac{\partial}{\partial x}
+ v_{y} \frac{\partial}{\partial y}+ \omega_{\vth}
\frac{\partial}{\partial \vth}$;

\item a positive definite scalar product $\Phi : V(\M) \times_{\M}
V(\M) \to \Euclideo$, acting on the fibers of $V(\M)$. It is
usually called the vertical metric and it takes intrinsically into
account the mass properties of the system. Using again the local
coordinates $(t, x, y, \vth)$, the vertical metric is expressed by
the positive definite matrix $G \= diag(m,m,A)$ where $m$ is the
mass of the disk and $A$ its inertia momentum. A standard
calculation shows that, using coordinates $(t, \xi, \eta, \vth)$
with $\xi =kx + y , \eta = kx-y$ coordinates, the matrix
expression of the vertical metric is no more diagonal and it is
transformed in
\begin{eqnarray}
\Gamma \= \left(
\begin{array}{ccc}
\dfrac{m(1+k^2)}{4k^2} & \dfrac{m(1-k^2)}{4k^2} & 0\\ \\
\dfrac{m(1-k^2)}{4k^2} & \dfrac{m(1+k^2)}{4k^2} & 0 \\ \\
0 & 0 & A
\end{array} \right) \, ;
\end{eqnarray}

\item the class $\H_{\M}$ of the frames of reference of the system
(without any assumption of rigidity), that is the set of global
fibred sections ${\bf h}_{\M} : \M \to \J$. Using local
coordinates, the elements of $\H_{\M}$ have the form ${\bf h}_{\M}
= \frac{\partial}{\partial t} + H^{x} \frac{\partial}{\partial x}
+ H^{y} \frac{\partial}{\partial y}+ H^{\vth}
\frac{\partial}{\partial \vth}$;

\item the pair of additional positional constraints $\SO, \ST$
describing the walls of the corner. The subbundle $i_1:\SO \to \M$
of $\M$ representing $\SO$ can be described by the cartesian
representation $kx-y=0$ or by the parametric representation given
by the immersion $(t,x,\vth) \rightsquigarrow (t,x,kx,\vth)$. The
subbundle $i_2:\ST \to \M$ of $\M$ representing $\ST$ can be
described by the cartesian representation $kx+y=0$ or by the
parametric representation given by the immersion $(t,x,\vth)
\rightsquigarrow (t,x,-kx,\vth)$. The single constraints $\SO,
\ST$ determine the multiple constraint $\SOT = \SO \cap \ST$ where
$i_{12}:\SOT \to \M$ is the subbundle of $\M$ described by the
cartesian representation $x=y=0$ or by the parametric
representation given by the immersion $(t,\vth) \rightsquigarrow
(t,0,0,\vth)$. The system is {\it in contact} with one or both the
constraints if its space--time configuration belongs to $\SO, \ST$
or $\SOT$. Each of the subbundles $\SO, \ST, \SOT$ determines its
first jet bundle $\JSO, \JST, \JSOT$ of the absolute velocities
tangent respectively to $\SO, \ST, \SOT$ (all of them affine
subbundles of $\J$), and its vertical vector bundle $V(\SO),
V(\ST), V(\SOT)$ (all of them vector subbundles of $V(\M)$);

\item the so called contact bundles
\begin{eqnarray*}
\begin{array}{ccc}
i_1^*(\J)  & \hookrightarrow & \J \\
i_2^*(\J) & \hookrightarrow & \J \\
i_{12}^*(\J)  & \hookrightarrow & \J \\ \\
i_1^*(V(\M)) & \hookrightarrow & V(\M) \\
i_2^*(V(\M)) & \hookrightarrow & V(\M) \\
i_{12}^*(V(\M)) & \hookrightarrow & V(\M),
\end{array}
\end{eqnarray*}
pull--back bundles of the absolute velocities $\J$ and impulses
$V(\M)$ respectively on $\SO,\ST,\SOT$ and representing all the
possible absolute velocities and impulses of the system when the
system is in contact with respectively $\SO,\ST,\SOT$. Of course,
being $\SOT =\SO\cap \ST$, we have $i_{12}^*(\J) = i_{1}^*(\J)
\cap i_{2}^*(\J)$ and $i_{12}^*(V(\M)) = i_{1}^*(V(\M)) \cap
i_{2}^*(V(\M))$. Since we suppose the disk in contact with both
$\SO,\ST$, then the geometric context where the multiple
contact/impact of the disk can be framed is given by
$i_{12}^*(\J)$ and $i_{12}^*(V(\M))$;

\item the set of projection operators determined by the vertical
metric $\Phi$. In particular, with obvious notation, these
projection operators determine natural splits
\begin{eqnarray*}
\begin{array}{ccccccc}
i_1^*(V(\M)) &\=& V(\SO) \oplus V^{\perp}(\SO) & \Leftrightarrow &
{\bf v} &\=&
{\bf v}^{\|}_1 + \velort_1({\bf p})\\
i_2^*(V(\M)) &\=& V(\ST) \oplus V^{\perp}(\ST) & \Leftrightarrow &
{\bf v} &\=&
{\bf v}^{\|}_2 + \velort_2({\bf p})\\
i_{12}^*(V(\M)) &\=& V(\SOT) \oplus V^{\perp}(\SOT) &
\Leftrightarrow &
{\bf v} &\=& {\bf v}^{\|}_{12} + \velort_{12}({\bf p})\\ \\
i_1^*(\J) &\=& \JSO \oplus V^{\perp}(\SO) & \Leftrightarrow & {\bf
p} &\=&
{\bf p}_1 + \velort_1({\bf p})\\
i_2^*(\J) &\=& \JST \oplus V^{\perp}(\ST) & \Leftrightarrow & {\bf
p} &\=&
{\bf p}_2 + \velort_2({\bf p})\\
i_{12}^*(\J) &\=& \JSOT \oplus V^{\perp}(\SOT) & \Leftrightarrow &
{\bf p} &\=& {\bf p}_{12} + \velort_{12}({\bf p}).
\end{array}
\end{eqnarray*}
Being $\SO, \ST$ of codimension $1$, also the orthogonal vector
subbundles $V^{\perp}(\SO), V^{\perp}(\ST)$ has codimension $1$,
so that we can introduce two vertical vectors $\uelortO, \uelortT$
such that $V^{\perp}(\SO)=Lin(\uelortO),
V^{\perp}(\ST)=Lin(\uelortT)$. The impact (with $\SO, \ST$ or
both) nature of an absolute velocity ${\bf p} \in i_{12}^*(\J)$ is
determined by the (suitably chosen) sign of the scalar products
$\Phi(\velort_1({\bf p}),\uelortO), \Phi(\velort_2({\bf
p}),\uelortT)$;

\item the subclasses $\H_{\SO}, \H_{\ST}, \H_{\SOT}$ of the frames
of reference of $\H_{\M}$ that are tangent to $\SO, \ST, \SOT$
respectively. Using local coordinates, the elements of $\H_{\SO},
\H_{\ST}, \H_{\SOT}$ have the forms
\begin{eqnarray*}
\begin{array}{l}
{\bf h}_{\SO} = \dfrac{\partial}{\partial t} + H^{x}
\dfrac{\partial}{\partial x} + k \, H^{x}
\dfrac{\partial}{\partial
y}+ H^{\vth} \dfrac{\partial}{\partial \vth} \\ \\
{\bf h}_{\ST} = \dfrac{\partial}{\partial t} + H^{x}
\dfrac{\partial}{\partial x} - k \, H^{x}
\dfrac{\partial}{\partial y}+ H^{\vth}
\dfrac{\partial}{\partial \vth} \\ \\
{\bf h}_{\SOT} = \dfrac{\partial}{\partial t} \phantom{+ H^{x}
\dfrac{\partial}{\partial x} - k \, H^{x}
\dfrac{\partial}{\partial y}} + H^{\vth} \dfrac{\partial}{\partial
\vth}
\end{array}
\end{eqnarray*}
The three subclasses represent the set of frames that can be
considered at rest with $\SO, \ST, \SOT$ respectively, that is the
sole frames for which the conservation of kinetic energy could
have an invariant meaning (see
\cite{Pasquero2005uni,Pasquero2016Multiple}).
\end{itemize}

\subsection{Constitutive law for ideal impact}

The frame independent description of an impulsive dynamic problem
in this geometric context consists in determining an element
${\bf{p}}_R \in i_{12}^*(\J)$, the right--velocity, once an
element ${\bf{p}}_L \in i_{12}^*(\J)$, the left--velocity, is
known. Taking into account the action of $V(\M)$ as modelling
vector bundle of the affine bundle $\J$, this is equivalent to the
assignment of an impulse ${\bf I }={\bf I}({\bf{p}}_L) \in
i^*_{12}(V(\M))$ such that ${\bf{p}}_R = {\bf{p}}_L +
{\bf{I}}({\bf{p}}_L)$. A constitutive law is then an assignment
\begin{eqnarray}
\begin{array}{rcccl}
{\bf I}_{const} & : & i^*_{12}(\J) & \to & i^*_{12}(V(\M)) \\
&&{\bf{p}}_L & \rightsquigarrow & {\bf I}_{const}({\bf{p}}_L) \, .
\end{array}
\end{eqnarray}

The ideal constitutive law for multiple impact presented in
(\cite{Pasquero2016Multiple}) is based on three assumptions:
\begin{itemize}
\item[a)] the preservation of kinetic energy of the system before
and after the contact/impact with the constraints in every frame
of reference for which the requirement has a clear meaning;

\item[b)] in absence of additional information about the
constraints, the constraints involved in a multiple contact/impact
cannot be discerned;

\item[c)] in case of contact of the disk with both the constraints
$\SO, \ST$ but of impact with only one of the constraints, the
constitutive characterization must coincide with the usual ideal
characterization of a single constraint.
\end{itemize}

Let the coefficient $\lambda_{ideal}$ be defined by
\begin{eqnarray*}
\lambda_{ideal} \= -2 \dfrac{\Phi\left( \velort_{12}({\bf p}_L),
\velort_1({\bf p}_L) + \velort_2({\bf p}_L) \right)}{\Phi\left(
\velort_1({\bf p}_L) + \velort_2({\bf p}_L), \velort_1({\bf p}_L)
+ \velort_2({\bf p}_L) \right)} :
\end{eqnarray*}
recalling that the three assumption listed above are insufficient
to determine univocally the constitutive characterization in case
of multiple contacts/impacts, the simplest non trivial ideal
constitutive characterization of possible multiple impacts
presented in \cite{Pasquero2016Multiple} applied to the disk in
the corner is:
\begin{eqnarray}\label{CC_TeoricaIdeale}
\hskip-2truecm
\begin{array}{lclll}
\mbox{no impact} & \Rightarrow  & {\bf I}({\bf{p}}_L) \= 0 & \Leftrightarrow & {\bf p} = {\bf p}_L = {\bf p}_R \\ \\

\mbox{impact with } \SO &
\Rightarrow & {\bf I}({\bf{p}}_L) \= -2 \, \velort_1({\bf p}_L) & \Leftrightarrow & {\bf p} = {\bf p}_L -2 \, \velort_1({\bf p}_L)\\ \\

\mbox{impact with } \ST &
\Rightarrow & {\bf I}({\bf{p}}_L) \= -2 \, \velort_2({\bf p}_L) & \Leftrightarrow & {\bf p} = {\bf p}_L -2 \, \velort_2({\bf p}_L) \\ \\

\mbox{multiple impact} & \Rightarrow & {\bf I}({\bf{p}}_L) \=
\lambda \, \left( \velort_1({\bf p}_L) + \velort_2({\bf
p}_L)\right) & \Leftrightarrow & {\bf p} = {\bf p}_L +
\lambda_{ideal} \, \left( \velort_1({\bf p}_L) + \velort_2({\bf
p}_L)\right)
\end{array}
\end{eqnarray}

\medskip

In complete analogy, let $\veps \in [0,1)$ and let the coefficient
$\lambda_{non-ideal}$ be defined by
\begin{eqnarray*}
\lambda_{non-ideal} \= -(1+\veps) \dfrac{\Phi\left(
\velort_{12}({\bf p}_L), \velort_1({\bf p}_L) + \velort_2({\bf
p}_L) \right)}{\Phi\left( \velort_1({\bf p}_L) + \velort_2({\bf
p}_L), \velort_1({\bf p}_L) + \velort_2({\bf p}_L) \right)} :
\end{eqnarray*}
the simplest non trivial non--ideal constitutive characterization
of possible multiple impacts applied to the disk in the corner is:
\begin{eqnarray}\label{CC_TeoricaNonIdeale}
\hskip-2truecm
\begin{array}{lclll}
\mbox{no impact} & \Rightarrow  & {\bf I}({\bf{p}}_L) \= 0 & \Leftrightarrow & {\bf p} = {\bf p}_L = {\bf p}_R \\ \\

\mbox{impact with } \SO &
\Rightarrow & {\bf I}({\bf{p}}_L) \= -2 \, \velort_1({\bf p}_L) & \Leftrightarrow & {\bf p} = {\bf p}_L -(1+\veps) \, \velort_1({\bf p}_L)\\ \\

\mbox{impact with } \ST &
\Rightarrow & {\bf I}({\bf{p}}_L) \= -2 \, \velort_2({\bf p}_L) & \Leftrightarrow & {\bf p} = {\bf p}_L -(1+\veps) \, \velort_2({\bf p}_L) \\ \\

\mbox{multiple impact} & \Rightarrow & {\bf I}({\bf{p}}_L) \=
\lambda \, \left( \velort_1({\bf p}_L) + \velort_2({\bf
p}_L)\right) & \Leftrightarrow & {\bf p} = {\bf p}_L +
\lambda_{non-ideal} \, \left( \velort_1({\bf p}_L) +
\velort_2({\bf p}_L)\right)
\end{array}
\end{eqnarray}

\medskip

It is however clear that the assignment of a new velocity ${\bf
p}_{new} \= {\bf p}_L + {\bf I}({\bf{p}}_L)$ with the rule
described above does not ensure that the system subject to ${\bf
p}_{new}$ does not impact again with the constraints. It is then
necessary to construct an iterative procedure that applies the
rule until the velocity obtained does not give an impact with one
or both the constraints and can be considered the right velocity
${\bf p}_R$ of the system after the impact. This of course opens
the problem of the termination analysis of the algorithm discussed
above.

\subsection{Coordinate expressions of the iterative rule}

For the convenience of the Reader, we list now the local
coordinate expressions (using the cartesian representation of
$\SO, \ST$ in the fibred coordinates $(t,x,y,\vth)$ for $\M$) of
the main objects described in this appendix and used to obtain the
iterative rule described in the paper. Taking into account that we
focus our attention to multiple contact/impact, let ${\bf p} \in
i_{12}^*(\J)$ be an absolute velocity of the system when the
system is in contact with both $\SO, \ST$ (so that $\pi({\bf p})
\in \SOT$). Then:
\begin{eqnarray}\label{velocita}
\hskip -3.4truecm {\bf p} \= \dfrac{\partial}{\partial t} +
\dot{x} \dfrac{\partial}{\partial x} + \dot{y}
\dfrac{\partial}{\partial y}+ \dot{\vth} \dfrac{\partial}{\partial
\vth};
\end{eqnarray}
\begin{eqnarray}\label{velort}
\begin{array}{l}
\velort_1({\bf p}) \= \dfrac{k}{1+k^2} \, (k\dot{x} - \dot{y})
\dfrac{\partial}{\partial x} - \dfrac{1}{1+k^2} \, (k\dot{x} -
\dot{y}) \dfrac{\partial}{\partial y}; \\ \\

\velort_2({\bf p}) \= \dfrac{k}{1+k^2} \, (k\dot{x} + \dot{y})
\dfrac{\partial}{\partial x} + \dfrac{1}{1+k^2} \, (k\dot{x} +
\dot{y}) \dfrac{\partial}{\partial y}; \\ \\

\velort_{12}({\bf p}) \= \dot{x} \, \dfrac{\partial}{\partial x} +
\dot{y} \, \dfrac{\partial}{\partial y};
\end{array}
\end{eqnarray}
\begin{eqnarray}\label{uelort}
\hskip-3.4truecm
\begin{array}{l} \uelortO \= \dfrac{k}{1+k^2} \,
\dfrac{\partial}{\partial x} -
\dfrac{1}{1+k^2} \, \dfrac{\partial}{\partial y}; \\ \\
\uelortT \= \dfrac{k}{1+k^2} \, \dfrac{\partial}{\partial x} +
\dfrac{1}{1+k^2} \, \dfrac{\partial}{\partial y};
\end{array}
\end{eqnarray}
\begin{eqnarray}\label{scalprod}
\hskip-3.4truecm
\begin{array}{l}
\Phi(\velort_1({\bf p}),\uelortO) \= \dfrac{m}{1+k^2} \, (k\dot{x} - \dot{y})\\ \\

\Phi(\velort_2({\bf p}),\uelortT) \= \dfrac{m}{1+k^2} \, (k\dot{x}
+ \dot{y})
\end{array}
\end{eqnarray}
The coordinates expressions of the various geometrical objects
describing the problem and the rules assigning the constitutive
characterization of the contact/impact suggest two simple remarks:
\begin{itemize}
\item[1)] since the coordinate expressions (\ref{velort}) do not
involve terms pertaining $\frac{\partial}{\partial \vth}$, and the
rules (\ref{CC_TeoricaIdeale},\ref{CC_TeoricaNonIdeale}) assigning
the absolute velocity after the impact involves only the
orthogonal velocities (\ref{velort}), then the component along
$\frac{\partial}{\partial \vth}$ of ${\bf p}$ is not changed by
the impact. This is coherent with the absence of friction of the
disk with $\SO$ and $\ST$ implied by the ideality of the contact.
Moreover, this justifies the fact that the iterative rule can then
be expressed using the coordinates $\dot{x}, \dot{y}$ alone;

\item[2)] taking into account Fig. 1, the conditions for ${\bf p}$
to be an impact velocity for $\SO, \ST$ are then
\begin{eqnarray*}
\Phi(\velort_1({\bf p}),\uelortO)>0 \, , \quad \Phi(\velort_2({\bf
p}),\uelortT)>0
\end{eqnarray*}
respectively. Together with (\ref{scalprod}), this justifies the
description of the zones $\Z_0, \Z_1, \Z_2, \Z_{12}$ in terms of
the signs of $k\dot{x} + \dot{y}, k\dot{x} - \dot{y}$.
\end{itemize}

\end{document}